\documentclass[twoside,11pt,paperA4]{article}
\usepackage{amsmath,amssymb,graphicx,amsthm}
\usepackage{mathtools}
\usepackage[utf8]{inputenc}
\usepackage{anysize}
\usepackage{enumerate}
\marginsize{2cm}{2cm}{2cm}{3cm}

\DeclareMathOperator\supp{supp}

\usepackage{hyperref} 
\hypersetup{colorlinks=true, citecolor=blue, urlcolor  = blue} 

\usepackage{tikz}
\usepackage{colortbl}
\usepackage{stmaryrd} 

\theoremstyle{definition}
  
  \newtheorem{remark}{Remark}

\theoremstyle{plain}
  \newtheorem{theorem}{Theorem}[section]

\usepackage{dsfont}
\usepackage{multirow}

\usepackage{caption}
\usepackage{subcaption}

\usepackage{lineno} 
\title{Wave interactions and stability of Riemann solutions for a nonautonomous Chromatography-type system of Langmuir isotherm}
\author{
Richard De la cruz\thanks{Universidad Pedagógica y Tecnológica de Colombia, School of Mathematics and Statistics, Tunja 150003, Colombia. 
\href{mailto:richard.delacruz@uptc.edu.co}{richard.delacruz@uptc.edu.co} } 
,\, 
Rakib Mondal\thanks{Indian Institute of Technology Bombay, Department of Mathematics, Powai, Mumbai, Maharashtra 400076, India.
\href{mailto:rkbmon90@gmail.com}{rkbmon90@gmail.com} }
,\, and
Wladimir Neves\thanks{Universidade Federal do Rio de Janeiro, Instituto de Matemática, Cidade Universitária, 21945-970, Rio de Janeiro, Rio de Janeiro, Brazil.
\href{mailto:wladimir@im.ufrj.br}{wladimir@im.ufrj.br} }
}
\date{
}

\begin{document}
\maketitle

\begin{abstract}
We investigate the wave interactions and stability of Riemann solutions for a nonautonomous chromatography-type system of Langmuir isotherm with time-dependent damping and flux. The system models two-component chromatographic separation with a time-dependent saturation capacity $n(t)$, leading to a nonautonomous hyperbolic system of balance laws. We consider a perturbed Riemann problem with piecewise constant initial data having two jump discontinuities at $x = \pm\epsilon$, and construct the global weak solution by analyzing all possible wave interactions, both classical (shock waves, rarefaction waves, contact discontinuities) and nonclassical (delta shock waves). We prove that as $\epsilon \to 0$, the solution of the perturbed Riemann problem converges to the solution of the corresponding Riemann problem in the space of Radon measures, establishing the stability of Riemann solutions under small perturbations of the initial data. To the best of our knowledge, this is the first instance of wave interaction and stability analysis for a nonautonomous chromatography-type system with time-dependent coefficients. Numerical experiments using a Lax-Friedrichs type scheme illustrate the wave interaction structure, the profiles at selected times, and the asymptotic convergence as $\epsilon \to 0$.
\end{abstract}

\textbf{Keywords:} Chromatography-type system; Nonlinear wave interactions; Nonclassical waves; Stability analysis; Nonautonomous system

\section{Introduction}

The hyperbolic system of balance laws takes center stage in the field of partial differential equations arising in many physical problems in science and engineering, including industrial problems and chemical processes. One of the most important such examples is the chromatography process, widely used in science, engineering, and industry, to separate a chemical mixture of two or more solutes into its components. In a chromatographic column, two or more components pass through a reactor filled with solid particles, which adsorb different amounts of components, leading to the distribution of different concentrations that move down the reactor at different rates. The most extensively used adsorption isotherm is the Langmuir adsorption isotherm, first introduced by Langmuir \cite{Langmuir_1916} for the ideal localized single layer system. Since then, it has been widely used by researchers for both single solute and multiple solute systems; see, e.g., \cite{glueckauf_1946, glueckauf_1949, Marco_Mazzotti_2009, Marco_Mazzotti_2010} and \cite{De_boer_1953, Rhee_1970} for a more detailed discussion.\\

In chromatographic modeling, the assumption of constant adsorption capacity is a common simplification that may not fully account for variations in operating conditions. In industrial and laboratory settings, the effective adsorption capacity of the stationary phase may vary over time due to several physical mechanisms, including temperature variations, gradual activation of adsorption sites, and progressive changes in the accessibility of the porous structure. To model these effects in chromatography column, the authors in \cite{Rcruz_RM_WN} introduced a time-dependent saturation capacity $n(t)$ and developed a new nonautonomous system of balance laws in which the source term $\sigma(t) = n'(t)/n(t)$ arises intrinsically as the logarithmic rate of change of the adsorption capacity. In \cite{Rcruz_RM_WN}, they proposed and studied the following nonautonomous chromatography system of Langmuir isotherm:
\begin{equation} \label{sys}
\begin{cases}
v_t + \left(\dfrac{v}{h(t)+v}\right)_x = -\sigma(t)v,\\[6pt]
w_t + \left(\dfrac{w}{h(t)+v}\right)_x = -\sigma(t)w,
\end{cases}
\quad (x,t) \in \mathbb{R} \times \mathbb{R}_+,
\end{equation}
where $h(t) = e^{-\int_0^t \sigma(s)ds}$, $v(x,t) \geq 0$ and $w(x,t) \in \mathbb{R}$. In \cite{Rcruz_RM_WN}, the authors constructed explicit Riemann solutions for \eqref{sys}, distinguishing cases depending on the initial data, and established the existence of delta-shock wave solutions as the vanishing viscosity limit in the singular case $v_- = 0$.\\

The stability of Riemann solutions under perturbations of the initial data is a fundamental question in the theory of hyperbolic systems of conservation laws. For classical autonomous systems, this has been studied extensively through wave interaction analysis; see, e.g., \cite{shen_JMAA_2010, M_Sun_AML} for the homogeneous chromatography system, and \cite{GUO_2014, Q_Zhang_ZAMP, Cheng_JMAA_2011} for related systems. The key idea is to consider a perturbed Riemann problem with two jump discontinuities and to analyze the interactions of the elementary waves emanating from each discontinuity, showing that the global solution converges to the Riemann solution as the perturbation parameter $\epsilon \to 0$. For nonautonomous systems with time-dependent source terms, however, the wave interaction analysis is significantly more complex, since the waves are no longer self-similar and their speeds depend explicitly on time through 
$h(t)$; see \cite{DELACRUZ_JDE, Rcruz_RM_WN} for related nonautonomous systems. However, these references does not explore the wave interaction problems. In fact, there is no result available for wave interactions including delta shock wave for these nonautonomous systems. Hence, this motivates us to investigate the wave interactions to establish the stability of Riemann solutions.\\

In this paper, we investigate the wave interactions and stability of Riemann solutions for the nonautonomous chromatography system \eqref{sys}. Specifically, we consider the perturbed Riemann problem with piecewise constant initial data having two jump discontinuities at $x = \pm\epsilon$:
\begin{equation} \label{PRP}
(v(x,0), w(x,0)) = \begin{cases}
(v_-, w_-), & x < -\epsilon,\\
(v_\sim, w_\sim), & -\epsilon < x < \epsilon,\\
(v_+, w_+), & x > \epsilon,
\end{cases}
\end{equation}
where $(v_\pm, w_\pm)$, $(v_\sim, w_\sim)$ are arbitrary constant states and $\epsilon > 0$ is a small perturbation parameter. We construct the global weak solution of \eqref{sys}--\eqref{PRP} by analyzing all possible wave interactions, both classical (shock waves, rarefaction waves, contact discontinuities) and nonclassical (delta shock waves), and prove that as $\epsilon \to 0$ the solution converges to the Riemann solution of \eqref{sys} in the space of Radon measures \cite{Ambrosio_SIAM}, establishing the stability of Riemann solutions under small perturbations of the initial data.\\

The main contributions of this paper are the following. We provide a complete classification of all wave interaction cases for the nonautonomous system \eqref{sys}, including both classical and nonclassical interactions involving delta shock waves. Moreover, we prove that the Riemann solutions of \eqref{sys} are stable under small perturbations of the initial data, in the sense of Theorem~\ref{main_thm}. Finally, we provide numerical evidence of the wave interaction structure and the asymptotic stability, using a Lax-Friedrichs type scheme \cite{Rcruz_RM_WN} and varying the perturbation parameter $\epsilon$.\\

To the best of our knowledge, this is the first instance of wave interaction and stability analysis for a nonautonomous chromatography-type system with time-dependent damping and flux, extending the results of \cite{shen_JMAA_2010, M_Sun_AML} to the nonautonomous setting. This contribution bridges the gap between physically relevant modeling of chromatographic processes with time-dependent adsorption capacity \cite{Rcruz_RM_WN} and the modern nonlinear wave theory for hyperbolic systems with source terms \cite{Rcruz_RM_WN, DELACRUZ_JDE}.\\

The remainder of the paper is organized as follows. In Section~\ref{sec:prelim}, we recall the elementary waves and Riemann solutions for \eqref{sys} from \cite{Rcruz_RM_WN}. In Section~\ref{sec:Wint}, we analyze all wave interactions for the perturbed Riemann problem \eqref{PRP} and prove the main stability result, Theorem~\ref{main_thm}. In Section~\ref{sec:numerics}, we present numerical experiments illustrating the wave interaction structure, the profiles at selected times, and the asymptotic convergence as $\epsilon \to 0$. Finally, we conclude the results and discuss some future directions in Section \ref{sec:con}.

\section{Preliminaries} \label{sec:prelim}

\noindent By virtue of the following variable transformation \cite{Rcruz_RM_WN} 
\begin{equation*} 
    (V(x,t),W(x,t))=\Big(\frac{v(x,t)}{h(t)},\frac{w(x,t)}{h(t)}\Big),
\end{equation*}
the system \eqref{sys} is transformed into a conservative form
\begin{equation} \label{AuxSys}
    \begin{cases}
        V_t + \frac{1}{h(t)} (\frac{V}{1+V})_x = 0,\\
        W_t + \frac{1}{h(t)} (\frac{W}{1+V})_x = 0.
    \end{cases}
\end{equation}
We reformulate the system \eqref{AuxSys} for smooth $(V,W)$ into the quasilinear form
 \begin{equation*} 
     \begin{pmatrix}
         V\\ W
     \end{pmatrix}_t+ 
     \begin{pmatrix}
         \frac{1}{h(1+V)^2} & 0\\
         -\frac{W}{h(1+V)^2} & \frac{1}{h(1+V)}
     \end{pmatrix}
     \begin{pmatrix}
         V\\ W
     \end{pmatrix}_x = \begin{pmatrix}
         0\\0
     \end{pmatrix}
 \end{equation*}
to obtain the eigenvalues and corresponding right eigenvectors of \eqref{AuxSys} as follows
\begin{equation}\nonumber
    \lambda_1=\frac{1}{(1+V)^2 h(t)}, \quad \lambda_2=\frac{1}{(1+V)h(t)},
\end{equation}
and respectively
\begin{equation}\nonumber
    \mathbf{r}_1=\begin{pmatrix}
        V\\W
    \end{pmatrix}, \quad
    \mathbf{r}_2=\begin{pmatrix}
        0\\1
    \end{pmatrix}.
\end{equation}
Observe that $\nabla\lambda_1\cdot \mathbf{r}_1\neq 0$ and $\nabla\lambda_2\cdot \mathbf{r}_2= 0$ for all $(V,W)$ with $V> 0$, therefore system \eqref{AuxSys} is strictly hyperbolic with $\lambda_1$ is genuinely nonlinear and $\lambda_2$ is linearly degenerate, for $V> 0$. Hence, the associated elementary waves are rarefaction or shock waves corresponding to $\lambda_1$ characteristic field and contact discontinuity corresponding to $\lambda_2$. 

For a given left state $(v_-, w_-)$, the set of states $(V, W)$ that can be connected on the right by a 1-rarefaction wave is as follows
\begin{equation*} 
    R_1(v_-,w_-): 
    \begin{cases}
        \xi = \frac{1}{(1+V)^2},\\
        \frac{W}{V}=\frac{w_-}{v_-}, \quad V \le v_-,
    \end{cases}
\end{equation*}
where the similarity variable $\xi$ is given by (see \cite{Rcruz_RM_WN})
\begin{equation*}
    \xi = \frac{x}{\int_0^t \frac{1}{h(s)}ds}.
\end{equation*}

For a given fixed $(v_-, w_-)$, the set of states $(V, W)$ that are connected to $(v_-, w_-)$ on the right by a 1-shock wave, is given by
\begin{equation*}\label{1-shock}
    S_1(v_-,w_-):
    \begin{cases}
        \frac{dx}{dt} = \frac{1}{(1+V)(1+v_-)h(t)},\\
        \frac{W}{V} = \frac{w_-}{v_-}, \quad V>v_-.
    \end{cases}
\end{equation*}

On the other hand, the contact discontinuity corresponding to $\lambda_2$, namely $J_2(v_-,w_-)$ starting from $(v_-,w_-)$ can be expressed as follows 
\begin{equation*}\label{ContDisc}
    J_2(v_-,w_-):
    \begin{cases}
        \frac{dx}{dt}=\frac{1}{(1+v_-)h(t)}= \frac{1}{(1+V)h(t)},\\
         V=v_-.
    \end{cases}
\end{equation*}

According to the reference \cite{Rcruz_RM_WN}, the Riemann problem for the system \eqref{sys} with the initial data
\begin{equation} \label{RP}
    (v(x,0), w(x,0))= \begin{cases}
        (v_-, w_-),\quad \text{if}\; x<0,\\
        (v_+, w_+),\quad \text{if}\; x>0,
    \end{cases}
\end{equation}
has a solution consisting of 1-rarefaction wave followed by 2-contact discontinuity, i.e., $R_1+J_2$, which is given by
\begin{equation*}
    (v(x,t),w(x,t)) =
    \begin{cases}
        h(t)(v_-,w_-),  &\mbox{if } x < g(t; v_-),\\
        h(t)(\sqrt{\frac{\int_0^t \frac{1}{h(s)}ds}{x}}-1, \frac{w_-}{v_-}(\sqrt{\frac{\int_0^t \frac{1}{h(s)}ds}{x}}-1)), &\mbox{if } g(t; v_-) \le x \le g(t; v_+),\\
        h(t)(v_+,\frac{v_+ w_-}{v_-}), &\mbox{if } g(t; v_+) < x < g_2(t;v_+),\\
        h(t)(v_+,w_+), &\mbox{if } x > g_2(t;v_+),
    \end{cases}
\end{equation*}
where $g(t;v) = \frac{1}{(1+v)^2}\int_0^t \frac{1}{h(s)}ds$ and $g_2(t;v) = \frac{1}{1+v}\int_0^t \frac{1}{h(s)}ds$,
if $v_+<v_-$. If $0<v_-<v_+$, the solution to the Riemann problem \eqref{sys} and \eqref{RP} consists of a 1-shock wave followed by 2-contact discontinuity, i.e., $S_1+J_2$, given by
\begin{equation}\label{S+J}
    (v(x,t),w(x,t)) =
    \begin{cases}
        h(t)(v_-,w_-),  &\mbox{if } x < g_1(t; v_-,v_+),\\
        h(t)(v_+,\frac{v_+ w_-}{v_-}), &\mbox{if } g_1(t; v_-, v_+) < x < g_2(t;v_+),\\
        h(t)(v_+,w_+), &\mbox{if } x > g_2(t;v_+),
    \end{cases}
\end{equation}
    where $g_1(t;v_1,v_2)=\frac{1}{(1+v_1)(1+v_2)}\int_0^t \frac{1}{h(s)}ds$.

Finally, if $0=v_-<v_+$, the solution to the Riemann problem \eqref{sys} and \eqref{RP} takes the form of a delta-shock wave as follows:
    \begin{equation*} \label{deltaS}
    (v(x,t),w(x,t)) =
    \begin{cases}
        (0, h(t) w_-), &\mbox{if } x < \theta \int_0^t \frac{1}{h(s)} ds,\\
        (h(t)v_+, \alpha(t)\delta(x-\theta \int_0^t \frac{1}{h(s)} ds)), &\mbox{if } x = \theta \int_0^t \frac{1}{h(s)} ds,\\
        (h(t) v_+, h(t) w_+), &\mbox{if } x > \theta \int_0^t \frac{1}{h(s)} ds,
    \end{cases}
\end{equation*}
in the sense of distribution, in which $x(t)=\theta \int_0^t \frac{1}{h(s)} ds$ is the delta shock curve, $v_\delta(t)=h(t)v_+$ is the value of $v$ on $x=x(t)$ and $\alpha(t)=h(t)\frac{w_-v_+}{1+v_+}\int_0^t \frac{1}{h(s)} ds$, $\theta=\frac{1}{1+\frac{1}{h(t)}v_\delta}$.

\medskip

\noindent \textbf{Notations and spaces}
\medskip

\noindent To prepare the detailed analysis in the following section, let us first list some essential notations and define the spaces relevant to our study. We denote $\mathbb{R}_+^2 = \mathbb{R} \times (0,+\infty)$ and $\overline{\mathbb{R}_+^2} = \mathbb{R} \times [0,+\infty)$. A brief review of the left- and right-hand delta functions, which we will use extensively, is essential; for more details, see \cite{M_Nedeljkov_08}. Let $\mathbb{R}_+^2$ be divided into finitely many disjoint nonempty open sets $\Omega_i$, $i=1,\dots,n$ in $\mathbb{R}_+^2$, with piecewise smooth boundary curves $\Gamma_i$, $i=1,\dots,n$, which satisfies $\Omega_i \cap \Omega_j = \varnothing$ if $i \neq j$, and $\bigcup\limits_{i=1}^n \overline{\Omega_i}=\overline{\mathbb{R}_+^2}$, where $\overline{\Omega_i}$ denotes the closure of $\Omega_i$. Suppose, $\mathcal{C}(\overline{\Omega_i})$ be the space of continuous bounded functions from $\overline{\Omega_i}$ into $\mathbb{R}$, equipped with the uniform $L^\infty$ norm, and $\mathcal{M}(\overline{\Omega_i})$ be the space of measures on $\overline{\Omega_i}$. Let us now consider the Cartesian product spaces as follows 
\begin{equation*}
    \mathcal{C}_\Gamma:=\prod\limits_{i=1}^n \mathcal{C}(\overline{\Omega_i}), \quad \mathcal{M}_\Gamma:=\prod\limits_{i=1}^n \mathcal{M}(\overline{\Omega_i}).
\end{equation*}
Let $G=(G_1,\dots, G_n)\in \mathcal{C}_\Gamma$ and $D= (D_1,\dots, D_n)\in \mathcal{M}_\Gamma$. Then we define the product of an element of $\mathcal{C}_\Gamma$ with an element of $\mathcal{M}_\Gamma$ using the usual product of a continuous function and a measure as $G\cdot D := (G_1 D_1, \dots, G_n D_n)$, where $G_i D_i$ ($i=1, \dots, n$) denotes usual product of a continuous function $G_i$ and a measure $D_i$. Furthermore, any $D_i\in \mathcal{M}(\overline{\Omega_i})$ (a measure on $\overline{\Omega_i}$) can be treated as a measure on the whole space $\overline{\mathbb{R}_+^2}$ together with a support in $\overline{\Omega_i}$. This implies the mapping $m: \mathcal{M}_\Gamma \rightarrow \mathcal{M}(\overline{\mathbb{R}_+^2})$ defined by 
\begin{equation*}
    m(D):=\sum\limits_{i=1}^n D_i, \quad D_i\in \mathcal{M}(\overline{\Omega_i}), \; i=1,\dots,n
\end{equation*}
is well defined. For more clarity, we consider a typical example when $\mathbb{R}_+^2$ is divided by a piecewise smooth curve $x=\Upsilon(t)$ into two open sets $\Omega_1$ and $\Omega_2$. Then $\Omega_1 \cap \Omega_2=\varnothing$, $\overline{\mathbb{R}_+^2}= \overline{\Omega_1}\cup \overline{\Omega_2}$, and the delta measure $\delta(x-\Upsilon(t))$ supported on the curve $x=\Upsilon(t)$ can splits into two components as $D^-\in \mathcal{M}(\overline{\Omega_1})$ and $D^+\in \mathcal{M}(\overline{\Omega_2})$ in a non-unique way as follows
\begin{equation*}
    \delta(x-\Upsilon(t))= \alpha^-(t) D^- + \alpha^+(t) D^+= m(\alpha^-(t) D^-, \alpha^+(t) D^+), \quad\text{and}\; \alpha^-(t)+ \alpha^+(t)=1.
\end{equation*}


\section{Wave interactions} \label{sec:Wint}
In this section, we investigate the global weak entropy solution by means of the local Riemann solutions for the perturbed Riemann problem with two jump discontinuities of system \eqref{sys} as \eqref{PRP} 
where $(v_\pm, w_\pm), (v_\thicksim, w_\thicksim)\in [0,\infty)\times \mathbb{R}$ are arbitrary constant states, and $\epsilon>0$ is a small perturbation parameter. To construct the global solution of the initial value problem \eqref{sys} with \eqref{PRP}, we encounter the interactions of Riemann solutions in finite time, both classical and non-classical waves emanating from the points $-\epsilon$ and $\epsilon$. In what follows, we discuss several cases of wave interactions for all possible choices of the step-like initial data \eqref{PRP} and consequently, establish the global solutions for the Cauchy problem. 
More precisely, we shall prove the main result of this article.
\begin{theorem} \label{main_thm}
    The Cauchy problem for \eqref{sys} with the arbitrary step-like perturbed initial data \eqref{PRP} admits a global weak solution in the sense of distribution. Moreover, as the perturbation parameter $\epsilon\to 0$, the solution of \eqref{sys} and \eqref{PRP} converges to the solution of the corresponding Riemann problem \eqref{sys} and \eqref{RP} in the space of Radon measures $\mathcal{M}_{loc}(\overline{\mathbb{R}_+^2})$, which suggests that the Riemann solutions for the nonautonomous Chromatography-type system \eqref{sys} are stable under small perturbations of initial data.
\end{theorem}

\begin{remark}
    In the case of the interaction between classical elementary waves, the limit in Theorem \ref{main_thm} becomes $L^1$ limit. In the subsequent sections, below, we establish the proof of the Theorem \eqref{main_thm}, which is the main finding of the article. To the best of our knowledge, this is the first instance of wave interactions and stability results for a nonautonomous system of balance laws.
\end{remark}

\begin{figure}[ht!]  
    \centering
\includegraphics[width=0.7\linewidth]{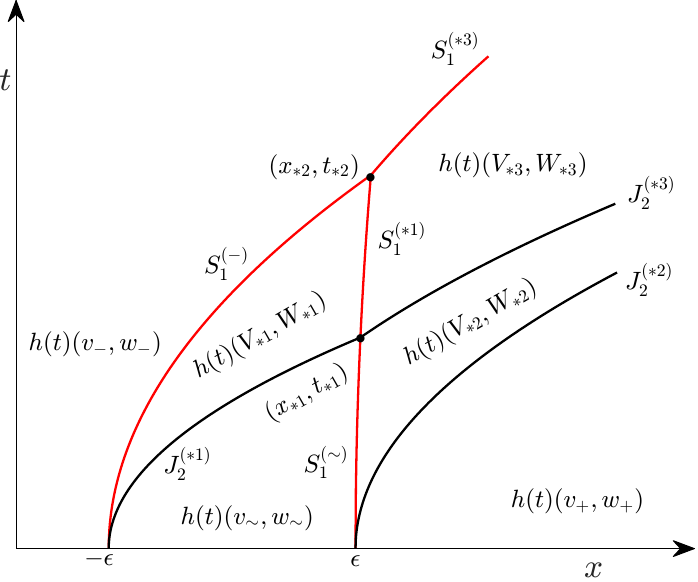}
\caption{Wave interactions and solution to the perturbed Riemann problem \eqref{sys} and \eqref{PRP} when $0<v_-<v_\thicksim<v_+$.}
\label{Fig_case1}
\end{figure}

\subsection{\textbf{Case 1:} $0<v_-<v_\thicksim<v_+$.} \label{Sect3.1}

From \eqref{S+J}, for sufficiently small time $t$, the solution of \eqref{sys} and \eqref{PRP} can be expressed as follows
\begin{equation}\label{case1}
\begin{aligned}
    (v_- h(t), w_- h(t))+ S_1^{(-)}+(V_{\ast 1}h(t), W_{\ast 1}h(t))+J_2^{(\ast 1)}+ (v_\thicksim h(t), w_\thicksim h(t))+S_1^{(\thicksim)}\\ + (V_{\ast 2}h(t), W_{\ast 2}h(t)) + J_2^{(\ast 2)}+ (v_+ h(t), w_+ h(t)),
    \end{aligned}
\end{equation}
where $S_1^{(-)}$ and $J_2^{(\ast 1)}$ are respectively the 1-shock wave and the 2-contact discontinuity starting from the point $(-\epsilon,0)$, connecting the left state $(v_-h(t),w_-h(t))$ to the middle state $(v_\thicksim h(t), w_\thicksim h(t))$ from left to right, and separating by the intermediate state $(V_{\ast1}h(t), W_{\ast1}h(t))$ (see Figure \ref{Fig_case1}). Moreover, $S_1^{(\thicksim)}$ and $J_2^{(\ast 2)}$ are the 1-shock wave and 2-contact discontinuity starting from $(\epsilon,0)$, which connect the states $(v_\thicksim h(t), w_\thicksim h(t))$ and $(v_+ h(t), w_+ h(t))$ from left to right with the intermediate state $(V_{\ast2}h(t), W_{\ast2}h(t))$, respectively. Here, the intermediate states $(V_{\ast 1}h(t), W_{\ast 1}h(t))$ and $(V_{\ast 2}h(t), W_{\ast 2}h(t))$ are determined as 
\begin{equation} \label{interm_V1,W1}
    (V_{\ast 1}h(t), W_{\ast 1}h(t))=\left(v_\thicksim h(t), \frac{v_\thicksim w_-}{v_-}h(t)\right),\quad (V_{\ast 2}h(t), W_{\ast 2}h(t))= \left(v_+ h(t), \frac{v_+ w_\thicksim}{v_\thicksim} h(t)\right).
\end{equation}
However, note that the solution \eqref{case1} is valid only for sufficiently small values of $t$ as the discontinuity arises in a finite time due to the interaction of these waves. So, to construct the global solution for any time $t$, one has to deal with these interaction problems. In fact, the speeds of propagation of $J_2^{(\ast 1)}$ and $S_1^{(\thicksim)}$ are respectively given by $\sigma_2^{(\ast 1)}=\frac{1}{(1+v_\thicksim)h}$, and $s_1^{(\thicksim)}= \frac{1}{(1+v_\thicksim)(1+V_{\ast 2})h}=\frac{1}{(1+v_\thicksim)(1+v_+)h}$. Then, one obtains 
\begin{equation*}
    \sigma_2^{(\ast 1)}- s_1^{(\thicksim)}= \frac{v_+}{(1+v_\thicksim)(1+v_+)h}>0,
\end{equation*}
which implies $J_2^{(\ast 1)}$ must interact $S_1^{(\thicksim)}$ in a finite time, say $t=t_{\ast1}>0$. We denote $x(t_{\ast1}):= x_{\ast1}$. Then, the first point of interaction $(x_{\ast1}, t_{\ast1})$ can be determined from
\begin{equation}\label{IntPnt_x1t1}
    \begin{aligned}
        &x_{\ast1}+ \epsilon= \frac{1}{1+v_\thicksim}\int_0^{t_{\ast1}} \frac{1}{h(s)}ds,\\
       & x_{\ast1}-\epsilon = \frac{1}{(1+v_\thicksim)(1+v_+)}\int_0^{t_{\ast1}} \frac{1}{h(s)}ds.
    \end{aligned}
\end{equation}
Now, this interaction leads to a new Riemann problem at $(x_{\ast1}, t_{\ast1})$ with $(V_{\ast 1}h(t), W_{\ast 1}h(t))$ and $(V_{\ast 2}h(t), W_{\ast 2}h(t))$ as the left hand and right hand states, respectively. Then, we have $$V_{\ast1}=v_\thicksim < v_+=V_{\ast2},$$
which implies the new Riemann problem at $(x_{\ast1}, t_{\ast1})$ can be connected by a 1-shock wave, namely $S_1^{(\ast 1)}$ and a 2-contact discontinuity, namely $J_2^{(\ast 3)}$ with the intermediate state $(V_{\ast 3}h(t), W_{\ast 3}h(t))$, given by
\begin{equation}\label{int_ast3}
    (V_{\ast 3}, W_{\ast 3})= \left(V_{\ast 2}, \frac{V_{\ast 3} W_{\ast 1}}{V_{\ast 1}}\right)= \left(v_+,\frac{v_+ w_-}{v_-}\right).
\end{equation}
Moreover, the interaction between $J_2^{(\ast 1)}$ and $S_1^{(\thicksim)}$ produce a new 1-shock wave $S_1^{(\ast 1)}$ propagating from $(x_{\ast1}, t_{\ast1})$ with speed 
\begin{equation}\label{s1ast1}
s_1^{(\ast 1)}= \frac{1}{(1+V_{\ast 3})(1+V_{\ast 1})h}= \frac{1}{(1+v_+)(1+v_\thicksim)h},
\end{equation}
and a new 2-contact discontinuity $J_2^{(\ast 3)}$ propagating from $(x_{\ast1}, t_{\ast1})$ with speed 
\begin{equation}\label{sigma2ast}
    \sigma_2^{(\ast 3)}= \frac{1}{(1+V_{\ast 3})h}= \frac{1}{(1+V_{\ast 2})h}= \frac{1}{(1+v_+)h}.
\end{equation}

Now, in order to examine whether $J_2^{(\ast 3)}$ interacts with $J_2^{(\ast 2)}$ or not, let us trace their speeds. In fact, from the expression of 2-contact discontinuity, the speed of $J_2^{(\ast 2)}$, namely $\sigma_2^{(\ast 2)}$, is given by 
\begin{equation*}
    \sigma_2^{(\ast 2)}= \frac{1}{(1+V_{\ast 2})h}= \frac{1}{(1+v_+)h},
\end{equation*}
which coincides with the speed of $J_2^{(\ast 3)}$. Hence, the contact discontinuities $J_2^{(\ast 3)}$ and $J_2^{(\ast 2)}$ are parallel, and will not interact. On the other hand, the interaction between $S_1^{(-)}$ and $J_2^{(\ast 3)}$ is not possible, since the speed of the wave on the left-hand side, i.e., $S_1^{(-)}$, is smaller than the speed of $J_2^{(\ast 3)}$. Thus, the 1-shock wave $S_1^{(-)}$ will never catch the contact discontinuity $J_2^{(\ast 3)}$. But, there is a possibility that $S_1^{(-)}$ might interact with $S_1^{(\ast 1)}$. Let us compare their propagating speeds. The speed of propagation of $S_1^{(-)}$ is denoted as $s_1^{(-)}$, and given by
\begin{equation*}
    s_1^{(-)}= \frac{1}{(1+V_{\ast 1})(1+v_-)h}= \frac{1}{(1+v_\thicksim)(1+v_-)h}.
\end{equation*}
Then, we have
\begin{equation*}
\begin{aligned}
    s_1^{(-)}- s_1^{(\ast 1)}&= \frac{1}{(1+v_\thicksim)(1+v_-)h} - \frac{1}{(1+v_+)(1+v_\thicksim)h}\\
    &= \frac{1}{(1+v_\thicksim)h} \frac{v_+ - v_-}{(1+v_-)(1+v_+)} >0\\
    \implies & s_1^{(-)} > s_1^{(\ast 1)}.
    \end{aligned}
\end{equation*}
This shows that $S_1^{(-)}$ must interact with $S_1^{(\ast 1)}$ in a later time, say $t=t_{\ast 2} >0$. Let us set $x(t_{\ast 2}):= x_{\ast 2}$. Then, the point of interaction $(x_{\ast 2}, t_{\ast 2})$ are given by
\begin{equation*}\label{IntPt_x2t2}
    \begin{aligned}
        &x_{\ast 2}+ \epsilon= \frac{1}{(1+v_\thicksim)(1+v_-)}\int_0^{t_{\ast 2}} \frac{1}{h(s)}ds,\\
       & x_{\ast 2}- x_{\ast 1} = \frac{1}{(1+v_+)(1+v_\thicksim)}\int_{t_{\ast 1}}^{t_{\ast 2}} \frac{1}{h(s)}ds.
    \end{aligned}
\end{equation*}
Moreover, this interaction gives rise to a new Riemann problem at $(x_{\ast 2}, t_{\ast 2})$ with the initial data having the left-hand and right-hand states respectively, as follows $(v_- h(t), w_- h(t))$ and $(V_{\ast 3}h(t), W_{\ast 3}h(t))$, where $(V_{\ast 3}, W_{\ast 3})= \left(v_+,\frac{v_+ w_-}{v_-}\right)$ (follow Figure \ref{Fig_case1}). Since, $V_{\ast 3}= v_+> v_-$ and $W_{\ast 3}= \frac{v_+ w_-}{v_-}$, it follows from the Theorem on classical Riemann solution that states $(v_- h(t), w_- h(t))$ and $(V_{\ast 3}h(t), W_{\ast 3}h(t))$ must be connected by a single shock wave, namely, $S_1^{(\ast 3)}$. That is, the new Riemann problem at $(x_{\ast 2}, t_{\ast 2})$ is solved by a 1-shock wave $S_1^{(\ast 3)}$ propagating from $(x_{\ast 2}, t_{\ast 2})$ with speed, say $s_1^{(\ast 3)}$, given by
\begin{equation*}\label{Speed_S3}
    s_1^{(\ast 3)}= \frac{1}{(1+V_{\ast 3})(1+v_-)h}= \frac{1}{(1+v_+)(1+v_-)h}.
\end{equation*}

Comparing the speeds of $S_1^{(\ast 3)}$ and $J_2^{(\ast 3)}$ yields
\begin{equation*}
    s_1^{(\ast 3)}- \sigma_2^{(\ast 3)}= \frac{-v_-}{(1+v_-)(1+v_+)h}<0,
\end{equation*}
then $s_1^{(\ast 3)} < \sigma_2^{(\ast 3)}$, and therefore no further interaction is possible for time $t> t_{\ast 2}$. Hence, the solution to the Cauchy problem \eqref{sys} with \eqref{PRP} for time $t> t_{\ast 2}$ is of the form
\begin{equation*}\label{Sol_Case1}
\begin{aligned}
    (v_- h(t), w_- h(t))+ S_1^{(\ast 3)}+(V_{\ast 3}h(t), W_{\ast 3}h(t))+J_2^{(\ast 3)}+ (V_{\ast 2}h(t), W_{\ast 2}h(t)) + J_2^{(\ast 2)}+ (v_+ h(t), w_+ h(t)),
    \end{aligned}
\end{equation*}
where $V_{\ast 3}= V_{\ast 2}= v_+$, $W_{\ast 2}= \frac{v_+ w_\thicksim}{v_\thicksim}$, $W_{\ast 3}= \frac{v_+ w_-}{v_-}$, the shock and contact discontinuity curves $S_1^{(\ast 3)}$ and $J_2^{(\ast 3)}$ propagating from $(x_{\ast 2}, t_{\ast 2})$ and $(x_{\ast 1}, t_{\ast 1})$, respectively are given as follows
\begin{equation*}
    \begin{aligned}
        &S_1^{(\ast 3)}: x- x_{\ast 2} = \frac{1}{(1+v_+)(1+v_-)} \int_{t_{\ast 2}}^t \frac{1}{h(s)}ds,\\
        &J_2^{(\ast 3)}: x- x_{\ast 1}= \frac{1}{(1+v_+)} \int_{t_{\ast 1}}^t \frac{1}{h(s)}ds.
    \end{aligned}
\end{equation*}
A schematic representation of the wave interactions for this case is depicted in Figure \ref{Fig_case1}.

Now, letting the perturbed parameter $\epsilon\to 0$, the points of local Riemann problems, i.e., the initial discontinuities $(-\epsilon, 0)$, $(\epsilon, 0)$, and the point of interactions $(x_{\ast1}, t_{\ast1})$, $(x_{\ast 2}, t_{\ast 2})$ converge to the origin $(0, 0)$ in $x$-$t$ plane, only the constant states $(v_-, w_-)$ and $(v_+, w_+)$ remains. Moreover, the initial data \eqref{PRP} tends to the following Riemann initial data \eqref{RP} at $(0,0)$, i.e., 
\begin{equation*}
    (v(x,0), w(x,0))= \begin{cases}
        (v_-, w_-),\quad &\text{if}\; x< 0,\\
        (v_+, w_+),\quad &\text{if}\; x>0,
    \end{cases}
\end{equation*}
as $\epsilon\to 0$. Now the question is whether the solution of the interaction problem \eqref{sys} and \eqref{PRP} converges to the solution of the corresponding Riemann problem \eqref{sys} and \eqref{RP}. In fact, the contact discontinuities $J_2^{(\ast 2)}$ and $J_2^{(\ast 3)}$ propagating respectively from $(\epsilon, 0)$ and $(x_{\ast1}, t_{\ast1})$ must coincides, and converge to the following 2-contact discontinuity propagating from $(0,0)$
\begin{equation*}
    J_2: x= \frac{1}{(1+v_+)} \int_0^t \frac{1}{h(s)}ds.
\end{equation*}
Moreover, as $\epsilon\to 0$, the shock curves $S_1^{(-)}$, $S_1^{(\thicksim)}$, $S_1^{(\ast 1)}$ and $S_1^{(\ast 3)}$ coincide and converge to the following 1-shock curve propagating from $(0,0)$ 
\begin{equation*}
    S_1: x= \frac{1}{(1+v_+)(1+v_-)} \int_{0}^t \frac{1}{h(s)}ds.
\end{equation*}
Hence, the solution to the Cauchy problem \eqref{sys} and \eqref{PRP} converges exactly to the solution of the Riemann problem \eqref{sys} and \eqref{RP} as the perturbation parameter $\epsilon\to 0$. 

On the other hand, as $t\to \infty$, i.e., for sufficiently large time $t> t_{\ast 2}$, the solution to the perturbed Riemann problem \eqref{sys} and \eqref{PRP} can be expressed as  
\begin{equation*}
\begin{aligned}
    (v_- h(t), w_- h(t))+ S_1^{(\ast 3)}+(V_{\ast 3}h(t), W_{\ast 3}h(t))+J_2^{(\ast 3)}+ (v_+ h(t), w_+ h(t)).
    \end{aligned}
\end{equation*}
This is exactly the solution to the Riemann problem \eqref{sys} and \eqref{RP}, which implies that as $t\to \infty$, the solution to the Cauchy problem \eqref{sys} and \eqref{PRP} is governed by the corresponding initial Riemann problem with left and right states $(v_-, w_-)$ and $(v_+,w_+)$, respectively. Hence, in this case, we conclude that the Riemann solution of \eqref{sys} with \eqref{RP} is stable globally with respect to the small perturbation of \eqref{RP}.

\begin{remark}
   Note that in Figures \ref{Fig_case2.1}-\ref{Fig_case7} for rest of the cases, we denote the intermediate states $h(t)(V_{\ast i}, W_{\ast i})$ ($i=1,2,3$) by $(\ast i)$  and the states $h(t)(v_l, w_l)$ ($l=-, \thicksim, +$) by $(l)$ for simplicity.
\end{remark}

\subsection{\textbf{Case 2:} $0<v_\thicksim<v_\pm$.}

Since the Riemann problem at $x=-\epsilon$ with left and right states $(v_-,w_-)$ and $(v_\thicksim, w_\thicksim)$ satisfy $v_\thicksim< v_-$, then solution of Riemann problem is connected by 1-rarefaction wave $R_1^{(-)}$ and 2-contact discontinuity $J_2^{(\ast 1)}$ propagating from $(-\epsilon,0)$. Indeed, the solution is given by
\begin{equation*}
    (v_- h(t), w_- h(t)) + R_1^{(-)} + (V_{\ast 1}h(t), W_{\ast 1} h(t))+ J_2^{(\ast 1)}+ (v_\thicksim h(t), w_\thicksim h(t)), 
\end{equation*}
where the 1-rarefaction wave $R_1^{(-)}$ is given by
\begin{equation} \label{R1(-)}
    R_1^{(-)}: 
    \begin{cases}
        \xi = \frac{1}{(1+V)^2},\\
        \frac{W}{V}=\frac{w_-}{v_-}, \\
        x(0)= -\epsilon, \quad V_{\ast 1} \le V \le v_-, 
    \end{cases}
\end{equation}
and the intermediate state $(V_{\ast 1}h(t), W_{\ast 1} h(t))$ satisfies
\begin{equation*}
    (V_{\ast 1}h(t), W_{\ast 1} h(t))= h(t)\Big(v_\thicksim,\frac{v_\thicksim w_-}{v_-}\Big).
\end{equation*}
The Riemann problem at $x=\epsilon$ having the left state $(v_\thicksim, w_\thicksim)$ and the right state $(v_+, w_+)$ satisfy $0<v_\thicksim<v_+$. Thus, the solution to this Riemann problem consists of a 1-shock wave $S_1^{(\thicksim)}$ and a 2-contact discontinuity $J_2^{(\ast 2)}$, i.e., 
\begin{equation*}
    (v_\thicksim h(t), w_\thicksim h(t))+S_1^{(\thicksim)} + (V_{\ast 2}h(t), W_{\ast 2}h(t)) + J_2^{(\ast 2)}+ (v_+ h(t), w_+ h(t)),
\end{equation*}
the intermediate state $(V_{\ast 2}h(t), W_{\ast 2}h(t))$ is given by 
\begin{equation*}
    (V_{\ast 2}h(t), W_{\ast 2}h(t))= \left(v_+ h(t), \frac{v_+ w_\thicksim}{v_\thicksim} h(t)\right).
\end{equation*}
Thus, the solution of Cauchy problem \eqref{sys} and \eqref{PRP} is given by
\begin{equation*}\label{case2}
\begin{aligned}
    (v_- h(t), w_- h(t)) + R_1^{(-)} + (V_{\ast 1}h(t), W_{\ast 1} h(t))+ J_2^{(\ast 1)}+ (v_\thicksim h(t), w_\thicksim h(t))+ S_1^{(\thicksim)}\\+ (V_{\ast 2}h(t), W_{\ast 2}h(t)) + J_2^{(\ast 2)}+ (v_+ h(t), w_+ h(t)),
    \end{aligned}
\end{equation*}
which exists for a sufficiently small time $t$, since local Riemann solutions may interact with each other in a finite time. In fact, previously in the Case 1, we proved that $J_2^{(\ast 1)}$ interacts with $S_1^{(\thicksim)}$ at $(x_{\ast 1}, t_{\ast 1})$, $t_{\ast 1}>0$, given in \eqref{IntPnt_x1t1}, and the interaction produce new 1-shock wave $S_1^{(\ast 1)}$ and 2-contact discontinuity $J_2^{(\ast 3)}$, both propagating from $(x_{\ast1}, t_{\ast1})$ with speeds given in \eqref{s1ast1} and \eqref{sigma2ast}, respectively, separating by the intermediate state $(V_{\ast 3}, W_{\ast 3})$, given in \eqref{int_ast3}.

\begin{figure}[ht!]  
    \centering
\includegraphics[width=0.75\linewidth]{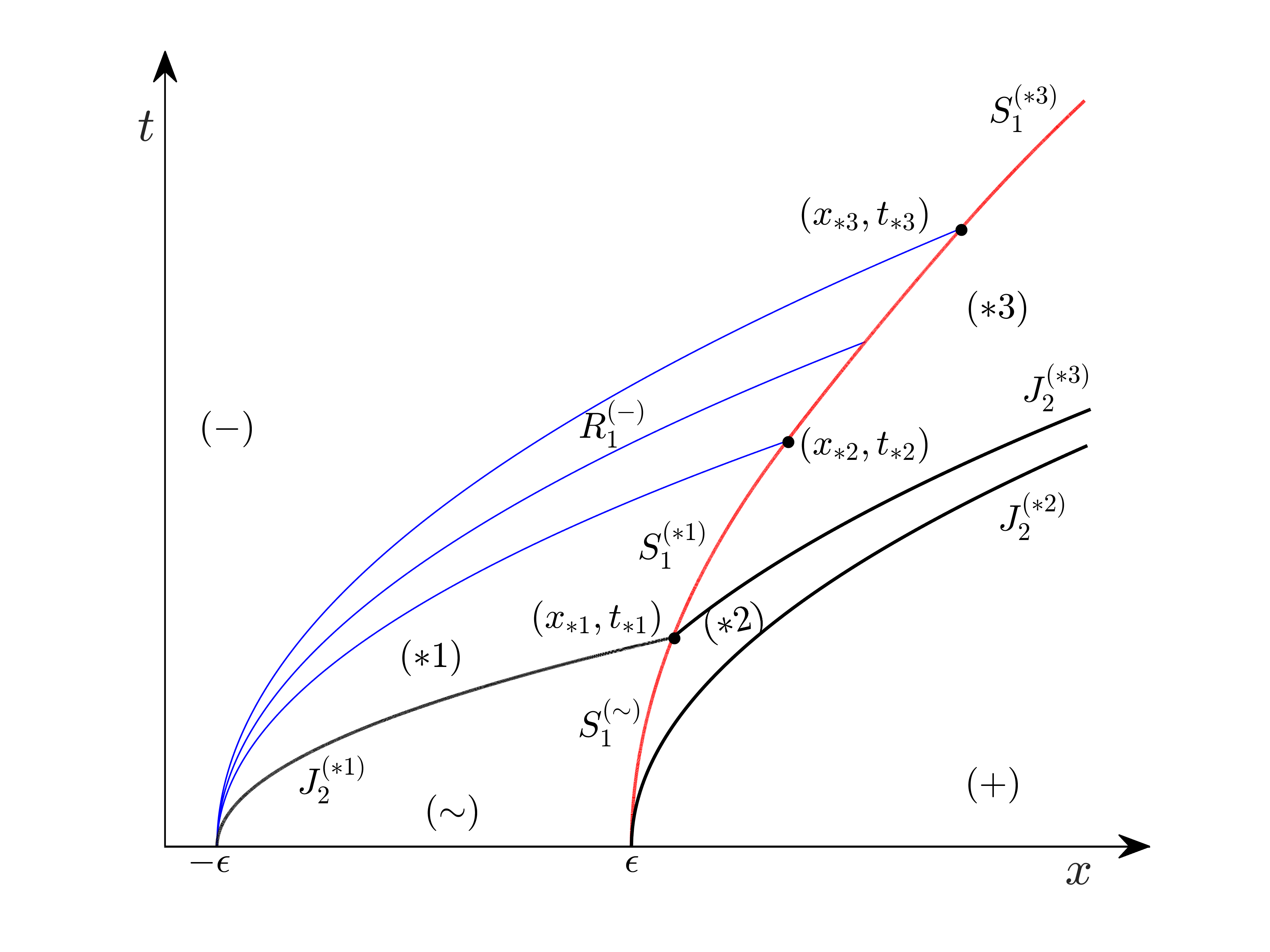}
\caption{Wave interactions and solution to the perturbed Riemann problem for the subcase 2.1, i.e., when $0<v_\thicksim<v_-<v_+$.}
\label{Fig_case2.1}
\end{figure}

From \eqref{s1ast1}, we have the speed of 1-shock $S_1^{(\ast 1)}$ as follows
\begin{equation*}
    s_1^{(\ast 1)}= \frac{1}{(1+V_{\ast 3})(1+V_{\ast 1})h}= \frac{1}{(1+v_+)(1+v_\thicksim)h}.
\end{equation*}
On the other hand, the wave front of $R_1^{(-)}$, i.e., the head of the 1-rarefaction wave $R_1^{(-)}$ is propagating with the speed of 1-characteristic, namely $\theta_1^{(\ast 1)}$, which is given by
\begin{equation*}
   \theta_1^{(\ast 1)} = \frac{1}{(1+V_{\ast 1})^2 h}.
\end{equation*}
Then, comparing the speeds $\theta_1^{(\ast 1)}$ and $s_1^{(\ast 1)}$, we obtain
\begin{equation*}
\begin{aligned}
    \theta_1^{(\ast 1)} -s_1^{(\ast 1)}&= \frac{1}{(1+V_{\ast 1})^2 h} - \frac{1}{(1+V_{\ast 3})(1+V_{\ast 1})h}\\
    &= \frac{1}{(1+v_{\thicksim})^2 h} - \frac{1}{(1+v_+)(1+v_\thicksim)h}\\
    & = \frac{v_+ - v_\thicksim}{(1+v_+)(1+v_\thicksim)^2 h}>0,
    \end{aligned}
\end{equation*}
which implies $\theta_1^{(\ast 1)} >s_1^{(\ast 1)}$. Therefore, the wave front of $R_1^{(-)}$ interact with the 1-shock $S_1^{(\ast 1)}$. Let the point of interaction be $(x_{\ast 2}, t_{\ast 2})$. Then, it satisfies the following set of equations
\begin{equation*}\label{x2,t2_case2}
\begin{aligned}
    &x_{\ast 2} +\epsilon= \frac{1}{(1+V_{\ast 1})^2} \int_0^{t_{\ast 2}} \frac{1}{h(s)}ds,\\
    &x_{\ast 2}- x_{\ast 1} = \frac{1}{(1+v_+)(1+v_\thicksim)}\int_{t_{\ast 1}}^{t_{\ast 2}} \frac{1}{h(s)}ds.
\end{aligned} 
\end{equation*}
Consequently, the 1-shock wave $S_1^{(\ast 1)}$ penetrates the 1-rarefaction wave $R_1^{(-)}$ for $t>t_{\ast 2}$. On the penetration region, the states $(V,W)$ lying on the rarefaction wave $R_1^{(-)}$ must also lie on the curve of discontinuity of $S_1^{(\ast 1)}$. Then for $t>t_{\ast 2}$ the curve of discontinuity of $S_1^{(\ast 1)}$, namely $x=x(t)$ satisfies
\begin{equation}\label{S1_t>t2}
    \begin{cases}
        \frac{dx}{dt}= \frac{1}{(1+V_{\ast 3})(1+V)h(t)},\\
        \frac{W}{V}=\frac{W_{\ast 3}}{V_{\ast 3}},\\
        x +\epsilon= \frac{1}{(1+V)^2} \int_0^{t} \frac{1}{h(s)}ds,\\
        x(t_{\ast 2})=x_{\ast 2},\; v_\thicksim=V_{\ast 1}\le V< V_{\ast 3}=v_+.
    \end{cases}
\end{equation}
From \eqref{S1_t>t2}, we obtain by a simple calculations
\begin{equation} \label{dV/dt}
    \frac{dV}{dt}= \frac{(V_{\ast 3} -V)(1+V)}{2(1+V_{\ast 3})h(t)\int_0^t\frac{ds}{h(s)}}>0,
\end{equation}
which implies the value of $V$ increases during the penetration process. Integrating \eqref{dV/dt} over $[t_{\ast 2}, t]$, one obtains
\begin{equation*}
    \frac{2(1+V_{\ast 3})}{1+V_{\ast 3}}\ln{\frac{(1+V)(V_{\ast 3} -V_{\ast 1})}{(1+V_{\ast 1})(V_{\ast 3}- V)}}= \ln{\int_0^{t}\frac{ds}{h(s)}} - \ln{\int_0^{t_{\ast 2}}\frac{ds}{h(s)}},
\end{equation*}
which implies 
\begin{equation}\label{t_S1_1}
    \int_0^{t}\frac{1}{h(s)}ds= \frac{(1+V)^2(V_{\ast 3} -V_{\ast 1})^2}{(1+V_{\ast 1})^2(V_{\ast 3}- V)^2} \int_0^{t_{\ast 2}}\frac{1}{h(s)} ds
\end{equation}
The curve $x = -\epsilon+\frac{1}{(1+V)^2} \int_0^{t} \frac{1}{h(s)}ds$ together with \eqref{t_S1_1} gives the shock curve $S_1^{(\ast 1)}$ during the penetration for $t>t_{\ast 2}$, on which $V$ varies between $V_{\ast 1} \le V< V_{\ast 3}$, and $x(t_{\ast 2})=x_{\ast 2}$. The shock curve $S_1^{(\ast 1)}$ can fully penetrate the rarefaction $R_1^{(-)}$ and it depends on the value of $V_{\ast 3}=v_+$. Thus, depending on the value of $v_+$, we have two subcases: Subcase 2.1. $v_-\le v_+$ and Subcase 2.2. $v_+< v_-$. Now, we discuss these subcases separately.

\bigskip

\noindent \textbf{Subcase 2.1.} $v_- < v_+$.

Since the state $(V, W)$ also lies the rarefaction wave $R_1^{(-)}$ then we have $V\le v_-$. This implies $V\le v_- < v_+$. Therefore, the shock $S_1^{(\ast 1)}$ completely overtakes the rarefaction wave $R_1^{(-)}$, and interact the tail (front back) of the rarefaction wave at a point, namely $(x_{\ast 3}, t_{\ast 3})$, and the interaction point is determined from the following equations
\begin{equation*}
    \begin{aligned}
        &x_{\ast 3} +\epsilon= \frac{1}{(1+v_-)^2} \int_0^{t_{\ast 3}} \frac{1}{h(s)}ds,\\
        & \int_0^{t_{\ast 3}}\frac{1}{h(s)}ds= \frac{(1+v_-)^2(v_+-v_\thicksim)^2}{(1+v_\thicksim)^2(v_+-v_-)^2}\int_0^{t_{\ast 2}}\frac{1}{h(s)}ds.
    \end{aligned}
\end{equation*}
Now, a new Riemann problem is formed at $(x_{\ast 3}, t_{\ast 3})$ with the constant states $(v_-, w_-)$ and $(V_{\ast 3}, W_{\ast 3})$ (see Figure \ref{Fig_case2.1}), which satisfies $V_{\ast 3}> v_-$. Thus, the states $(v_-, w_-)$ and $(V_{\ast 3}, W_{\ast 3})$ is connected by a 1-shock wave $S_1^{(\ast 3)}$, which propagates from $(x_{\ast 3}, t_{\ast 3})$ with speed $s_1^{(\ast 3)}$, given by 
\begin{equation*}
    s_1^{(\ast 3)}= \frac{dx}{dt}= \frac{1}{(1+V_{\ast 3})(1+v_-)h(t)}= \frac{1}{(1+v_+)(1+v_-)h(t)}.
\end{equation*}
Comparing the speeds of $S_1^{(\ast 3)}$ and $J_2^{(\ast 3)}$, one obtain
\begin{equation*}
    s_1^{(\ast 3)}- \sigma_2^{(\ast 3)} = \frac{1}{(1+v_+)(1+v_-)h}- \frac{1}{(1+v_+)h}=\frac{-v_-}{(1+v_+)(1+v_-)h}<0,
\end{equation*}
which implies $s_1^{(\ast 3)} < \sigma_2^{(\ast 3)}$, and therefore, $S_1^{(\ast 3)}$ will never interact $J_2^{(\ast 3)}$. Thus, there are no further interactions. Hence, the solution to the interaction problem \eqref{sys} and \eqref{PRP}, for $t>t_{\ast 3}$ is of the form
\begin{equation}\label{sol_case2.1}
    \begin{aligned}
    (v_- h(t), w_- h(t))+ S_1^{(\ast 3)}+(V_{\ast 3}h(t), W_{\ast 3}h(t))+J_2^{(\ast 3)}+ (V_{\ast 2}h(t), W_{\ast 2}h(t)) \\+ J_2^{(\ast 2)}+ (v_+ h(t), w_+ h(t)).
    \end{aligned}
\end{equation}
Now, if we take $\epsilon\to 0$, then the points $(x_{\ast i}, t_{\ast i})$ ($i=1,2,3$) of interactions converge to origin $(0,0)$. Moreover, the curve of contact discontinuities $J_2^{(\ast 2)}$ and $J_2^{(\ast 3)}$ coincides and converges to 
\begin{equation}\label{J_2}
    J_2: x= \frac{1}{(1+v_+)} \int_0^t \frac{1}{h(s)}ds,
\end{equation}
while the shock curves $S_1^{(\ast 1)}$, $S_1^{(\ast 3)}$ coincides and becomes a 1-shock curve propagating from $(0,0)$ as follows 
\begin{equation*}
    S_1: x= \frac{1}{(1+v_+)(1+v_-)} \int_{0}^t \frac{1}{h(s)}ds.
\end{equation*}
Thus, the solution \eqref{sol_case2.1} of the perturbed Riemann problem tends to 
\begin{equation*}
    (v_- h(t), w_- h(t))+ S_1+(V_{\ast 3}h(t), W_{\ast 3}h(t))+J_2+ (v_+ h(t), w_+ h(t)),
\end{equation*}
which is exactly the solution of the Riemann problem \eqref{sys} and \eqref{RP}. Hence, in this subcase, the Riemann solution is stable with respect to the small perturbation of the initial data.

\bigskip

\noindent \textbf{Subcase 2.2.} $v_- > v_+$.

\begin{figure}[ht!]  
    \centering
\includegraphics[width=0.75\linewidth]{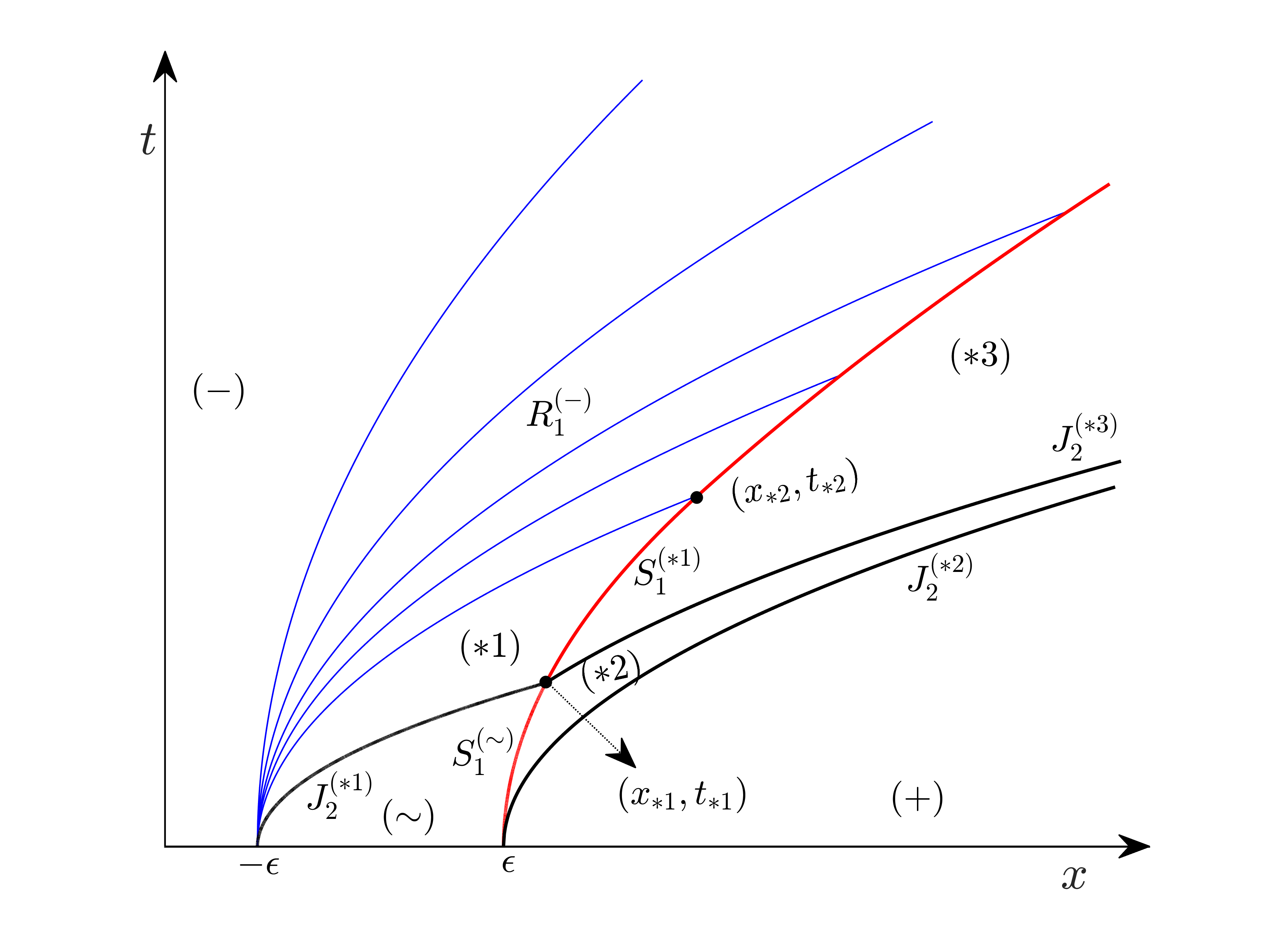}
\caption{Interactions of waves for the subcase 2.2, i.e., when $0<v_\thicksim<v_+<v_-$.}
\label{Fig_case2.2}
\end{figure}

In the penetration region $t>t_{\ast 2}$, the shock curve $S_1^{(\ast 1)}$ is given in \eqref{S1_t>t2} for $V_{\ast 1} \le V< V_{\ast 3}$. On account of \eqref{t_S1_1}, implies that $t\to \infty$ as $V\to V_{\ast 3}$, and thus in a finite time, the shock curve $S_1^{(\ast 1)}$ cannot fully penetrate the rarefaction wave $R_1^{(-)}$ since $V_{\ast 3}=v_+< v_-$ (see Figure \ref{Fig_case2.2}). Furthermore, in this subcase, the curve $x +\epsilon= \frac{1}{(1+V_{\ast 3})^2} \int_0^{t} \frac{1}{h(s)}ds$ becomes asymptote of $S_1^{(\ast 1)}$. Therefore, for sufficiently large $t> t_{\ast 2}$, the solution of the perturbed Riemann problem takes the form
\begin{equation}\label{sol_case2.2}
    (v_- h(t), w_- h(t))+ R_1^{(-)}+(V_{\ast 3}h(t), W_{\ast 3}h(t))+ J_2^{(\ast 3)}+ (V_{\ast 2}h(t), W_{\ast 2}h(t)) \\+ J_2^{(\ast 2)}+ (v_+ h(t), w_+ h(t)).
\end{equation}
Now, taking $\epsilon\to 0$, the contact discontinuities $J_2^{(\ast 3)}$ and $J_2^{(\ast 2)}$ propagating with the same speed must coincide and tend to a single contact discontinuity $J_2$, given in \eqref{J_2}. Moreover, the shock curve $S_1^{(\ast 1)}$ coincides with the wave front of $R_1^{(-)}$, and for sufficiently large time, the right portion of the rarefaction of the asymptote vanishes, and the left portion remains. More precisely, as $\epsilon\to 0$, the rarefaction wave 
$R_1^{(-)}$, given in \eqref{R1(-)} tends to 1-rarefaction wave $R_1$ as follows
\begin{equation} \label{R_1}
    R_1: 
    \begin{cases}
        \xi = \frac{1}{(1+V)^2},\\
        \frac{W}{V}=\frac{w_-}{v_-},\\
        x(0)=0, \quad V_{\ast 3} \le V \le v_-.
    \end{cases}
\end{equation}
Therefore, the solution \eqref{sol_case2.2} tends to the solution of the Riemann problem \eqref{sys} and \eqref{RP} with the following form
\begin{equation*}
    (v_- h(t), w_- h(t))+ R_1+(V_{\ast 3}h(t), W_{\ast 3}h(t))+J_2+ (v_+ h(t), w_+ h(t)).
\end{equation*}
Hence, the Riemann solution is stable in this case as well under the small perturbation of the initial data.

\subsection{\textbf{Case 3:} $0<v_-< v_\thicksim$ and $v_+< v_\thicksim$.}
\begin{figure}[ht!]  
    \centering
\includegraphics[width=0.75\linewidth]{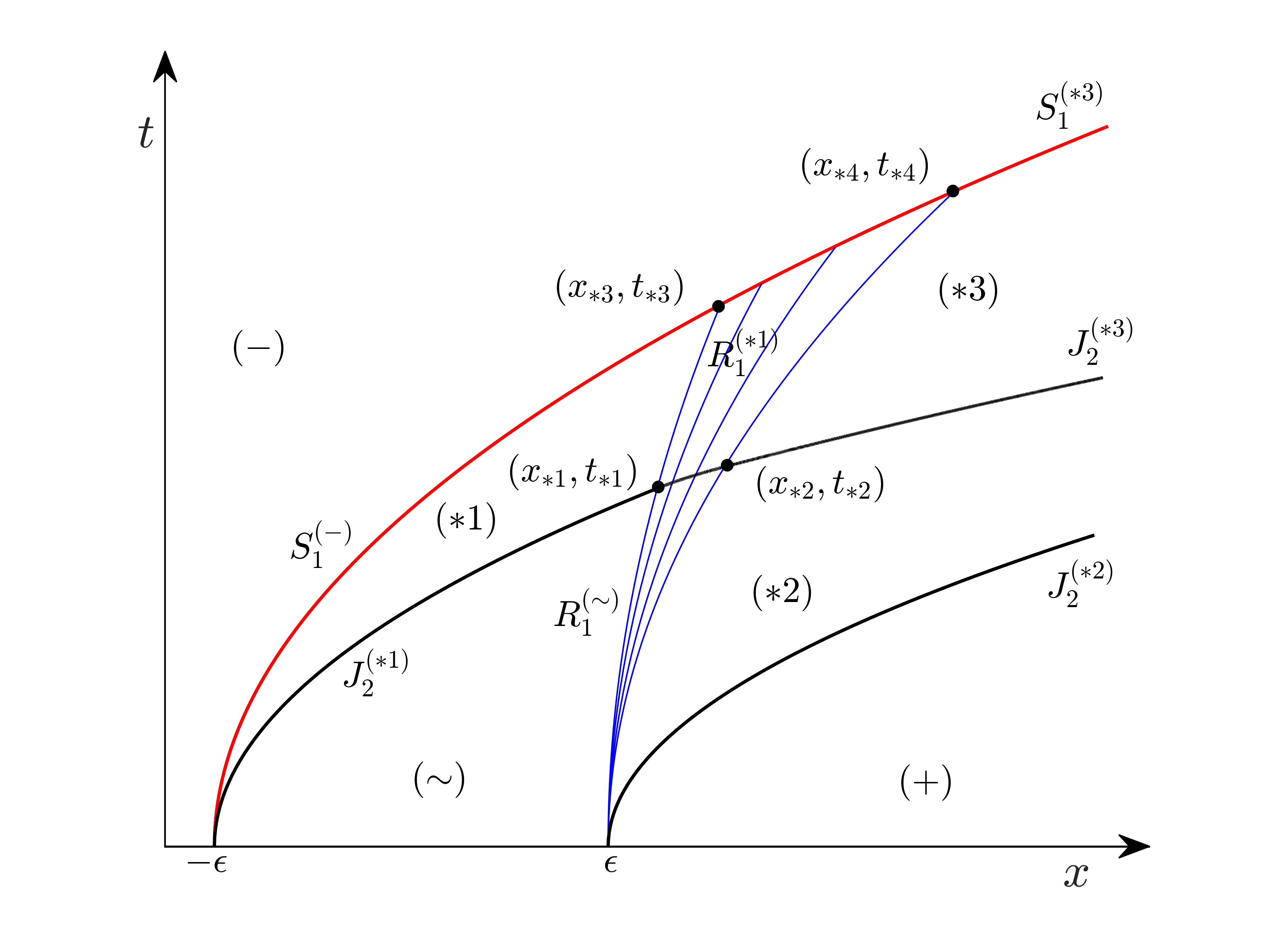}
\caption{Interactions of waves for the subcase 3.1, i.e., when $0< v_-< v_+<v_\thicksim$.}
\label{Fig_case3.1}
\end{figure}

\begin{figure}[ht!]  
    \centering
\includegraphics[width=0.75\linewidth]{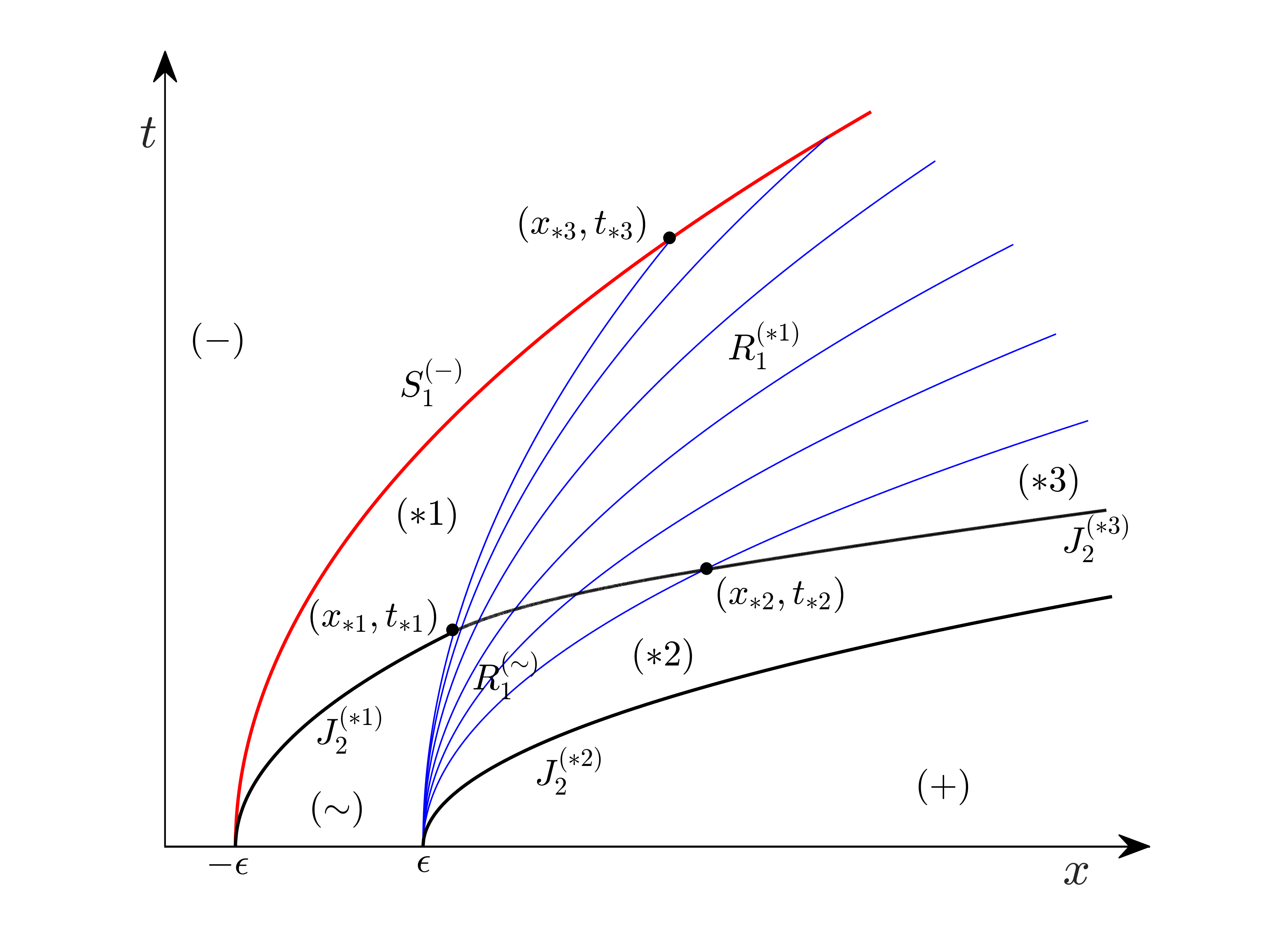}
\caption{Interactions of waves for the subcase 3.2, i.e., when $0< v_-<v_\thicksim$ and $v_+< v_-$.}
\label{Fig_case3.2}
\end{figure}

In this case, the solution to the local Riemann problem at $x=-\epsilon$ is connected by 
\begin{equation*}
    (v_- h(t), w_- h(t))+ S_1^{(-)}+(V_{\ast 1}h(t), W_{\ast 1}h(t))+J_2^{(\ast 1)}+ (v_\thicksim h(t), w_\thicksim h(t)),
\end{equation*}
while the solution to the Riemann problem at $x=\epsilon$ is given by 
\begin{equation*}
    (v_\thicksim h(t), w_\thicksim h(t))+ R_1^{(\thicksim)}+(V_{\ast 2}h(t), W_{\ast 2}h(t))+J_2^{(\ast 2)}+ (v_+ h(t), w_+ h(t)).
\end{equation*}
Here, $S_1^{(-)}$ and $J_2^{(\ast 1)}$ are 1-shock wave and 2-contact discontinuity with speeds $s_1^{(-)}= \frac{1}{(1+v_-)(1+v_\thicksim)h}$ and $\sigma_2^{(\ast 1)}=\frac{1}{(1+v_\thicksim)h}$ respectively, separated by the intermediate state $(V_{\ast 1}h(t), W_{\ast 1}h(t))$, given in \eqref{interm_V1,W1}. Moreover, $R_1^{(\thicksim)}$ is a 1-rarefaction wave and $J_2^{(\ast 2)}$ is a 2-contact discontinuity with speed $\sigma_2^{(\ast 2)}=\frac{1}{(1+v_+)h}$, separated by the intermediate state $(V_{\ast 2}h(t), W_{\ast 2} h(t))$, where
\begin{equation*} \label{R1(sim)}
    R_1^{(\thicksim)}: 
    \begin{cases}
        \xi = \frac{1}{(1+V)^2},\\
        \frac{W}{V}=\frac{w_\thicksim}{v_\thicksim}, \quad V_{\ast 2} \le V \le v_-,
    \end{cases}
\end{equation*}
and $(V_{\ast 2}h(t), W_{\ast 2} h(t))$ satisfies
\begin{equation*}
    (V_{\ast 2}h(t), W_{\ast 2} h(t))= h(t)\Big(v_+,\frac{v_+ w_\thicksim}{v_\thicksim}\Big).
\end{equation*}
Thus, for sufficiently small time $t$, the solution of the Cauchy problem \eqref{sys} and \eqref{PRP} takes the form
\begin{equation*}
\begin{aligned}
    (v_- h(t), w_- h(t))+ S_1^{(-)}+(V_{\ast 1}h(t), W_{\ast 1}h(t))+J_2^{(\ast 1)}+ (v_\thicksim h(t), w_\thicksim h(t))+ R_1^{(\thicksim)}\\+ (V_{\ast 2}h(t), W_{\ast 2}h(t)) + J_2^{(\ast 2)}+ (v_+ h(t), w_+ h(t)),
    \end{aligned}
\end{equation*}
which is depicted in Figures \ref{Fig_case3.1} and \ref{Fig_case3.2}. However, we are interested in a global solution for all $t$, therefore, it is necessary to investigate the wave interactions. In fact, in a finite time, $J_2^{(\ast 1)}$ interacts with the tail (wave back) of the rarefaction wave $R_1^{(\thicksim)}$, since $\sigma_2^{(\ast 1)}> \theta_1^{(\thicksim)}$, where $\theta_1^{(\thicksim)}$ denotes the speed of wave back of $R_1^{(\thicksim)}$, given by $\theta_1^{(\thicksim)}= \frac{1}{(1+v_\thicksim)^2 h}$. Let the point of interaction be $(x_{\ast 1}, t_{\ast 1})$. Then, we have
\begin{equation}\label{x1,t1_case3}
    \begin{aligned}
        &x_{\ast 1}+\epsilon= \frac{1}{(1+v_\thicksim)} \int_0^{t_{\ast 1}} \frac{1}{h(s)}ds,\\
        & x_{\ast 1} -\epsilon= \frac{1}{(1+v_\thicksim)^2} \int_0^{t_{\ast 1}} \frac{1}{h(s)}ds.
    \end{aligned}
\end{equation}
Thus, the contact discontinuity $J_2^{(\ast 1)}$ begins to penetrate the rarefaction wave $R_1^{(\thicksim)}$ at time $t_{\ast 1}$. Moreover, for time $t> t_{\ast 1}$, the curve of discontinuity $J_2^{(\ast 1)}$ enters the rarefaction region, and it penetrates the rarefaction over time. Therefore, discontinuity $J_2^{(\ast 1)}$ must change its speed during the process of penetration. Let $x=x(t)$ be the curve of discontinuity $J_2^{(\ast 1)}$ in the penetration region, and $(V_l,W_l)$, $(V_r,W_r)$ be respectively left and right states across the curve $x=x(t)$. Then, the penetrating curve $x=x(t)$ is determined from the following
\begin{equation}\label{J1_t>t1}
    \begin{cases}
        \frac{dx}{dt}= \frac{1}{(1+V_l)h(t)}= \frac{1}{(1+V_r)h(t)},\\
        x -\epsilon= \frac{1}{(1+V_r)^2} \int_0^{t} \frac{1}{h(s)}ds,\\
        \frac{W_r}{V_r}=\frac{w_\thicksim}{v_\thicksim}, \\
        \frac{W_l}{V_l}=\frac{W_{\ast 1}}{V_{\ast 1}},\\
        x(t_{\ast 1})=x_{\ast 1},\; v_+= V_{\ast 2}\le V_r \le v_\thicksim.
    \end{cases}
\end{equation}
Taking the derivative of the second equation of \eqref{J1_t>t1} with respect to $t$ and using the first equation leads to 
\begin{equation}\label{dVr/dt}
    \frac{dV_r}{dt} = -\frac{V_r(1+V_r)}{2h(t)\int_0^t\frac{1}{h(s)}ds} <0,
\end{equation}
and thus, $V_r$ decreases over time during penetration. Integrating \eqref{dVr/dt} over $[t_{\ast 1}, t]$ yields
\begin{equation*}
    2\Big(\ln{\frac{1+V}{V}}- \ln{\frac{1+v_\thicksim}{v_\thicksim}}\Big) = \ln{\int_0^t\frac{1}{h(s)}ds} - \ln{\int_0^{t_{\ast 1}}\frac{1}{h(s)}ds},
\end{equation*}
which on account of \eqref{x1,t1_case3} implies
\begin{equation}\label{exp_t}
    \int_0^{t}\frac{1}{h(s)}ds= 2\epsilon v_\thicksim \Big(\frac{1+V}{V}\Big)^2, \quad V_{\ast 2}\le V\le v_\thicksim.
\end{equation}

Hence, the contact discontinuity $J_2^{(\ast 1)}$ completely penetrates the rarefaction wave at a finite time, say $t_{\ast 2}$. The time $t_{\ast 2}$, that is, the time of interaction of $J_2^{(\ast 1)}$ with the wave front of rarefaction fan is determined from \eqref{exp_t} by taking $V\to V_{\ast 2}$, and the corresponding $x$-value, namely $x_{\ast 2}$ can be obtained from $x_{\ast 2}=\epsilon+ \frac{1}{(1+V_{\ast 2})^2} \int_0^{t_{\ast 2}} \frac{1}{h(s)}ds$. Therefore, for time $t> t_{\ast 2}$, we have a new contact discontinuity, denoted by $J_2^{(\ast 3)}$ propagating with speed $\sigma_2^{(\ast 3)}= \frac{1}{(1+V_{\ast 2})h}=\frac{1}{(1+v_+)h}$. Moreover, above the penetration region, we have a new rarefaction wave, denoted by $R_1^{(\ast 1)}$. Let us denote the new intermediate state between $R_1^{(\ast 1)}$ and $J_2^{(\ast 3)}$ by $(V_{\ast 3}, W_{\ast 3})$. Then, $V_{\ast 3}=V_{\ast 2}=v_+$ and $W_{\ast 3}= \frac{V_{\ast 2}W_{\ast 1}}{V_{\ast 1}}= \frac{v_+ w_-}{v_-}$. Since $\sigma_2^{(\ast 3)}= \sigma_2^{(\ast 2)}$, i.e., the speeds of contact discontinuities $J_2^{(\ast 3)}$ and $J_2^{(\ast 2)}$ are the same, then they are parallel to each other and will not interact forever. 

On the other hand, the speed of the tail (wave back) of $R_1^{(\ast 1)}$ is given by $\theta_1^{(\ast 1)}= \frac{1}{(1+V_{\ast 1})^2 h}= \frac{1}{(1+v_\thicksim)^2 h}$, and the speed of $S_1^{(-)}$ is $s_1^{(-)}= \frac{1}{(1+v_-)(1+v_\thicksim)h}$. Comparing one obtains 
\begin{equation*}
    s_1^{(-)}- \theta_1^{(\ast 1)}= \frac{v_\thicksim- v_-}{(1+v_-)(1+v_\thicksim)^2 h}> 0,
\end{equation*}
and therefore, $S_1^{(-)}$ must interact with $R_1^{(\ast 1)}$ in a finite time. Let $(x_{\ast 3}, t_{\ast 3})$ be the point of interaction between $S_1^{(-)}$ and the wave back of $R_1^{(\ast 1)}$. Then, we have 
\begin{equation*}
\begin{aligned}
&x_{\ast 3}+ \epsilon = \frac{1}{(1+V_{\ast 1})(1+v_-)}\int_{0}^{t_{\ast 3}} \frac{1}{h(s)}ds,\\
    &x_{\ast 3} - x_{\ast 1}= \frac{1}{(1+V_{\ast 1})^2} \int_{t_{\ast 1}}^{t_{\ast 3}} \frac{1}{h(s)}ds.
\end{aligned} 
\end{equation*}
Again, we have the following two subcases depending on the values of $v_-$ and $v_+$. These two subcases are analogous to the subcases in Case 2. We will not discuss these subcases in detail now.

\noindent \textbf{Subcase 3.1.} $v_-< v_+$.

The 1-shock curve $S_1^{(-)}$ begins to penetrate the rarefaction wave $R_1^{(\ast 1)}$ at time $t_{\ast 3}$. The curve, namely $x=x(t)$ of $S_1^{(-)}$ in the penetration region, can be derived in a similar manner for $V_{\ast 3}< V\le V_{\ast 1}$ and $v_-< V$, where $(V, W)$ is any state on $x=x(t)$. So, in this subcase, we have $v_-< v_+=V_{\ast 3}< V\le V_{\ast 1}$. Thus, 1-shock curve $x=x(t)$ must interact with the wave front of $R_1^{(\ast 1)}$ in a finite time, say $t_{\ast 4}$, can be derived analogously. This means the 1-shock wave $S_1^{(-)}$ fully penetrate the rarefaction wave in time $t_{\ast 4}$, and produces a new 1-shock wave $S_1^{(\ast 3)}$ propagating with speed $s_1^{\ast 3}= \frac{1}{(1+v_-)(1+V_{\ast 3})h}$ for $t> t_{\ast 4}$; see Figure \ref{Fig_case3.1}. Hence, in this subcase, the solution to the Cauchy problem \eqref{sys} and \eqref{PRP} for $t> t_{\ast 4}$ is of the form
\begin{equation}\label{sol_case3.1}
    \begin{aligned}
    (v_- h(t), w_- h(t))+ S_1^{(\ast 3)}+(V_{\ast 3}h(t), W_{\ast 3}h(t))+J_2^{(\ast 3)}+ (V_{\ast 2}h(t), W_{\ast 2}h(t)) \\+ J_2^{(\ast 2)}+ (v_+ h(t), w_+ h(t)).
    \end{aligned}
\end{equation}
As $\epsilon \to 0$, the solution \eqref{sol_case3.1} becomes the solution to the Riemann problem \eqref{sys} and \eqref{RP}. Hence, the Riemann solution is stable.

\bigskip

\noindent \textbf{Subcase 3.2.} $v_-> v_+$.

In this subcase, the shock curve $S_1^{(-)}$ penetrates the rarefaction wave in the region $V_{\ast 3}=v_+<v_-< V\le V_{\ast 1}$. Thus, it cannot penetrate the whole rarefaction wave; see Figure \ref{Fig_case3.2}. In fact, the 1-characteristic curve inside the rarefaction fan $R_1^{(\ast 1)}$ with speed $\frac{1}{(1+v_-)^2 h(t)}$ becomes the asymptote of the shock curve $x=x(t)$. So, the shock curve cannot penetrate the right part of the rarefaction fan beside the asymptote. Therefore, for sufficiently large time $t>t_{\ast 3}$, the solution of the perturbed Riemann problem \eqref{sys} and \eqref{PRP} can be expressed as follows
\begin{equation}\label{sol_case3.2}
    (v_- h(t), w_- h(t))+ R_1^{(\ast 1)}+(V_{\ast 3}h(t), W_{\ast 3}h(t))+ J_2^{(\ast 3)}+ (V_{\ast 2}h(t), W_{\ast 2}h(t)) + J_2^{(\ast 2)}+ (v_+ h(t), w_+ h(t)).
\end{equation}
Similarly, as $\epsilon\to 0$, we can show that the solution \eqref{sol_case3.2} tends to the corresponding Riemann solution of \eqref{sys} and \eqref{RP}. As a conclusion, the Riemann solution is stable under the perturbation.

\subsection{\textbf{Case 4:} $v_+< v_\thicksim< v_-$.}

\begin{figure}[ht!]  
    \centering
\includegraphics[width=0.75\linewidth]{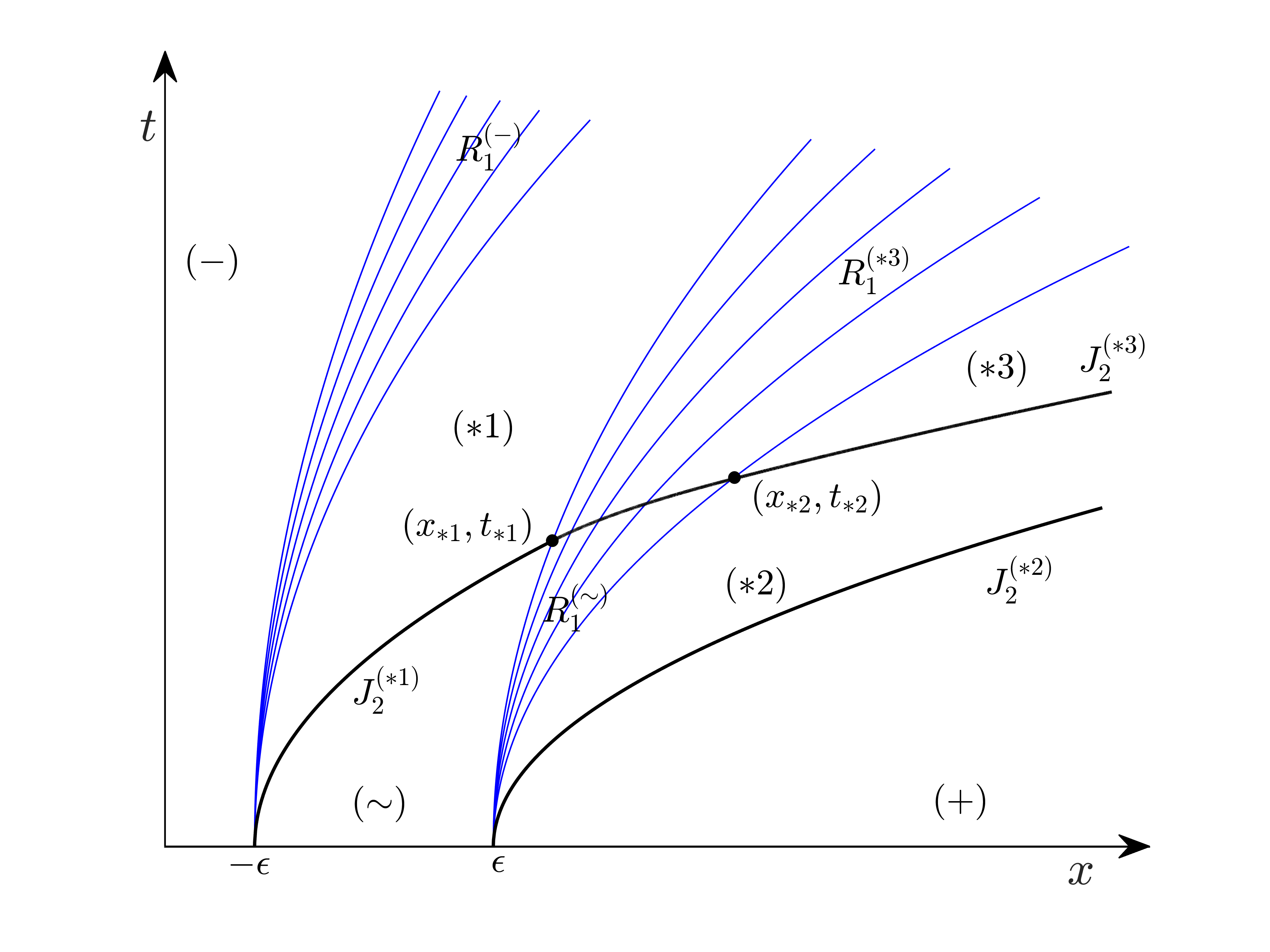}
\caption{Interactions of waves for the case 4, i.e., when $v_+< v_\thicksim< v_-$.}
\label{Fig_case4}
\end{figure}

Here, for a sufficiently small time $t$, the solution to the Cauchy problem \eqref{sys} and \eqref{PRP} can be expressed as follows
\begin{equation*}
\begin{aligned}
    (v_- h(t), w_- h(t))+ R_1^{(-)}+(V_{\ast 1}h(t), W_{\ast 1}h(t))+J_2^{(\ast 1)}+ (v_\thicksim h(t), w_\thicksim h(t))+ R_1^{(\thicksim)}\\+ (V_{\ast 2}h(t), W_{\ast 2}h(t)) + J_2^{(\ast 2)}+ (v_+ h(t), w_+ h(t)),
    \end{aligned}
\end{equation*}
where $R_1^{(-)}$ and $R_1^{(\thicksim)}$ are rarefaction waves, respectively given by  
\begin{equation*}
    R_1^{(-)}: 
    \begin{cases}
        \xi = \frac{1}{(1+V)^2},\\
        \frac{W}{V}=\frac{w_-}{v_-},\\
        x(0)=-\epsilon, \quad V_{\ast 1} \le V \le v_-, 
    \end{cases}
    \quad R_1^{(\thicksim)}: 
    \begin{cases}
        \xi = \frac{1}{(1+V)^2},\\
        \frac{W}{V}=\frac{w_\thicksim}{v_\thicksim},\\
        x(0)=\epsilon, \quad V_{\ast 2} \le V \le v_\thicksim,
    \end{cases}
\end{equation*}
and $(V_{\ast 1}, W_{\ast 1})=(v_\thicksim,\frac{v_\thicksim w_-}{v_-})$ and $(V_{\ast 2}, W_{\ast 2})=(v_+,\frac{v_+ w_\thicksim}{v_\thicksim})$. 

Analogous to the previous Case 3, it can be shown that the contact discontinuity $J_2^{(\ast 1)}$ interacts and penetrates the whole rarefaction fan $R_1^{(\thicksim)}$ (see Figure \ref{Fig_case4}). Consequently, after the penetration, it produces a new contact discontinuity $J_1^{(\ast 3)}$ propagating from $(x_{\ast 2}, t_{\ast 2})$ and a new rarefaction wave $R_1^{(\ast 3)}$ separated by an intermediate state $(V_{\ast 3}, W_{\ast 3})h(t)$. The speeds of these contact discontinuities $J_1^{(\ast 3)}$ and $J_1^{(\ast 2)}$ are same, and therefore, they are parallel. Furthermore, the rarefaction waves $R_1^{(-)}$ and $R_1^{(\ast 3)}$ will never interact, as the speed of the head (wave front) of $R_1^{(-)}$ coincides with the speed of the tail (wave back) of $R_1^{(\ast 3)}$, that is $\frac{1}{(1+V_{\ast 1})^2 h(t)}$. Thus, there is no further interaction for $t> t_{\ast 2}$.
Hence, for any time $t> t_{\ast 2}$, the solution to the perturbed Riemann problem is of the following form
\begin{equation}\label{sol_case4}
\begin{aligned}
    (v_- h(t), w_- h(t))+ R_1^{(-)}+(V_{\ast 1}h(t), W_{\ast 1}h(t))+ R_1^{(\ast 3)}+(V_{\ast 3}h(t), W_{\ast 3}h(t))\\+ J_2^{(\ast 3)}+ (V_{\ast 2}h(t), W_{\ast 2}h(t))+ J_2^{(\ast 2)}+ (v_+ h(t), w_+ h(t)).
    \end{aligned}
\end{equation}

When $\epsilon\to 0$, all the points $(-\epsilon, 0)$, $(\epsilon, 0)$, $(x_{\ast 1}, t_{\ast 1})$ and $(x_{\ast 2}, t_{\ast 2})$ tend to $(0,0)$. Moreover, as $\epsilon\to 0$, the rarefaction waves $R_1^{(-)}$ and $R_1^{(\ast 3)}$ coincides and tends to a 1-rarefaction wave $R_1$ propagating from the origin, given in \eqref{R_1}, and the contact discontinuities $J_2^{(\ast 3)}$ and $J_2^{(\ast 2)}$ tends to a contact discontinuity $J_2$, given in \eqref{J_2}.
Therefore, taking $\epsilon\to 0$, the perturbed Riemann solution \eqref{sol_case4} converges to 
\begin{equation*}
    (v_- h(t), w_- h(t))+ R_1+(V_{\ast 3}h(t), W_{\ast 3}h(t))+J_2+ (v_+ h(t), w_+ h(t)),
\end{equation*}
which is a Riemann solution for \eqref{sys}. Hence, the solution of the Riemann problem \eqref{sys} and \eqref{RP} is stable.

\begin{remark}
    In all these four cases 1-4, we constructed the solutions of the Cauchy problem \eqref{sys} and \eqref{PRP} by analyzing the classical wave interactions, and consequently, proved that the Riemann solutions for \eqref{sys} are stable for all cases 1-4. Note that there will be no further cases of interactions of classical waves. Moreover, to establish the solution of the Cauchy problem for any values of initial data \eqref{PRP}, we need to consider other cases of interactions where at least one local Riemann problem in \eqref{PRP} is solved by a nonclassical delta shock wave, i.e., at least one of $v_-, v_\thicksim$ is $0$. In the subsequent sections, we shall now discuss all possible such cases to complete the proof of Theorem \ref{main_thm}.
\end{remark}

\subsection{\textbf{Case 5:} $0=v_\thicksim< v_\pm$.}

In this case, the first local Riemann problem at $(-\epsilon, 0)$ is solved by a 1-rarefaction wave $R_1^{(-)}$ and a 2-contact discontinuity $J_2^{(\ast 1)}$, moreover, the solution is as follows
\begin{equation*}
    (v_- h(t), w_- h(t)) + R_1^{(-)} + (V_{\ast 1}h(t), W_{\ast 1} h(t))+ J_2^{(\ast 1)}+ (0, w_\thicksim h(t)), 
\end{equation*}
where the intermediate state $(V_{\ast 1}h(t), W_{\ast 1} h(t))= h(t)\Big(v_\thicksim,\frac{v_\thicksim w_-}{v_-}\Big)=(0,0)$. Moreover, the speeds of $J_2^{(\ast 1)}$ and the wave front of $R_1^{(-)}$ coincides, that is $\sigma_2^{(\ast 1)}=\frac{1}{h(t)}$. Thus, the wave configuration $R_1^{(-)}+ J_2^{(\ast 1)}$ is a composite wave, as both the wave front and $J_2^{(\ast 1)}$ are supported on the same curve $x+\epsilon = \int_0^t \frac{1}{h(s)}ds$. Let us denote this composite wave as $R_1^{(-)}J_2^{(\ast 1)}$.

On the other hand, the local Riemann problem at $(\epsilon, 0)$ is solved by a delta shock wave, namely $\delta S$ as follows
\begin{equation*}
    (0, h(t) w_\thicksim)+ \delta S + (h(t) v_+, h(t) w_+),
\end{equation*}
where $\delta S$ is the delta shock wave supported on the curve $x-\epsilon=\theta \int_0^t \frac{1}{h(s)} ds$, on which the state $(v(x,t),w(x,t))$ is given by $(v(x,t),w(x,t))= (v_\delta(t), \alpha(t)\delta(x- \epsilon-\theta \int_0^t \frac{1}{h(s)} ds))$, with the strength $\alpha(t)=h(t)\frac{w_\thicksim v_+}{1+v_+}\int_0^t \frac{1}{h(s)} ds$, $\theta=\frac{1}{1+\frac{1}{h(t)}v_\delta}$, $v_\delta(t)=h(t)v_+$. 
Then, in this case, for a small enough time $t>0$, the solution to the perturbed Riemann problem is of the form 
\begin{equation*}\label{case5}
    (v_- h(t), w_- h(t))+ R_1^{(-)}J_2^{(\ast 1)} +(0, h(t) w_\thicksim)+ \delta S + (h(t) v_+, h(t) w_+).
\end{equation*}

The speed of propagation of $\delta S$ is $\sigma_\delta= \frac{\theta}{h(t)}= \frac{1}{(1+v_+)h(t)}$. Comparing the speeds of $J_2^{(\ast 1)}$ and $\delta S$, one obtains $\sigma_2^{(\ast 1)}- \sigma_\delta= \frac{v_+}{(1+v_+)h}>0$, and therefore, they must interact in finite time, say $t_{\ast 1}$. Let $x_{\ast 1}=x(t_{\ast 1})$. Then, the point of interaction $(x_{\ast 1}, t_{\ast 1})$ satisfies 
\begin{equation}\label{x1,t1_case5}
    \begin{aligned}
        &x_{\ast 1}+\epsilon= \int_0^{t_{\ast 1}} \frac{1}{h(s)}ds,\\
        & x_{\ast 1} -\epsilon= \frac{1}{(1+v_+)} \int_0^{t_{\ast 1}} \frac{1}{h(s)}ds.
    \end{aligned}
\end{equation}
Exploiting \eqref{x1,t1_case5}, we obtain $x_{\ast 1}= \epsilon\frac{2+v_+}{v_+}$, and $t_{\ast 1}$ is determined by $\int_0^{t_{\ast 1}} \frac{1}{h(s)}ds= 2\epsilon\frac{1+v_+}{v_+}$. Then the strength of the delta shock wave $\delta S$ at the point $(x_{\ast 1}, t_{\ast 1})$ is given by 
\begin{equation*}\label{alpha_t1}
    \alpha(t_{\ast 1})=h(t_{\ast 1})\frac{w_\thicksim v_+}{1+v_+}\int_0^{t_{\ast 1}} \frac{1}{h(s)} ds= 2\epsilon w_\thicksim h(t_{\ast 1}).
\end{equation*}
Since $v_-\neq 0$, then for time $t>t_{\ast 1}$, the solution does not contain a new delta shock wave. As done in \cite{Liu&Smoller, Shen&Sun_09}, let us approximate the rarefaction wave $R_1^{(-)}$ by a set of nonphysical shock waves. Thus, at $(x_{\ast 1}, t_{\ast 1})$, we have a new local Riemann problem as follows 
\begin{equation} \label{rp_x1,t1}
    v\big\lvert_{t=t_{\ast 1}}= \Bigg\{\begin{array}{lr}
        V h(t), \quad &x< x_{\ast 1},\\
        v_+ h(t), \quad &x> x_{\ast 1},
    \end{array},\quad
    w\big\lvert_{t=t_{\ast 1}}= \Bigg\{\begin{array}{lr}
        W h(t), \quad &x< x_{\ast 1},\\
        w_+ h(t), \quad &x> x_{\ast 1},
    \end{array}\Bigg\} + \alpha(t_{\ast 1})\delta_{(x_{\ast 1}, t_{\ast 1})}.
\end{equation}
The left state $(V,W)$ in \eqref{rp_x1,t1}, which is continuously varying, is given by 
\begin{equation*}
    \begin{cases}
        x+\epsilon = \frac{1}{(1+V)^2}\int_0^{t} \frac{1}{h(s)}ds,\\
        \frac{W}{V}=\frac{w_-}{v_-}, \quad 0< V < v_-,
    \end{cases}
\end{equation*}
which implies 
\begin{equation*}
    (V, W)= \bigg(V, \frac{V w_-}{v_-}\bigg)= \left(\sqrt{\frac{\int_0^{t} \frac{1}{h(s)}ds}{x+\epsilon}}-1, \frac{w_-}{v_-}\Bigg(\sqrt{\frac{\int_0^{t} \frac{1}{h(s)}ds}{x+\epsilon}}-1\Bigg)\right).
\end{equation*}

\begin{figure}[ht!]  
    \centering
\includegraphics[width=0.75\linewidth]{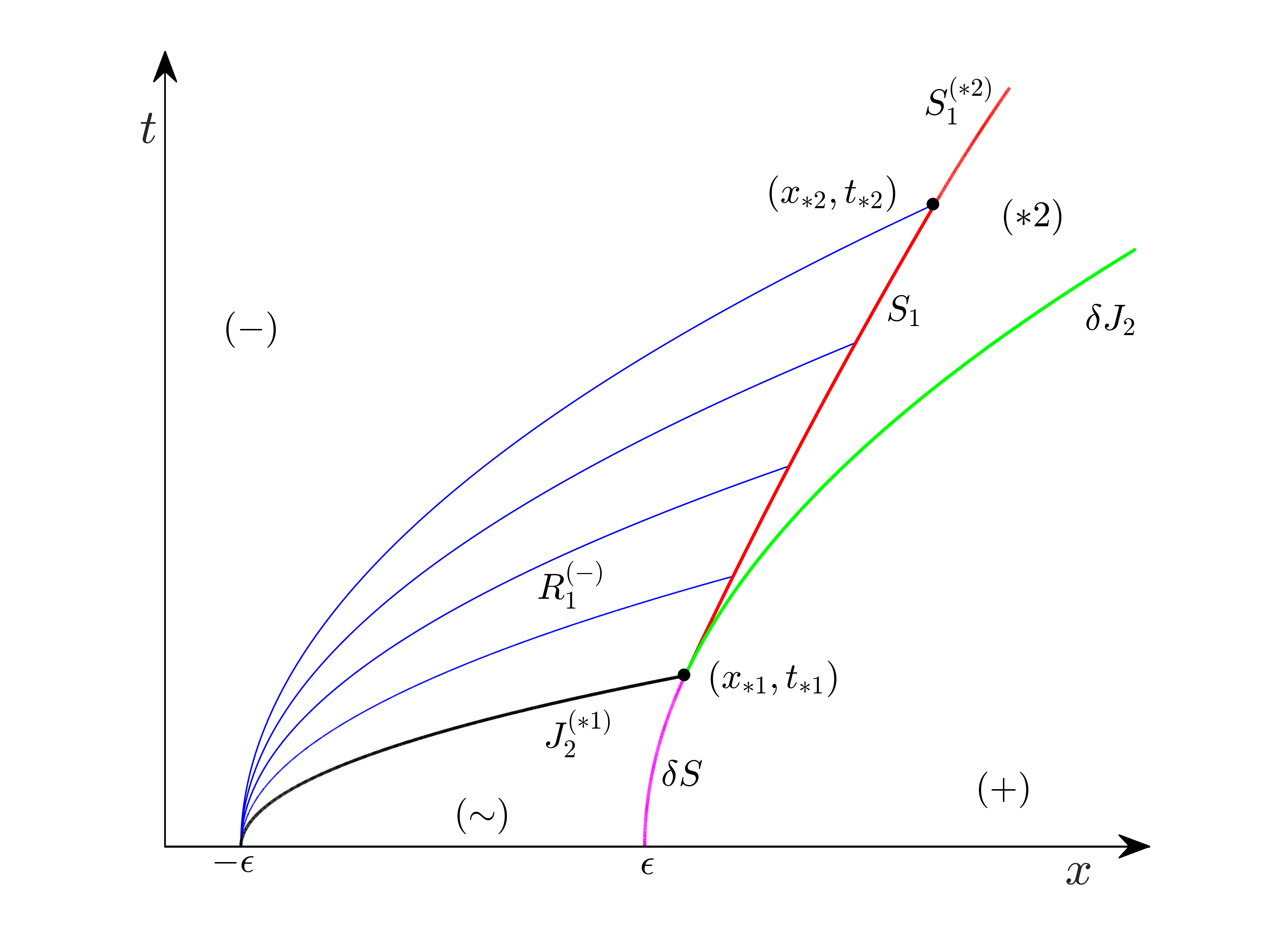}
\caption{Interactions of waves for the subcase 5.1, i.e., when $0= v_\thicksim< v_-<v_+$.}
\label{Fig_case5.1}
\end{figure}

Now, we claim that after the interaction, i.e., for $t> t_{\ast 1}$, the delta shock wave $\delta S$ splits into a shock wave supported on a curve $\Upsilon_1$ and a delta contact discontinuity supported on a curve $\Upsilon_2$, where we need to determine the curves $\Upsilon_i$ ($i=1,2$) as well. 

Indeed, the solution of \eqref{sys} and \eqref{rp_x1,t1}, i.e., the Riemann problem at $(x_{\ast 1}, t_{\ast 1})$ can be established in the following form
\begin{equation}\label{sol_rp_x1t1}
\begin{aligned}
    &v(x,t)=  \Bigg\{\begin{array}{lr}
        V h(t), \quad &x< x(t),\\
        v_+ h(t), \quad &x> x(t),
    \end{array},\\
    &w(x,t)= \left\{\begin{array}{lr}
        \frac{V w_-}{v_-} h(t), \quad &x< x(t),\\
        \frac{v_+ w_-}{v_-}h(t), \quad &x(t)<x < x_\delta(t),\\
        w_+ h(t), \quad &x> x_\delta(t),
    \end{array}\right\} + \alpha(t_{\ast 1})h(t)\delta(x-x_\delta(t)),
    \end{aligned}
\end{equation}
where $\Upsilon_1: x=x(t)$ is the shock curve in the local neighbourhood of $(x_{\ast 1}, t_{\ast 1})$ such that $x(t)=x_{\ast 1}+ \frac{1}{(1+v_+)(1+V)} \int_{t_{\ast 1}}^t \frac{1}{h(s)}ds$, and $\Upsilon_2: x=x_\delta(t)$ can be expressed as  $x_\delta(t)= \epsilon + \frac{1}{(1+v_+)}\int_0^{t} \frac{1}{h(s)}ds$, $t> t_{\ast 1}$. So, to prove \eqref{sol_rp_x1t1} is the solution of the initial value problem \eqref{sys} and \eqref{rp_x1,t1}, we need to verify that it satisfies the following weak formulation of \eqref{sys} in the sense of distribution (see \cite{Rcruz_RM_WN})
\begin{subequations} \label{weak_form}
\begin{align}
     &\langle v,\varphi_t \rangle + \langle \frac{v}{h+v}, \varphi_x \rangle = \langle \sigma v, \varphi \rangle,\label{wkfm_1}\\
            &\langle w,\varphi_t \rangle + \langle \frac{w}{h+v}, \varphi_x \rangle = \langle \sigma w, \varphi \rangle \label{wkfm_2},
\end{align}
\end{subequations}
for any $\varphi\in C_0^\infty(\mathbb{R}\times \mathbb{R}_+)$. If $\supp \varphi$ does not cut $\{(x,t): x=x_\delta(t), t> t_{\ast 1}\}$, i.e., $\supp \varphi \cap\{(x,t): x=x_\delta(t), t> t_{\ast 1}\}=\emptyset$, then one can easily verify that \eqref{sol_rp_x1t1} satisfies the weak formulation \eqref{weak_form}. On the other hand, if $\supp \varphi \cap\{(x,t): x=x_\delta(t), t> t_{\ast 1}\}\neq \emptyset$, this case needs to be treated carefully as \eqref{sol_rp_x1t1} contains a weighted Dirac delta measure across $\Upsilon_2$. However, the equation \eqref{wkfm_1} does not contain $w$, so verifying \eqref{wkfm_1} is standard and can easily be checked. To verify \eqref{wkfm_2}, we insert \eqref{sol_rp_x1t1} on the left-hand side of it and use Green's theorem to yield
\begingroup
\allowdisplaybreaks
\begin{align*}
    &\langle w,\varphi_t \rangle + \langle \frac{w}{h+v}, \varphi_x \rangle - \langle \sigma w, \varphi \rangle\\
    = &\int_{t_{\ast 1}}^\infty\int_{x(t)}^\infty \left(w \varphi_t + \frac{w}{h(t)+v}\varphi_x\right) dx dt + \int_{t_{\ast 1}}^\infty \left(\alpha(t_{\ast 1})h(t)\varphi_t+ \frac{\alpha(t_{\ast 1})h(t)}{h(t)+v_+ h(t)}\varphi_x\right) dt \\
    &- \int_{t_{\ast 1}}^\infty\int_{x(t)}^\infty \sigma(t) w\varphi \;dx dt - \int_{t_{\ast 1}}^\infty \sigma(t)\alpha(t_{\ast 1})h(t)\varphi \; dt\\
    =&\int_{t_{\ast 1}}^\infty\int_{x(t)}^{x_\delta(t)} \left(\frac{v_+ w_-}{v_-}h(t) \varphi_t + \frac{\frac{v_+ w_-}{v_-}h(t)}{h(t)+v_+ h(t)}\varphi_x\right) dx dt +\int_{t_{\ast 1}}^\infty\int_{x_\delta(t)}^\infty \left(w_+ h(t) \varphi_t + \frac{w_+ h(t)}{h(t)+v_+ h(t)}\varphi_x\right) dx dt  \\ 
            &+\int_{t_{\ast 1}}^\infty \alpha(t_{\ast 1})h(t)\frac{d\varphi}{dt} dt - \int_{t_{\ast 1}}^\infty\int_{x(t)}^{x_\delta(t)} \sigma(t)\frac{v_+ w_-}{v_-}h(t)\varphi \;dx dt \\
            &- \int_{t_{\ast 1}}^\infty\int_{x_\delta(t)}^\infty \sigma(t)w_+ h(t)\varphi \;dx dt - \int_{t_{\ast 1}}^\infty \sigma(t)\alpha(t_{\ast 1})h(t)\varphi \; dt\\
            =&\oint -\frac{v_+ w_-}{v_-}h(t) \varphi \;dx + \frac{v_+ w_-}{v_-(1+v_+)}\varphi \;dt - \oint -w_+ h(t)\varphi \;dx + \frac{w_+}{1+v_+}\varphi \;dt - \int_{t_{\ast 1}}^\infty \alpha(t_{\ast 1})\frac{d h(t)}{dt} \varphi dt\\
            &- \int_{t_{\ast 1}}^\infty \sigma(t)\alpha(t_{\ast 1})h(t)\varphi \; dt \\
            =& \int_{t_{\ast 1}}^\infty \left(-\alpha(t_{\ast 1})\frac{d h(t)}{dt}-\sigma(t)\alpha(t_{\ast 1})h(t)+\Big(w_+- \frac{v_+ w_-}{v_-}\Big)\frac{1}{1+v_+}+\Big(\frac{v_+ w_-}{v_-(1+v_+)}-\frac{w_+}{1+v_+}\Big)\right)\varphi \; dt\\
            =& 0,
    \end{align*}
    \endgroup
holds in a local small neighbourhood of $\Upsilon_2$ and for all $\varphi\in C_0^\infty(\mathbb{R}\times \mathbb{R}_+)$ such that $\supp \varphi \cap\{(x,t): x=x_\delta(t), t> t_{\ast 1}\}\neq \emptyset$. Therefore, \eqref{sol_rp_x1t1} satisfies the weak formulation \eqref{weak_form} for any $\varphi$. Hence, \eqref{sol_rp_x1t1} is a weak solution of \eqref{sys} and \eqref{rp_x1,t1}.

Note that the solution \eqref{sol_rp_x1t1} containing a weighted delta measure supported on a contact discontinuity curve $x=x_\delta(t)$, therefore, we call this solution the delta contact discontinuity (see \cite{M_Nedeljkov_08, C_Shen_10_NON}), and denote by $\delta J_2$. Thus, for $t> t_{\ast 1}$, the delta shock wave $\delta S$ splits into a shock wave, namely $S_1$ and a delta contact discontinuity $\delta J_2$ separated by the intermediate state $(V_{\ast 2}h(t), W_{\ast 2}h(t))=(v_+ h(t), \frac{v_+ w_-}{v_-}h(t))$. Furthermore, the delta contact discontinuity $\delta J_2$ propagates on the right from the point $(x_{\ast 1}, t_{\ast 1})$ with speed $\sigma_\delta=\frac{1}{(1+v_+)h(t)}$, and will not interact with the rarefaction wave, where the strength of $\delta J_2$ is $\alpha(t_{\ast 1})h(t)$. The shock curve $S_1$ begins to penetrate the rarefaction wave, and the state $(V, W)$ on it becomes a variable state during the penetration. Then, in the region of penetration, the shock curve with a variable state $(V, W)$ is given by
\begin{equation*}\label{S1_t>t1}
    \begin{cases}
        \frac{dx}{dt}= \frac{1}{(1+v_+)(1+V)h(t)},\\
        \frac{W}{V}=\frac{w_-}{v_-},\\
        x +\epsilon= \frac{1}{(1+V)^2} \int_0^{t} \frac{1}{h(s)}ds,\\
        x(t_{\ast 1})=x_{\ast 1},\; 0\le V\le v_-.
    \end{cases}
\end{equation*}

Now, analogous to Case 4, we have two subcases depending on the values of $v_-, v_+$.

\bigskip

\noindent \textbf{Subcase 5.1.} $v_- < v_+$.

Here, the shock curve $S_1$ penetrates the whole rarefaction wave $R_1^{(-)}$, and interact the wave back (tail) of $R_1^{(-)}$ at a point $(x_{\ast 2}, t_{\ast 2})$ (see Figure \ref{Fig_case5.1}), which can be calculated from
\begin{equation*}
    \begin{aligned}
        &x_{\ast 2} +\epsilon= \frac{1}{(1+v_-)^2} \int_0^{t_{\ast 2}} \frac{1}{h(s)}ds,\\
        & \int_0^{t_{\ast 2}}\frac{1}{h(s)}ds= \frac{(1+v_-)^2(v_+)^2}{(v_+-v_-)^2}\int_0^{t_{\ast 1}}\frac{1}{h(s)}ds.
    \end{aligned}
\end{equation*}
After the penetration, a new shock wave is produced, namely $S_1^{(\ast 2)}$ propagating between the states $(v_-h(t), w_-h(t))$ and $(V_{\ast 2}h(t), W_{\ast 2}h(t))$ with speed $\frac{1}{(1+v_+)(1+v_-)h(t)}$. For sufficiently large $t>t_{\ast 2}$, the solution of initial value problem \eqref{sys} and \eqref{PRP} is as follows 
\begin{equation}\label{sol_case5.1}
    (v_- h(t), w_- h(t))+ S_1^{(\ast 2)}+(V_{\ast 2}h(t), W_{\ast 2}h(t))+ \delta J_2 + (v_+ h(t), w_+ h(t)).
\end{equation}


\begin{figure}[ht!]  
    \centering
\includegraphics[width=0.75\linewidth]{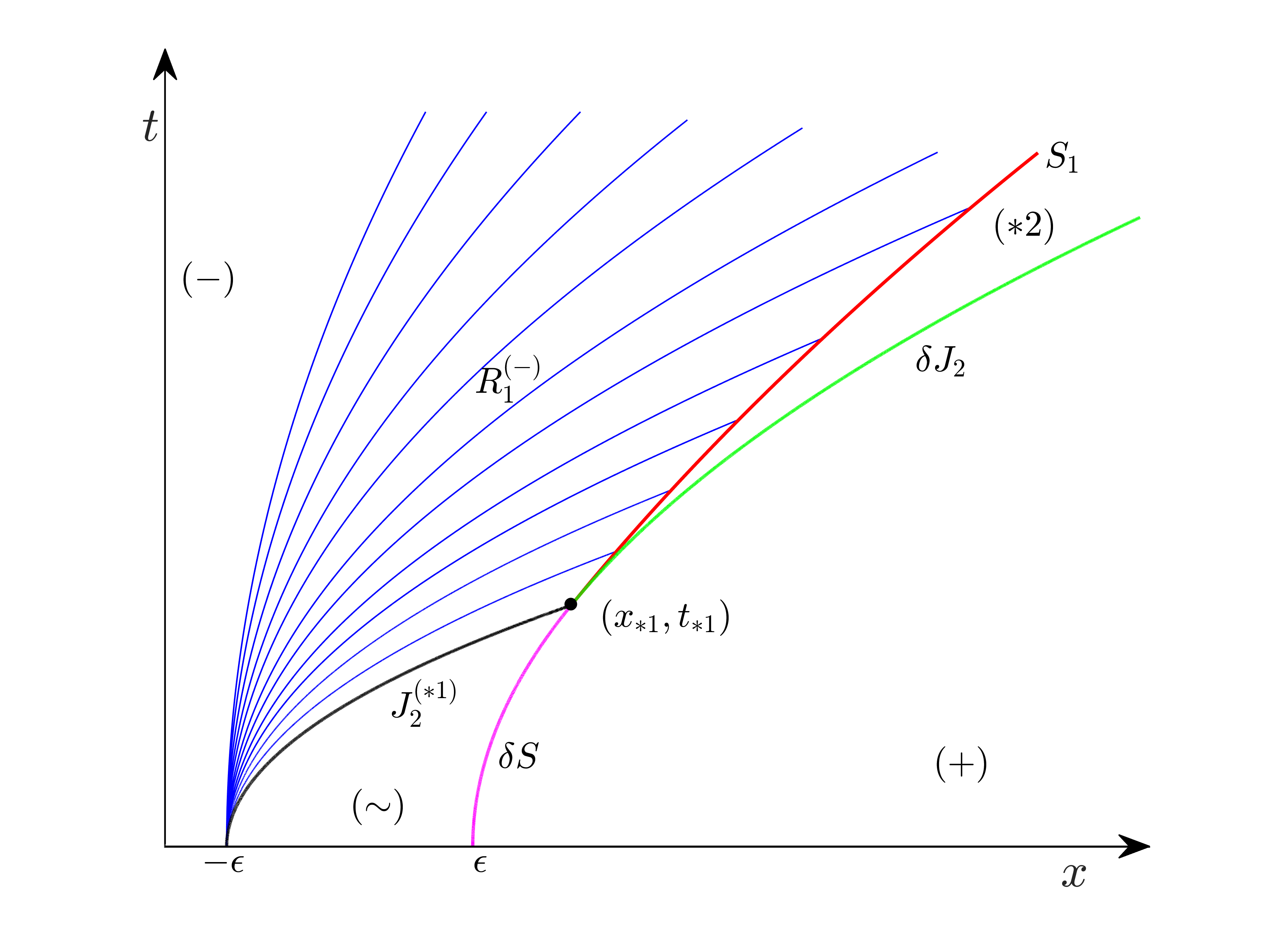}
\caption{Interactions of waves for the subcase 5.2, i.e., when $0= v_\thicksim< v_+<v_-$.}
\label{Fig_case5.2}
\end{figure}

\noindent \textbf{Subcase 5.2.} $v_- > v_+$.

In this subcase, the shock curve $S_1$ cannot fully penetrate the rarefaction wave; see Figure \ref{Fig_case5.2}, and the curve $x +\epsilon= \frac{1}{(1+v_+)^2} \int_0^{t} \frac{1}{h(s)}ds$ becomes its asymptote. Thus, for sufficiently large $t$, the solution of \eqref{sys} and \eqref{PRP} is of the form
\begin{equation}\label{sol_case5.2}
    (v_- h(t), w_- h(t))+ R_1^{(-)}+(V_{\ast 2}h(t), W_{\ast 2}h(t))+ \delta J_2 + (v_+ h(t), w_+ h(t)).
\end{equation}

As $\epsilon\to 0$, it is obvious that $\alpha(t_{\ast 1})\to 0$, and therefore, the delta contact discontinuity $\delta J_2$ converges to a contact discontinuity, namely $J_2$. Thus, as $\epsilon\to 0$, we can conclude that the solution \eqref{sol_case5.1} of the Cauchy problem \eqref{sys} and \eqref{PRP} converges to $S_1+ J_2$ when $v_-< v_+$ whereas the the solution \eqref{sol_case5.2} converges to $R_1+ J_2$ when $v_-> v_+$. Hence, the solution of the Riemann problem \eqref{sys} and \eqref{RP} is stable in this case as well.

\subsection{\textbf{Case 6:} $0=v_-< v_\thicksim< v_+$.}

\begin{figure}[ht!]  
    \centering
\includegraphics[width=0.75\linewidth]{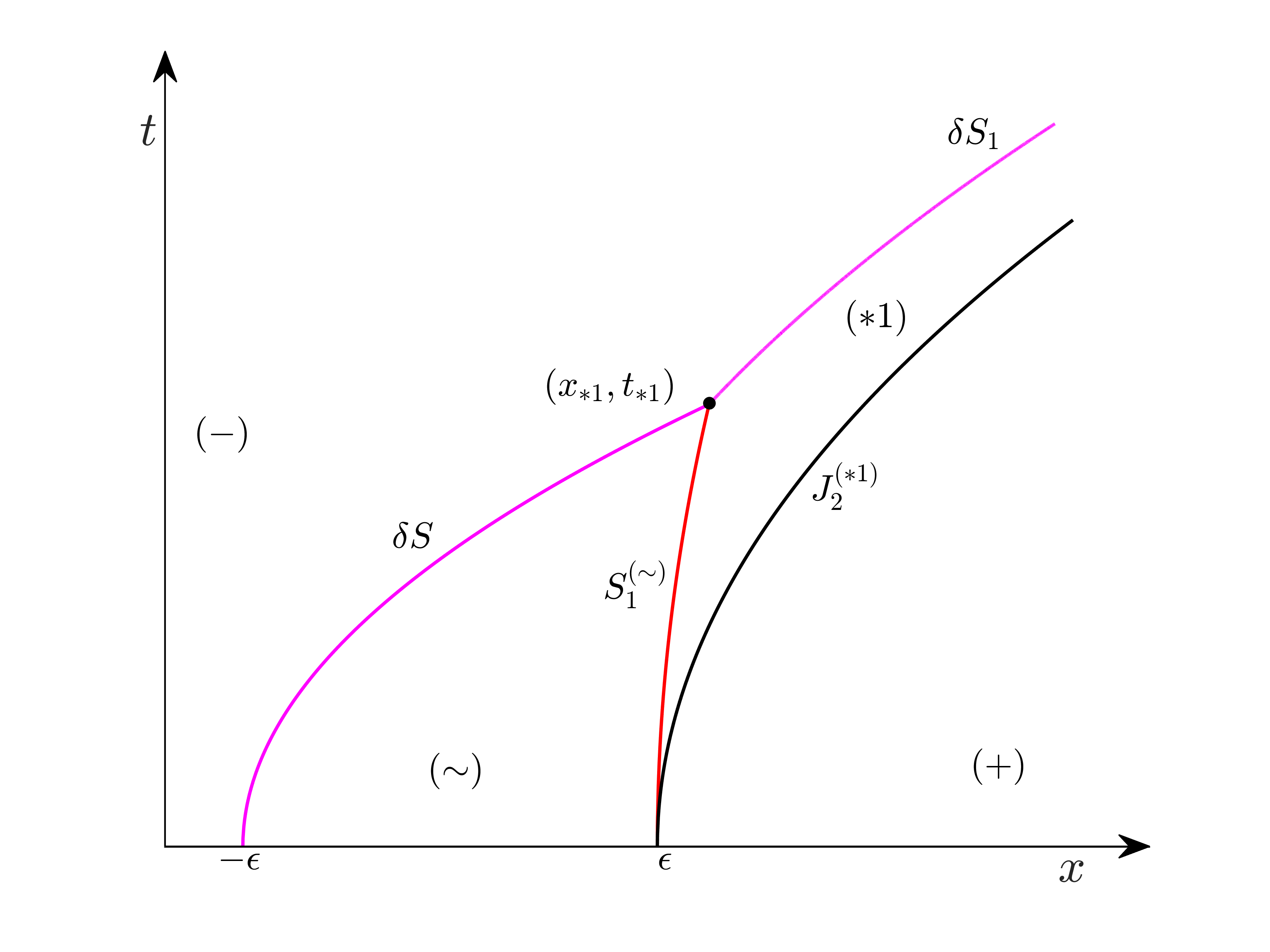}
\caption{Interactions of waves for the Case 6, i.e., when $0=v_-< v_\thicksim< v_+$.}
\label{Fig_case6}
\end{figure}

Since $v_-=0$, then the solution of the local Riemann problem at $(-\epsilon, 0)$ is connected by a delta shock wave $\delta S$, and is given by
\begin{equation*}
    (0, h(t) w_-)+ \delta S + (h(t) v_\thicksim, h(t) w_\thicksim),
\end{equation*}
where the delta shock wave $\delta S$ is  supported on the curve $x+\epsilon=\theta \int_0^t \frac{1}{h(s)} ds$, on which the state $(v(x,t),w(x,t))= (v_\delta(t), \alpha(t)\delta(x + \epsilon-\theta \int_0^t \frac{1}{h(s)} ds))$, with the strength $\alpha(t)=h(t)\frac{w_- v_\thicksim}{1+v_\thicksim}\int_0^t \frac{1}{h(s)} ds$, $\theta=\frac{1}{1+\frac{1}{h(t)}v_\delta}$, $v_\delta(t)=h(t)v_\thicksim$. Moreover, since $0<v_\thicksim< v_+$, then the solution to the local Riemann problem at $(\epsilon, 0)$ is connected by 1-shock wave $S_1^{(\thicksim)}$ and 2-contact discontinuity $J_2^{(\ast 1)}$ starting from $(\epsilon,0)$ as follows
\begin{equation*}
    (v_\thicksim h(t), w_\thicksim h(t)) + S_1^{(\thicksim)} + (V_{\ast 1}h(t), W_{\ast 1} h(t))+ J_2^{(\ast 1)}+ (v_+h(t), w_+ h(t)), 
\end{equation*}
where $(V_{\ast 1}, W_{\ast 1} )= (v_+, \frac{v_+ w_\thicksim}{v_\thicksim})$. Therefore, for sufficiently small $t>0$, the solution of Cauchy problem \eqref{sys} and \eqref{PRP} is of the form (see Figure \ref{Fig_case6})
\begin{equation*}
    (0, h(t) w_-)+ \delta S + (h(t) v_\thicksim, h(t) w_\thicksim) + S_1^{(\thicksim)} + (V_{\ast 1}h(t), W_{\ast 1} h(t))+ J_2^{(\ast 1)}+ (v_+h(t), w_+ h(t)).
\end{equation*}
 
The speeds of propagation for $\delta S$ and $S_1^{(\thicksim)}$ are respectively given by $\sigma_\delta= \frac{1}{(1+v_\thicksim)h(t)}$ and $s_1^{(\thicksim)}=\frac{1}{(1+v_\thicksim)(1+V_{\ast 1})h(t)}$. Then, clearly $\sigma_\delta > s_1^{(\thicksim)}$, and thus $\delta S$ interacts with $S_1^{(\thicksim)}$ at a point, say $(x_{\ast 1}, t_{\ast 1})$, where $t_{\ast 1}>0$ is a finite time and the interaction point $(x_{\ast 1}, t_{\ast 1})$ is calculated from
\begin{equation}\label{x1,t1_case6}
    \begin{aligned}
        &x_{\ast 1}+\epsilon= \frac{1}{(1+v_\thicksim)} \int_0^{t_{\ast 1}} \frac{1}{h(s)}ds,\\
        & x_{\ast 1} -\epsilon= \frac{1}{(1+v_\thicksim)(1+v_+)} \int_0^{t_{\ast 1}} \frac{1}{h(s)}ds.
    \end{aligned}
\end{equation}
Moreover, using \eqref{x1,t1_case6}, we derive the strength of the delta shock wave $\delta S$ at $(x_{\ast 1}, t_{\ast 1})$ is as follows $\alpha(t_{\ast 1}) =h(t_{\ast 1})\frac{w_- v_\thicksim}{1+v_\thicksim}\int_0^{t_{\ast 1}} \frac{1}{h(s)} ds= 2\epsilon w_- v_\thicksim \frac{(1+v_+)}{v_+} h(t_{\ast 1})$.

Now, at the point of interaction $(x_{\ast 1}, t_{\ast 1})$, we have a new local Riemann problem consisting of a weighted delta measure as follows
\begin{equation} \label{rp_x1,t1_case6}
    v\big\lvert_{t=t_{\ast 1}}= \Bigg\{\begin{array}{lr}
        v_- h(t), \quad &x< x_{\ast 1},\\
        V_{\ast 1} h(t), \quad &x> x_{\ast 1},
    \end{array},\quad
    w\big\lvert_{t=t_{\ast 1}}= \Bigg\{\begin{array}{lr}
        w_- h(t), \quad &x< x_{\ast 1},\\
        W_{\ast 1} h(t), \quad &x> x_{\ast 1},
    \end{array}\Bigg\} + \alpha(t_{\ast 1})\delta_{(x_{\ast 1}, t_{\ast 1})}.
\end{equation}
Since, $v_-=0$, then the Riemann problem \eqref{sys} and \eqref{rp_x1,t1_case6} is solved by a new delta shock wave, namely $\delta S_1$, which is of the form
\begin{equation}\label{sol_x1t1_case6}
v(x,t)= \Bigg\{\begin{array}{lr}
        v_- h(t), \quad &x< x_{\delta_1}(t),\\
        V_{\ast 1} h(t), \quad &x> x_{\delta_1}(t),
    \end{array},\quad
    w(x,t)= \Bigg\{\begin{array}{lr}
        w_- h(t), \quad &x<x_{\delta_1}(t),\\
        W_{\ast 1} h(t), \quad &x> x_{\delta_1}(t),
    \end{array}\Bigg\} + \alpha_-(t)D^- + \alpha_+(t)D^+,
\end{equation}
where $\alpha_1(t)D= \alpha_-(t)D^- + \alpha_+(t)D^+$ is a split delta function supported on the curve $x=x_{\delta_1}(t)$ with $x_{\delta_1}(t)=x_{\ast 1}+ \theta_{\delta_1} \int_{t_{\ast 1}}^{t}\frac{1}{h(s)}ds$, $\theta_{\delta_1}=\frac{1}{1+v_+}$, and $\sigma_{\delta_1}= \frac{\theta_{\delta_1}}{h(t)}$ is the propagating speed of $\delta S_1$. Although the delta measures $D^\pm$ supported on the curve $x=x_{\delta_1}(t)$, $D^-$ is the delta measure on $\overline{\mathbb{R}_+^2} \cap\{(x,t):x\le x_{\delta_1}(t)\}$, while $D^+$ is the delta measure on $\overline{\mathbb{R}_+^2} \cap\{(x,t):x\ge x_{\delta_1}(t)\}$.

Moreover, on account of \eqref{sol_x1t1_case6}, we obtain the weak derivatives as follows
\begin{equation}\label{w_t-case6}
\begin{aligned}
    w_t (x,t)=& w_-h^\prime(t) + (W_{\ast 1}h^\prime(t)- w_- h^\prime(t)) H - \frac{\theta_{\delta_1}}{h(t)} (W_{\ast 1}h(t)- w_-h(t)) \delta + (\alpha_-^\prime(t) + \alpha_+^\prime (t))\delta \\&- \frac{\theta_{\delta_1}}{h(t)} (\alpha_-(t) + \alpha_+ (t))\delta^\prime,
    \end{aligned}
\end{equation}
and 
\begin{equation}\label{flux_x-case6}
    \Big(\frac{w}{h(t)+v}\Big)_x (x,t)= \Big(\frac{W_{\ast 1}h(t)}{h(t)+V_{\ast 1}h(t)} - \frac{w_-h(t)}{h(t)+ v_- h(t)}\Big)\delta + \Big(\frac{\alpha_-(t)}{h(t)+ v_-h(t)}+ \frac{\alpha_+(t)}{h(t)+V_{\ast 1}h(t)} \Big) \delta^\prime,
\end{equation}
in the sense of distribution, where $H$ is the Heaviside function supported on $x=x_{\delta_1}(t)$, i.e., the value of $H$ is equal to $0$ on the left-hand side of $x=x_{\delta_1}(t)$, while $1$ on the right-hand side of this curve. Consequently, inserting \eqref{w_t-case6} and \eqref{flux_x-case6} into the second equation of \eqref{sys} yields
\begin{equation*}
    \Big(\frac{w_-}{1+v_+}- \frac{w_-}{1+v_-}+ \alpha_-^\prime(t)+ \alpha_+^\prime(t)+ \sigma(t)\big(\alpha_-(t)+ \alpha_+(t)\big)\Big)\delta + \Big(\frac{\alpha_-(t)}{h(t)+v_-h(t)}- \frac{\alpha_-(t)}{h(t)+v_+h(t)}\Big)\delta^\prime =0,
\end{equation*}
then comparing the coefficients of $\delta$ and $\delta^\prime$, we obtain
\begin{align}
    &\frac{w_-}{1+v_+}- \frac{w_-}{1+v_-}+ \alpha_-^\prime(t)+ \alpha_+^\prime(t)+ \sigma(t)\big(\alpha_-(t)+ \alpha_+(t)\big)=0, \label{coeff_1}\\
    &\frac{\alpha_-(t)}{h(t)}\Big(\frac{1}{1+v_-}- \frac{1}{1+v_+}\Big)=0. \label{coeffi_2}
\end{align}
The equation \eqref{coeffi_2} implies $\alpha_-(t)=0$. On account of \eqref{coeff_1} and $v_-=0$, it follows that
\begin{equation*}
    \alpha_+(t)= 2\epsilon w_- v_\thicksim\Big(\frac{1+v_+}{v_+}\Big) h(t) +  \Big(\frac{w_- v_+}{1+v_+}\Big) h(t)  \int_{t_{\ast 1}}^{t}\frac{1}{h(s)}ds.
\end{equation*}
Thus, the strength of the delta shock wave $\delta S_1$ is given by $\alpha_1(t)= \alpha_-(t)+ \alpha_+(t)= \alpha_+(t)$. Also note that, the propagating speeds of $\delta S_1$ and the contact discontinuity $J_2^{(\ast 1)}$ are same, i.e., $\frac{1}{(1+v_+)h(t)}$. Therefore, they will never interact, and no further interaction happens; refer to Figure \ref{Fig_case6}. Hence, for sufficiently large $t$, the solution of the Cauchy problem is as follows
\begin{equation} \label{sol_case6}
(0, h(t) w_-)+ \delta S_1 + (V_{\ast 1}h(t), W_{\ast 1} h(t))+ J_2^{(\ast 1)}+ (v_+h(t), w_+ h(t)).
\end{equation}

Meanwhile, observe that as $\epsilon\to 0$, $(x_{\ast 1}, t_{\ast 1})\to (0,0)$ and $\alpha(t_{\ast 1}) \to 0$. Moreover, $\lim\limits_{\epsilon\to 0}\alpha_1(t)= h(t) \Big(\frac{w_- v_+}{1+v_+}\Big)  \int_{0}^{t}\frac{1}{h(s)}ds$. Therefore, the solution \eqref{sol_case6} converges to a single delta shock wave with strength $h(t) \Big(\frac{w_- v_+}{1+v_+}\Big)  \int_{0}^{t}\frac{1}{h(s)}ds$, which is exactly the solution of the Riemann problem \eqref{sys} and \eqref{RP} in this case. This suggests that the Riemann solution is stable in this case as well.

\subsection{\textbf{Case 7:} $0=v_-< v_+< v_\thicksim$.}
\begin{figure}[ht!]  
    \centering
\includegraphics[width=0.75\linewidth]{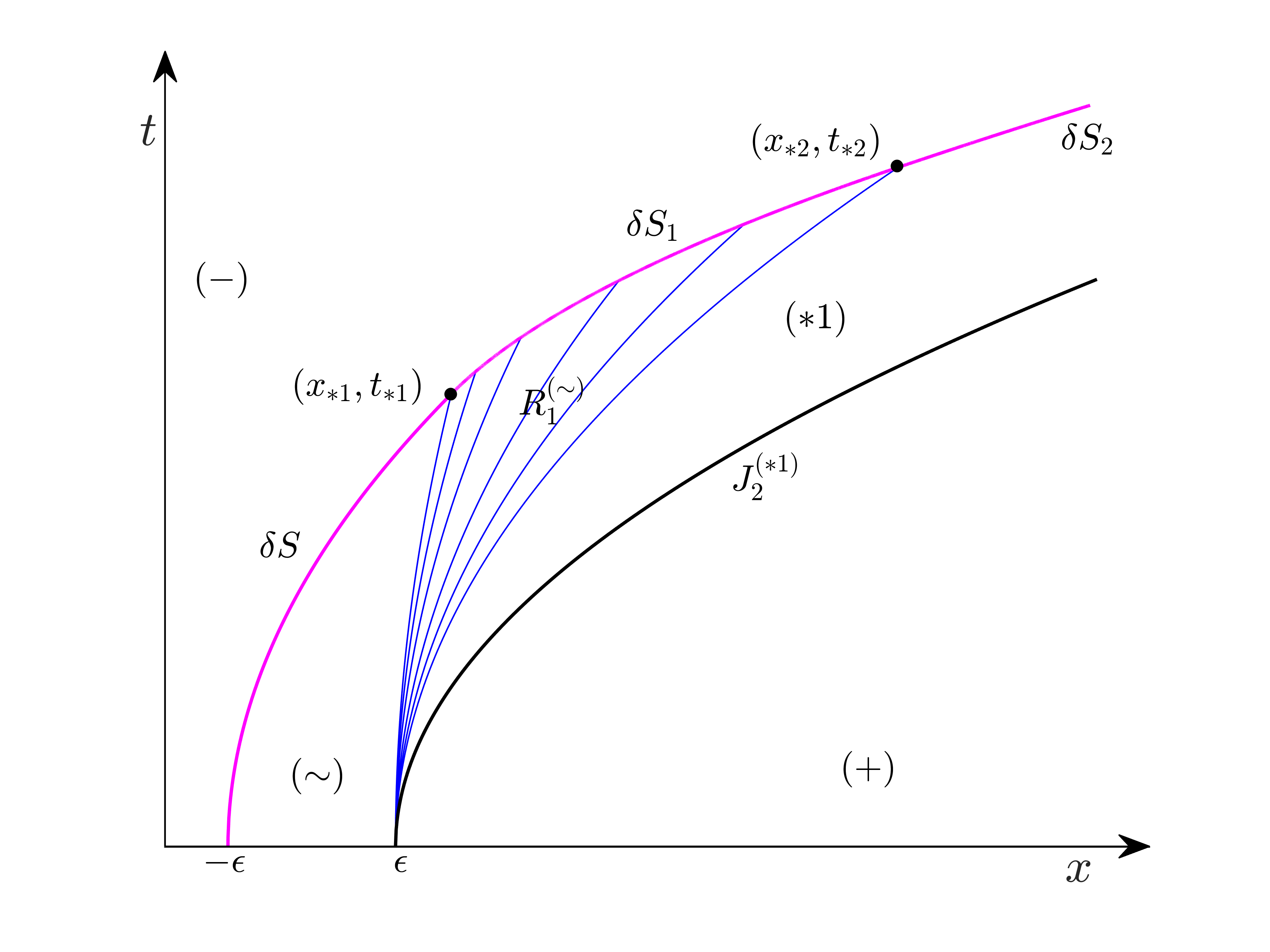}
\caption{Case 7: Interactions of waves when $0=v_-< v_+<v_\thicksim$.}
\label{Fig_case7}
\end{figure}

 In this case, the solution of the Cauchy problem can be expressed as follows 
 \begin{equation}\label{case7}
    (0, h(t) w_-)+ \delta S + (h(t) v_\thicksim, h(t) w_\thicksim) + R_1^{(\thicksim)} + (V_{\ast 1}h(t), W_{\ast 1} h(t))+ J_2^{(\ast 1)}+ (v_+h(t), w_+ h(t)),
\end{equation}
for sufficiently small $t$. In \eqref{case7}, $\delta S$ is the delta shock wave as in the previous case, while $R_1^{(\thicksim)}$ and $J_2^{(\ast 1)}$ are the rarefaction wave and contact discontinuity starting from $(\epsilon,0)$, with the intermediate state $(V_{\ast 1}h(t), W_{\ast 1} h(t))=(v_+, \frac{v_+ w_\thicksim}{v_\thicksim})h(t)$. The propagating speed of $\delta S$ is $\sigma_\delta= \frac{1}{(1+v_\thicksim)h(t)}$, which is larger than speed of wave back of $R_1^{(\thicksim)}$, i.e., $\frac{1}{(1+v_\thicksim)^2 h(t)}$, and therefore, $\delta S$ must catch the wave back in time $t_{\ast 1}>0$; see Figure \ref{Fig_case7}. Then, the point of interaction $(x_{\ast 1}, t_{\ast 1})$ is given by 
\begin{equation*}
    \begin{aligned}
        &x_{\ast 1}+\epsilon= \frac{1}{(1+v_\thicksim)} \int_0^{t_{\ast 1}} \frac{1}{h(s)}ds,\\
        & x_{\ast 1} -\epsilon= \frac{1}{(1+v_\thicksim)^2} \int_0^{t_{\ast 1}} \frac{1}{h(s)}ds.
    \end{aligned}
\end{equation*}
The strength of $\delta S$ at $(x_{\ast 1}, t_{\ast 1})$ is $\alpha(t_{\ast 1}) =h(t_{\ast 1})\frac{w_- v_\thicksim}{1+v_\thicksim}\int_0^{t_{\ast 1}} \frac{1}{h(s)} ds= 2\epsilon w_- (1+v_\thicksim) h(t_{\ast 1})$. From the point of interaction, the delta shock wave $\delta S$ starts penetrating the rarefaction wave $R_1^{(\thicksim)}$ for $t>t_{\ast 1}$, and produces a new delta shock wave during the penetration. We denote the delta shock wave by $\delta S_1$ in the penetration region. Then, the state on the left-hand side of $\delta S_1$ is $(0,h(t)w_-)$, whereas the right-hand state is $(Vh(t), W h(t))$ with $(V, W)$ being variable across the rarefaction fan $R_1^{(\thicksim)}$. Let us denote the delta shock curve for $\delta S_1$ as $\Upsilon_1:\{(x,t): x=x_{\delta_1}(t), t\ge t_{\ast 1}\}$. Then, the curve $\Upsilon_1$ can be determined from 
\begin{equation}\label{deltaS1_t>t1}
    \begin{cases}
       \sigma_{\delta_1}= \frac{dx}{dt}= \frac{1}{(1+V)h(t)},\\
        x -\epsilon= \frac{1}{(1+V)^2} \int_0^{t} \frac{1}{h(s)}ds,\\
        \frac{W}{V}=\frac{w_{\thicksim}}{v_{\thicksim}},\\
        x(t_{\ast 1})=x_{\ast 1},\; v_+\le V\le v_\thicksim.
    \end{cases}
\end{equation}

Then from the second and third equations of \eqref{deltaS1_t>t1}, we obtain the right state $(Vh(t), W h(t))$ across $\delta S_1$ as follows
\begin{equation}\label{right_st_V}
    (V, W)= \bigg(V, \frac{V w_\thicksim}{v_\thicksim}\bigg)= \left(\sqrt{\frac{\int_0^{t} \frac{1}{h(s)}ds}{x-\epsilon}}-1, \frac{w_\thicksim}{v_\thicksim}\Bigg(\sqrt{\frac{\int_0^{t} \frac{1}{h(s)}ds}{x-\epsilon}}-1\Bigg)\right).
\end{equation}
On account of first equation of \eqref{deltaS1_t>t1} and \eqref{right_st_V}, one obtains 
\begin{equation} \label{spd_upsilon1}
    \sigma_{\delta_1}= \frac{dx}{dt}= \frac{1}{h(t)}\sqrt{\frac{x-\epsilon}{\int_0^{t} \frac{1}{h(s)}ds}}, \quad x(t_{\ast 1})=x_{\ast 1},
\end{equation}
yields the delta shock curve $\Upsilon_1$ as follows 
\begin{equation}\label{upsilon1}
    \Upsilon_1: x=x_{\delta_1}(t)\equiv \epsilon + \Bigg(\sqrt{\int_0^t \frac{1}{h(s)}ds} -\sqrt{2\epsilon v_\thicksim}\Bigg)^2.
\end{equation}
Then, from \eqref{spd_upsilon1} and \eqref{upsilon1}, one obtains
\begin{equation*}
    \frac{d^2x}{dt^2}= \frac{\sigma(t)}{h(t)} \sqrt{\frac{x-\epsilon}{\int_0^{t} \frac{1}{h(s)}ds}} + \frac{1}{2(h(t))^2} \frac{\sqrt{2\epsilon v_\thicksim}}{\big(\int_0^{t} \frac{1}{h(s)}ds\big)^{3/2}} >0,
\end{equation*}
which implies the speed of propagation of the delta shock wave $\delta S_1$ increases in the penetration region. The delta shock wave $\delta S_1$ supported on $\Upsilon_1$ can be constructed using a weighted split delta function as follows 
\begin{equation}\label{deltaS1_case7}
v(x,t)= \Bigg\{\begin{array}{lr}
        0, \quad &x< x_{\delta_1}(t),\\
        V h(t), \quad &x> x_{\delta_1}(t),
    \end{array},\quad
    w(x,t)= \Bigg\{\begin{array}{lr}
        w_- h(t), \quad &x<x_{\delta_1}(t),\\
        W h(t), \quad &x> x_{\delta_1}(t),
    \end{array}\Bigg\} + \alpha_1^-(t)D_{\Upsilon_1}^- + \alpha_1^+(t)D_{\Upsilon_1}^+,
\end{equation}
where the split delta function $\alpha_1(t) D_{\Upsilon_1}= \alpha_1^-(t)D_{\Upsilon_1}^- + \alpha_1^+(t)D_{\Upsilon_1}^+$ is supported on $\Upsilon_1$, $\alpha_1(t)=\alpha_1^-(t)+ \alpha_1^+(t)$ represents the strength of $\delta S_1$, to be determined, and $V, W$ are given in \eqref{right_st_V}.
Then, from \eqref{deltaS1_case7}, we obtain in the sense of distribution that
\begin{equation}\label{w_t_case7}
    \begin{aligned}
        &w_t (x,t)= w_-h^\prime(t) + \Bigg(\frac{w_\thicksim}{v_\thicksim}\frac{1}{2\sqrt{(x-\epsilon)\int_0^t \frac{1}{h(s)}ds}}+ \frac{w_\thicksim}{v_\thicksim} \sqrt{\frac{\int_0^{t} \frac{1}{h(s)}ds}{x-\epsilon}} h^\prime(t) - \frac{w_\thicksim}{v_\thicksim} h^\prime(t) - w_- h^\prime(t)\Bigg) H \\&- \sigma_{\delta_1} \Bigg(\frac{w_\thicksim}{v_\thicksim}\Bigg(\sqrt{\frac{\int_0^{t} \frac{1}{h(s)}ds}{x-\epsilon}}-1\Bigg) h(t)- w_-h(t)\Bigg) \delta + \big((\alpha_1^{-})^\prime(t) + (\alpha_1^+)^\prime (t)\big)\delta - \sigma_{\delta_1} (\alpha_1^-(t) + \alpha_1^+ (t))\delta^\prime,
    \end{aligned}
\end{equation}
and 
\begin{equation}\label{flux_x-case7}
\begin{aligned}
    \Big(\frac{w}{h(t)+v}\Big)_x (x,t)= -\frac{w_\thicksim}{2v_\thicksim}\frac{1}{\sqrt{(x-\epsilon)\int_0^t\frac{1}{h(s)}ds}}H+ \Big(\frac{w_\thicksim}{v_\thicksim}\Big(1-\sqrt{\frac{x-\epsilon}{\int_0^t\frac{1}{h(s)}ds}}\Big)-w_-\Big)\delta \\ + \frac{\alpha_1^-(t)}{h(t)}\delta^\prime  + \frac{\alpha_1^+(t)}{h(t)}\sqrt{\frac{x-\epsilon}{\int_0^t\frac{1}{h(s)}ds}} \delta^\prime,
    \end{aligned}
\end{equation}
where $H, \delta, \delta^\prime$ are supported on $x=x_{\delta_1}(t)$, i.e., functions of $x-x_{\delta_1}(t)$. Inserting these \eqref{w_t_case7} and \eqref{flux_x-case7} into the second equation of \eqref{sys} and consequently, comparing the coefficients of $\delta$ and $\delta^\prime$, it follows that 
\begin{align}
    &-\sigma_{\delta_1} \Bigg(\frac{w_\thicksim}{v_\thicksim}\Bigg(\sqrt{\frac{\int_0^{t} \frac{1}{h(s)}ds}{x-\epsilon}}-1\Bigg) - w_-\Bigg)h(t) +(\alpha_1^-)^\prime(t) + (\alpha_1^+)^\prime(t)+ \Big(\frac{w_\thicksim}{v_\thicksim}\Big(1-\sqrt{\frac{x-\epsilon}{\int_0^t\frac{1}{h(s)}ds}}\Big)-w_-\Big)\nonumber\\& +(\alpha_1^-(t) +\alpha_1^+(t))\sigma(t)=0 , \label{coeff_1_case7}\\
   &\text{and}\; -\sigma_{\delta_1}(\alpha_1^-(t) +\alpha_1^+(t)) + \frac{\alpha_1^-(t)}{h(t)}+ \frac{\alpha_1^+(t)}{h(t)} \sqrt{\frac{x-\epsilon}{\int_0^t\frac{1}{h(s)}ds}} =0. \label{coeffi_2_case7}
\end{align}
Using \eqref{spd_upsilon1} and \eqref{coeffi_2_case7}, we obtain 
\begin{equation*}
    \frac{\alpha_1^-(t)}{h(t)} \Big(1-\sqrt{\frac{x-\epsilon}{\int_0^t\frac{1}{h(s)}ds}}\Big)=0,
\end{equation*}
which on account of \eqref{upsilon1} yields 
\begin{equation*}
    \frac{\alpha_1^-(t)}{h(t)} \sqrt{\frac{2\epsilon v_\thicksim}{\int_0^t\frac{1}{h(s)}ds}}=0,
\end{equation*}
therefore $\alpha_1^-(t)=0$. Thus, inserting \eqref{spd_upsilon1} and \eqref{upsilon1} into \eqref{coeff_1_case7} leads to
\begin{equation*}
    (\alpha_1^+)^\prime(t)+ \sigma(t) (\alpha_1^+)(t) = w_- \sqrt{\frac{2\epsilon v_\thicksim}{\int_0^t\frac{1}{h(s)}ds}},
\end{equation*}
which on integration over $[t_{\ast_1}, t]$ with $\alpha_1^+(t_{\ast_1})=\alpha_1(t_{\ast_1})=\alpha(t_{\ast_1})$ yields
\begin{equation*}
    \alpha_1^+(t)= 2w_- \sqrt{2\epsilon v_\thicksim} h(t)\sqrt{\int_0^t\frac{1}{h(s)}ds} - 2\epsilon w_-(1+v_\thicksim) h(t).
\end{equation*}
Hence, $\alpha_1(t)= \alpha_1^1(t)+\alpha_1^+(t)=\alpha_1^+(t)$ is the strength of the delta shock $\delta S_1$ during the penetration for $t> t_{\ast_1}$.

Now, when $V\to v_+$, the delta shock $\delta S_1$ overtakes $R_1^{(\thicksim)}$ at a point $(x_{\ast 2}, t_{\ast 2})$, with $t_{\ast 2}>t_{\ast 1}$. Then, we have $\sqrt{\frac{\int_0^{t_{\ast 2}} \frac{1}{h(s)}ds}{x_{\ast 2}-\epsilon}}-1= v_+$, and therefore, from \eqref{upsilon1}, it follows that $t_{\ast 2}$ can be determined by 
\begin{equation*}
    \int_0^{t_{\ast 2}} \frac{1}{h(s)}ds = 2\epsilon v_\thicksim \Big(\frac{1+v_+}{v_+}\Big)^2.
\end{equation*}
Thus, we obtain $x_{\ast 2}$ as follows $x_{\ast 2}= 2\epsilon \frac{v_\thicksim}{v_+^2}+\epsilon$. So, we determined the interaction point $(x_{\ast 2}, t_{\ast 2})$, where $\delta S_1$ fully penetrate the rarefaction wave. Consequently, it produces a new local Riemann problem at $(x_{\ast 2}, t_{\ast 2})$ as follows
\begin{equation*} 
    v\big\lvert_{t=t_{\ast 2}}= \Bigg\{\begin{array}{lr}
        0, \quad &x< x_{\ast 2},\\
        v_+ h(t), \quad &x> x_{\ast 2},
    \end{array},\quad
    w\big\lvert_{t=t_{\ast 2}}= \Bigg\{\begin{array}{lr}
        w_- h(t), \quad &x< x_{\ast 2},\\
        \frac{v_+ w_\thicksim}{v_\thicksim} h(t), \quad &x> x_{\ast 2},
    \end{array}\Bigg\} + \alpha_1(t_{\ast 2})\delta_{(x_{\ast 2}, t_{\ast 2})}.
\end{equation*}
This must be solved by a delta shock wave, denote it by $\delta S_2$. Moreover, the speed and the strength of $\delta S_2$ are, respectively, as follows
\begin{equation*}
    \sigma_{\delta_2}= \frac{1}{(1+v_+)h(t)}, \quad \alpha_2(t)= \frac{w_- v_+}{1+v_+} h(t) \int_{t_{\ast 2}}^t \frac{1}{h(s)}ds + \frac{\alpha_1(t_{\ast 2})}{h(t_{\ast 2})} h(t),
\end{equation*}
where $t>t_{\ast 2}$. The propagating speeds of $\delta S_2$ and $J_2^{(\ast 1)}$ coincide, so they never interact with each other; see Figure \ref{Fig_case7}. Hence, for $t>t_{\ast 2}$, the solution of the Cauchy problem \eqref{sys} and \eqref{PRP} is of the form as follows
\begin{equation} \label{sol_case7}
(0, h(t) w_-)+ \delta S_2 + (V_{\ast 1}h(t), W_{\ast 1} h(t))+ J_2^{(\ast 1)}+ (v_+h(t), w_+ h(t)).
\end{equation}

Now, as $\epsilon\to 0$, the solution \eqref{sol_case7} converges to a delta shock wave solution propagating from $(0,0)$ with the strength $\frac{w_- v_+}{1+v_+} h(t) \int_0^t \frac{1}{h(s)}ds$ of the corresponding Riemann problem \eqref{sys} and \eqref{RP}. In conclusion, the Riemann solution is stable in this case as well under small perturbations of the initial data.

\begin{remark}
    Note that the interaction between two delta shock waves, which might arise from the local Riemann problems, is not possible due to the condition of appearance of the delta shock and the nonnegativity of the state variable $v$. Moreover, there will be no further cases of interaction, as we discussed all possible choices of initial data \eqref{PRP}. Thus, combining all cases 1-7, we conclude that the proof of the Theorem \ref{main_thm} is complete.
\end{remark}


\section{Numerical evidences} \label{sec:numerics}

We discretize \eqref{sys} over the space-time domain $\mathbb{R} \times [0,T]$, with $T>0$, using a uniform spatial mesh size $\Delta x >0$ and a variable time step $\Delta t^n > 0$ satisfying an adaptive CFL condition. Let $x_i = i \Delta x$ and $x_{i\pm 1/2} = (i \pm \frac{1}{2}) \Delta x$ for $i \in \mathbb{Z}$, and define the discrete times recursively by $t^0 = 0$, $t^{n+1} = t^n + \Delta t^n$, $n \in \mathbb{N}$. 
The cell averages $v_i^n$ and $w_i^n$ approximate the mean values of $v$ and $w$ over the cell $(x_{i-1/2}, x_{i+1/2})$ at time $t^n$. 
We employ the Lax--Friedrichs type scheme with averaged source term \cite{Rcruz_RM_WN}:
\begin{equation*} \label{eq:LF}
\begin{aligned}
v_i^{n+1} &= \frac{v_{i-1}^n + v_{i+1}^n}{2} - \frac{\Delta t^n}{2 \Delta x} \left( \frac{v_{i+1}^n}{h^n + v_{i+1}^n} - \frac{v_{i-1}^n}{h^n + v_{i-1}^n} \right) - \frac{\Delta t^n}{2} \sigma(t^n) (v_{i-1}^n + v_{i+1}^n),\\
w_i^{n+1} &= \frac{w_{i-1}^n + w_{i+1}^n}{2} - \frac{\Delta t^n}{2 \Delta x} \left( \frac{w_{i+1}^n}{h^n + v_{i+1}^n} - \frac{w_{i-1}^n}{h^n + v_{i-1}^n} \right) - \frac{\Delta t^n}{2} \sigma(t^n) (w_{i-1}^n + w_{i+1}^n),
\end{aligned}
\end{equation*}
where $h^{n+1} = h^n \exp (-\sigma(t^n) \Delta t^n)$, so that $h^n >0$ for all $n \in \mathbb{N}$. 
The adaptive time step is given by
\begin{equation*}
\Delta t^n = \frac{\mathrm{CFL}}{\sigma^n + \frac{1}{\Delta x} \max\limits_{i} \left(\frac{1}{h^n + v_i^n}\right)},
\end{equation*}
with $0 < \mathrm{CFL} < 1$. We denote by $\Lambda(t)$ the integral $\int_0^t \frac{ds}{h(s)} = \int_0^t \exp \left(\frac{0.088\,s}{1+s} \right)ds$, and by $v^{\Delta x}(\cdot, t^n)$ and $w^{\Delta x}(\cdot, t^n)$ the piecewise constant functions reconstructed from the cell averages $v^n_i$ and $w^n_i$.

For all experiments in this section, following \cite{Rcruz_RM_WN}, we choose the time-gradually-degenerate damping coefficient $\sigma(t) = \frac{0.0880}{(1+t)^2}$, so that $h(t) = \exp\!\left(-0.0880\,\frac{t}{1+t}\right)$. 
This choice corresponds to the bounded growth regime identified in \cite{Rcruz_RM_WN}, in which $\int_0^\infty \sigma(s)\,ds = 0.0880 < \infty$ and both $n(t)$ and $h(t)$ stabilize at finite positive limits
\begin{equation*}
n(t) \xrightarrow{t\to\infty} n_\infty := n_0\,e^{0.0880} < \infty, 
\qquad h(t) \xrightarrow{t\to\infty} h_\infty := e^{-0.0880} 
\approx 0.9158 > 0.
\end{equation*}
This is the physically most relevant regime for chromatography, since the saturation loading capacity $n(t)$ remains strictly bounded for all $t \geq 0$, the flux denominator $h(t) + v$ does not degenerate, and the system retains the full competitive Langmuir structure asymptotically \cite{Rcruz_RM_WN}. 
Furthermore, the total growth $n_\infty/n_0 = e^{0.0880} \approx 1.092$, representing a $9.2\%$ increase in adsorption capacity, falls within the experimentally reported range of $5\%$ to $10\%$ moderate variations of the saturation capacity under controlled operating conditions \cite{Rcruz_RM_WN}.

The profiles displayed in the figures are obtained with mesh size $\Delta x = 0.001$ and $\mathrm{CFL} = 0.9$.

\subsection{Space-time evolution of the numerical solutions}

To validate the theoretical results established in Section~\ref{sec:Wint}, we present numerical solutions of the Cauchy problem \eqref{sys} and \eqref{PRP} obtained via the Lax-Friedrichs type scheme with CFL $= 0.9$ and $\sigma(t) = \frac{0.0880}{(1+t)^2}$, so that $h(t) = \exp \left( -0.0880 \frac{t}{1+t} \right)$.
The solutions are displayed as space-time heat maps, where the color intensity represents the value of $v(x,t)$ and $w(x,t)$ at each point $(x,t)$, allowing the wave interactions to be visualized as characteristic-like curves in the $(x,t)$ plane.
For each case we take $\epsilon = 1/2$, so the initial discontinuities are placed at $x = \pm 1/2$.
The initial data for all cases are summarized in different tables, except to the first case where we show each detail in the interaction. \\

{\bf Case 1: $0<v_-<v_\thicksim<v_+$.}

We consider the initial data
\begin{equation} \label{RPC1}
    (v_0(x),w_0(x))=
    \begin{cases}
        (1,2), &\mbox{if } x < -\frac{1}{2},\\
        (2,3), &\mbox{if } -\frac{1}{2} < x < \frac{1}{2},\\
        (4,1), &\mbox{if } x > \frac{1}{2},
    \end{cases}
\end{equation}
which satisfies $0 < v_- = 1 < v_\thicksim = 2 < v_+ = 4$, corresponding to Case~1 of Section~\ref{Sect3.1}. From Section~\ref{Sect3.1}, there are two interaction times $t_{*1}$ and $t_{*2}$.
For $0 \leq t < t_{*1}$, the solution of \eqref{sys} and \eqref{RPC1} is
\begin{equation*}
\exp \left( -0.0880 \frac{t}{1+t} \right) \times
\begin{cases}
(1,2), &\mbox{if } x < g_1(t;1,2),\\
(2,4), &\mbox{if } g_1(t;1,2) < x < g_2(t;2),\\
(2,3), &\mbox{if } g_2(t;2) < x < g_1(t;2,4),\\
(4,6), &\mbox{if } g_1(t;2,4) < x < g_2(t;4),\\
(4,1), &\mbox{if } x > g_2(t;4),
\end{cases}
\end{equation*}
where
\begin{align*}
g_1(t;1,2) &= -\frac{1}{2} + \frac{1}{6} \int_0^t \exp \left( 0.0880 
\frac{s}{1+s} \right) ds, \\
g_2(t;2) &= -\frac{1}{2} + \frac{1}{3} \int_0^t \exp \left( 0.0880 
\frac{s}{1+s} \right) ds, \\
g_1(t;2,4) &= \frac{1}{2} + \frac{1}{15} \int_0^t \exp \left( 0.0880 
\frac{s}{1+s} \right) ds,
\end{align*}
and
$$
g_2(t;4) = \frac{1}{2} + \frac{1}{5} \int_0^t \exp \left( 0.0880 
\frac{s}{1+s} \right) ds.
$$
The first interaction time $t_{*1}$ is determined by the condition $g_2(t_{*1};2) = g_1(t_{*1};2,4)$, which gives
\begin{equation*}
\int_0^{t_{*1}} \exp \left( 0.0880 \frac{s}{1+s} \right) ds = \frac{15}{4},
\end{equation*}
and by the Newton--Raphson method we obtain $t_{*1} \approx 3.5647$ and $x_{*1} = \frac{3}{4}$.
At $(x_{*1}, t_{*1})$, the contact discontinuity $J_2^{(*1)}$ meets the shock $S_1^{(\sim)}$, producing a new local Riemann problem with left state $(V_{*1}, W_{*1}) = (2, 4)$ and right state $(V_{*2}, W_{*2}) = (4, 6)$. 
Since $V_{*1} = 2 < V_{*2} = 4$, this is solved by a new shock $S_1^{(*1)}$ and a new contact discontinuity $J_2^{(*3)}$, with intermediate state $(V_{*3}, W_{*3}) = (4, 8)$.
Therefore, for $t_{*1} \leq t < t_{*2}$, the solution is
\begin{equation*}
\exp \left( -0.0880 \frac{t}{1+t} \right) \times
\begin{cases}
(1,2), &\mbox{if } x < g_1(t;1,2),\\
(2,4), &\mbox{if } g_1(t;1,2) < x < g_1^*(t;2,4),\\
(4,8), &\mbox{if } g_1^*(t;2,4) < x < g_2^*(t;4),\\
(4,6), &\mbox{if } g_2^*(t;4) < x < g_2(t;4),\\
(4,1), &\mbox{if } x > g_2(t;4),
\end{cases}
\end{equation*}
where
\begin{align*}
g_1^*(t;2,4) &= \frac{3}{4} + \frac{1}{15} \int_{t_{*1}}^t \exp \left( 0.0880 \frac{s}{1+s} \right) ds 
\mbox{ and } 
g_2^*(t;4) &= \frac{3}{4} + \frac{1}{5} \int_{t_{*1}}^t \exp \left( 0.0880 \frac{s}{1+s} \right) ds.
\end{align*}
Since $V_{*3} = V_{*2} = 4$, the contact discontinuities $J_2^{(*3)}$ and $J_2^{(*2)}$, given respectively by $g_2^*(t;4)$ and $g_2(t;4)$, satisfy $\frac{dg_2^*}{dt} = \frac{dg_2}{dt} = \frac{1}{5h(t)}$, for all $t > t_{*1}$, so they are parallel and will never interact. The second interaction time $t_{*2}$ is determined by the condition $g_1(t_{*2};1,2) = g_1^*(t_{*2};2,4)$, which gives
\begin{equation*}
\int_{t_{*1}}^{t_{*2}} \exp \left( 0.0880 \frac{s}{1+s} \right) ds = \frac{25}{4},
\end{equation*}
and by the Newton--Raphson method we obtain $t_{*2} \approx 9.3599$ and $x_{*2} = \frac{7}{6}$.

At $(x_{*2}, t_{*2})$, the shock $S_1^{(-)}$ meets the shock $S_1^{(*1)}$, producing a new local Riemann problem with left state $(v_-, w_-) = (1, 2)$ and right state $(V_{*3}, W_{*3}) = (4, 8)$. 
Since $v_- = 1 < V_{*3} = 4$, the solution structure is $S_1^{(*3)} + J_2$. 
The intermediate state between $S_1^{(*3)}$ and $J_2$ is 
\begin{equation*}
\left(V_{*3},\; V_{*3} \cdot \frac{w_-}{v_-}\right) = \left(4,\; 4 \cdot \frac{2}{1}\right) = (4, 8) = (V_{*3}, W_{*3}),
\end{equation*}
which coincides with the right state. Therefore the contact discontinuity $J_2$ is trivial and the solution consists of the single shock $S_1^{(*3)}$ with speed $\dfrac{dg_1^{**}}{dt} = \dfrac{1}{(1+v_-)(1+V_{*3})h(t)} = \dfrac{1}{10\,h(t)}$.

Therefore, for $t \geq t_{*2}$, the solution is
\begin{equation*}
\exp \left( -0.0880 \frac{t}{1+t} \right) \times
\begin{cases}
(1,2), &\mbox{if } x < g_1^{**}(t;1,4),\\
(4,8), &\mbox{if } g_1^{**}(t;1,4) < x < g_2^*(t;4),\\
(4,6), &\mbox{if } g_2^*(t;4) < x < g_2(t;4),\\
(4,1), &\mbox{if } x > g_2(t;4),
\end{cases}
\end{equation*}
where
\begin{equation*}
g_1^{**}(t;1,4) = \frac{7}{6} + \frac{1}{10} \int_{t_{*2}}^t 
\exp \left( 0.0880 \frac{s}{1+s} \right) ds.
\end{equation*}
Since $\frac{dg_1^{**}}{dt} = \frac{1}{10\,h(t)} < \frac{1}{5\,h(t)} = \frac{dg_2^*}{dt}$, for all $t > t_{*2}$, the shock $S_1^{(*3)}$ is always slower than $J_2^{(*3)}$, so no further interactions occur. 
As $\epsilon \to 0$, the interaction points $(x_{*1}, t_{*1})$ and $(x_{*2}, t_{*2})$ converge to the origin, the curves $g_2^*(t;4)$ and $g_2(t;4)$ merge into a single contact discontinuity $J_2$, and $g_1^{**}(t;1,4)$ becomes the shock $S_1$, so the solution converges to the Riemann solution $S_1 + J_2$ of \eqref{sys} and \eqref{RP}, which is consistent with the stability result of Theorem~\ref{main_thm}.

\begin{figure}[h]
    \centering
    \includegraphics[width=0.45\linewidth]{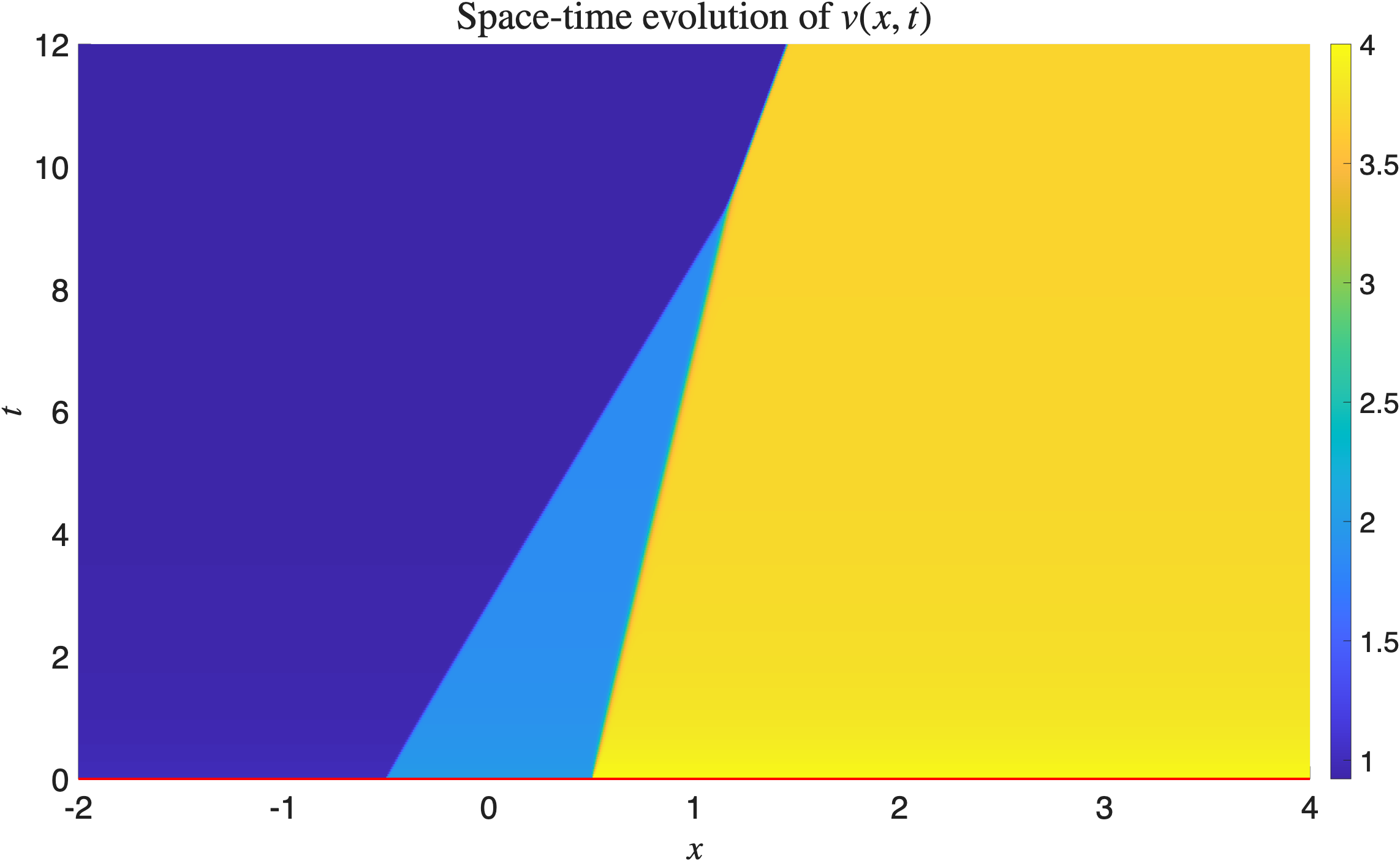}
    \includegraphics[width=0.45\linewidth]{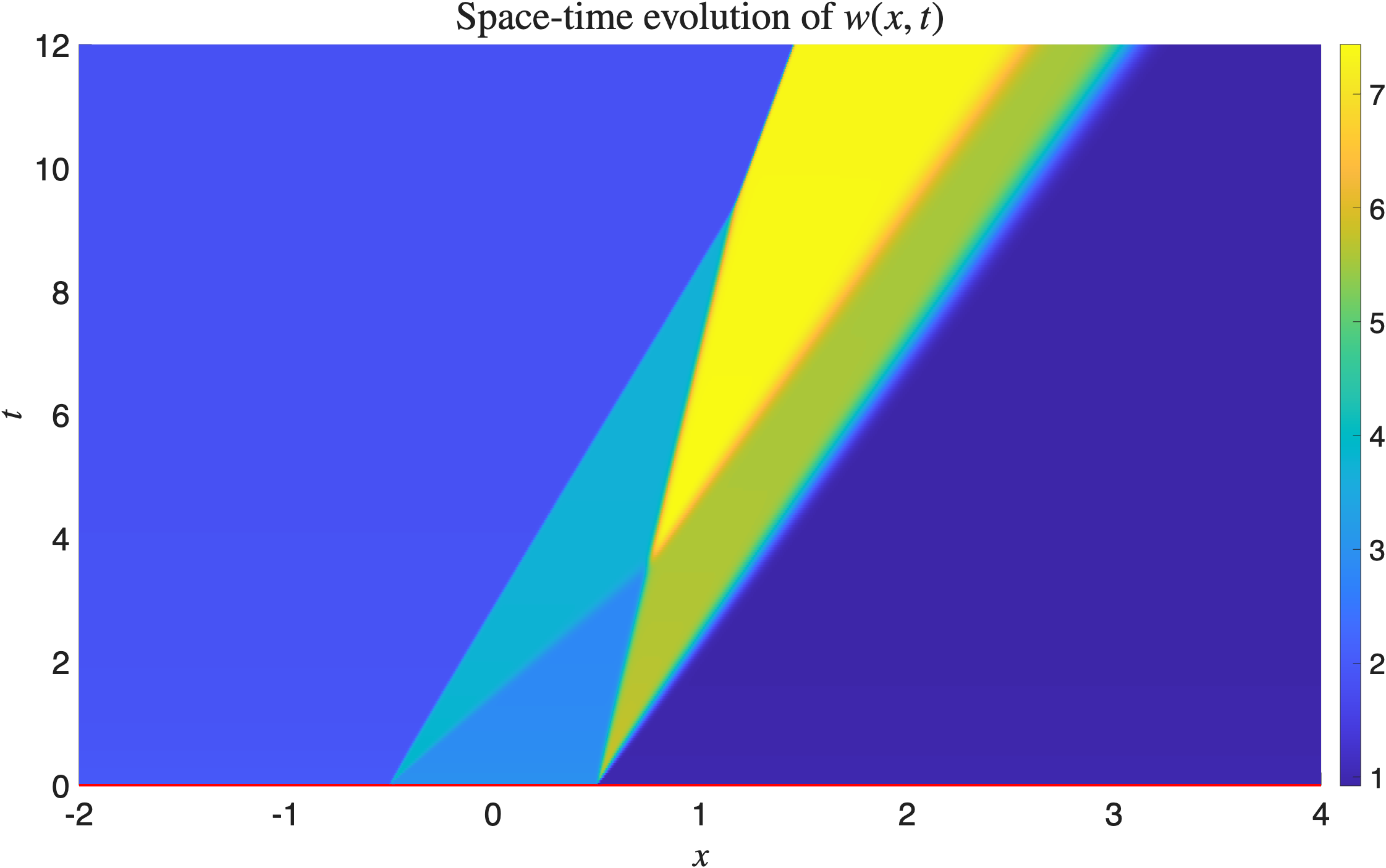}
    \caption{Space-time evolution of $v(x,t)$ (left) and $w(x,t)$ 
    (right) for Case~1 ($0 < v_- < v_\thicksim < v_+$).}
    \label{fig:case1}
\end{figure}

The numerical solution using the Lax-Friedrichs type scheme is shown in Figure~\ref{fig:case1}. In the left figure, we can observe the shock $S_1^{(-)}$ separating the states $v=1$ and $v=2$, and the shocks $S_1^{(\sim)}$ and $S_1^{(*1)}$ separating the states $v=2$ and $v=4$. The first interaction at $(x_{*1}, t_{*1}) \approx \left(\frac{3}{4}, 3.5647\right)$ and the second at $(x_{*2}, t_{*2}) \approx \left(\frac{7}{6}, 9.3599\right)$ are clearly visible as changes in the slope of the wave fronts. 
In the right figure, the intermediate states $W_{*1}=4$, $W_{*3}=8$ and $W_{*2}=6$ are distinguishable by the color scale, and the parallel contact discontinuities $J_2^{(*3)}$ and $J_2^{(*2)}$ are visible as the two parallel lines with the same slope for $t > t_{*1}$.

For the next cases of interaction waves, in similar way of this case, we use the Newton-Rapson method to approximate the times of interaction. Therefore, we only show the numerical solutions and report the times of the interaction.\\

{\bf Case 2: $0<v_\thicksim<v_-<v_+$.} (Subcase 2.1: $v_- < v_+$)

In Table~\ref{tb1} we show the initial data and the times and positions of the wave interactions. In this case, the local Riemann solution at $x = -\frac{1}{2}$ consists of a 1-rarefaction wave $R_1^{(-)}$ and a 2-contact discontinuity $J_2^{(*1)}$, with intermediate state $(V_{*1}, W_{*1}) = \left(1, \frac{2}{3}\right)$, while the local Riemann solution at $x = \frac{1}{2}$ consists of a 1-shock wave $S_1^{(\sim)}$ and a 2-contact discontinuity $J_2^{(*2)}$, with intermediate state $(V_{*2}, W_{*2}) = (5, 5)$.
The first interaction at $(x_{*1}, t_{*1})$ is between the contact discontinuity $J_2^{(*1)}$ and the shock wave $S_1^{(\sim)}$, producing a new local Riemann problem with left state $(V_{*1}, W_{*1}) = \left(1, \frac{2}{3}\right)$ and right state $(V_{*2}, W_{*2}) = (5, 5)$. Since $V_{*1} < V_{*2}$, this is solved by a new shock $S_1^{(*1)}$ with speed $\frac{1}{(1+V_{*1})(1+V_{*3})h(t)} = \frac{1}{12\,h(t)}$ and a new contact discontinuity $J_2^{(*3)}$ with speed $\frac{1}{(1+V_{*3})h(t)} = \frac{1}{6\,h(t)}$, with intermediate state $(V_{*3}, W_{*3}) = \left(5, \frac{10}{3}\right)$. Note that $S_1^{(*1)}$ and $J_2^{(*3)}$ propagate with the same speeds as $S_1^{(\sim)}$ and $J_2^{(*2)}$, respectively, but from the new point $(x_{*1}, t_{*1})$.
The second interaction at $(x_{*2}, t_{*2})$ is between the wave front of $R_1^{(-)}$, with speed $\frac{1}{(1+V_{*1})^2 h(t)} = \frac{1}{4\,h(t)}$, and the shock $S_1^{(*1)}$, with speed $\frac{1}{12\,h(t)}$. Since $\frac{1}{4} > \frac{1}{12}$, the wave front of $R_1^{(-)}$ catches $S_1^{(*1)}$ and the shock begins to cross the rarefaction wave $R_1^{(-)}$. Since $v_- = 3 < v_+ = 5$, the shock $S_1^{(*1)}$ fully crosses $R_1^{(-)}$, and the third interaction at $(x_{*3}, t_{*3})$ produces a new local Riemann problem with left state $(v_-, w_-) = (3, 2)$ and right state $(V_{*3}, W_{*3}) = \left(5, \frac{10}{3}\right)$. Since $\frac{w_-}{v_-} = \frac{2}{3} = \frac{W_{*3}}{V_{*3}}$, the contact discontinuity is trivial and the solution consists of the single shock $S_1^{(*3)}$ with speed $\frac{1}{(1+v_-)(1+V_{*3})h(t)} = \frac{1}{24\,h(t)}$, after which no further interactions occur.

\begin{table}[h]
\centering
\begin{tabular}{|c|c|c|} \hline
Initial data & Times of interaction & Space of interaction\\ \hline
$v_- = 3$, $v_\thicksim = 1$, $v_+ = 5$ & $t_{*1} \approx 2.3002$ 
& $x_{*1} = 0.7000$ \\ 
$w_- = 2$, $w_\thicksim = 1$, $w_+ = 3$ & $t_{*2} \approx 5.6582$ 
& $x_{*2} = 1.0000$ \\
& $t_{*3} \approx 88.3000$ & $x_{*3} = 5.5000$ \\ \hline
\end{tabular}
\caption{Initial data, times and positions of interactions for 
Case~2, Subcase~2.1 ($0 < v_\thicksim < v_- < v_+$).} 
\label{tb1}
\end{table}

\begin{figure}[h]
    \centering
    \includegraphics[width=0.45\linewidth]{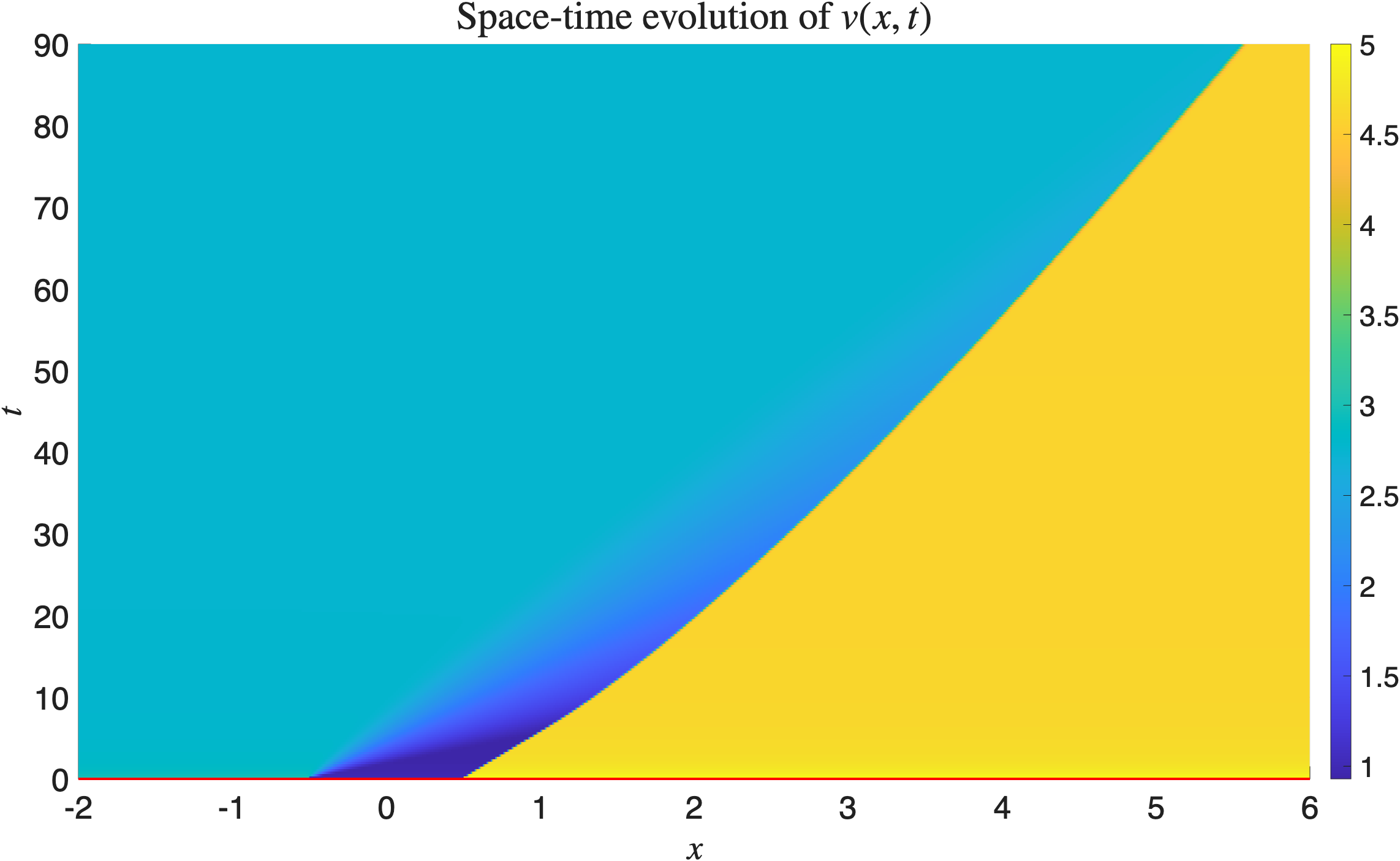}
    \includegraphics[width=0.45\linewidth]{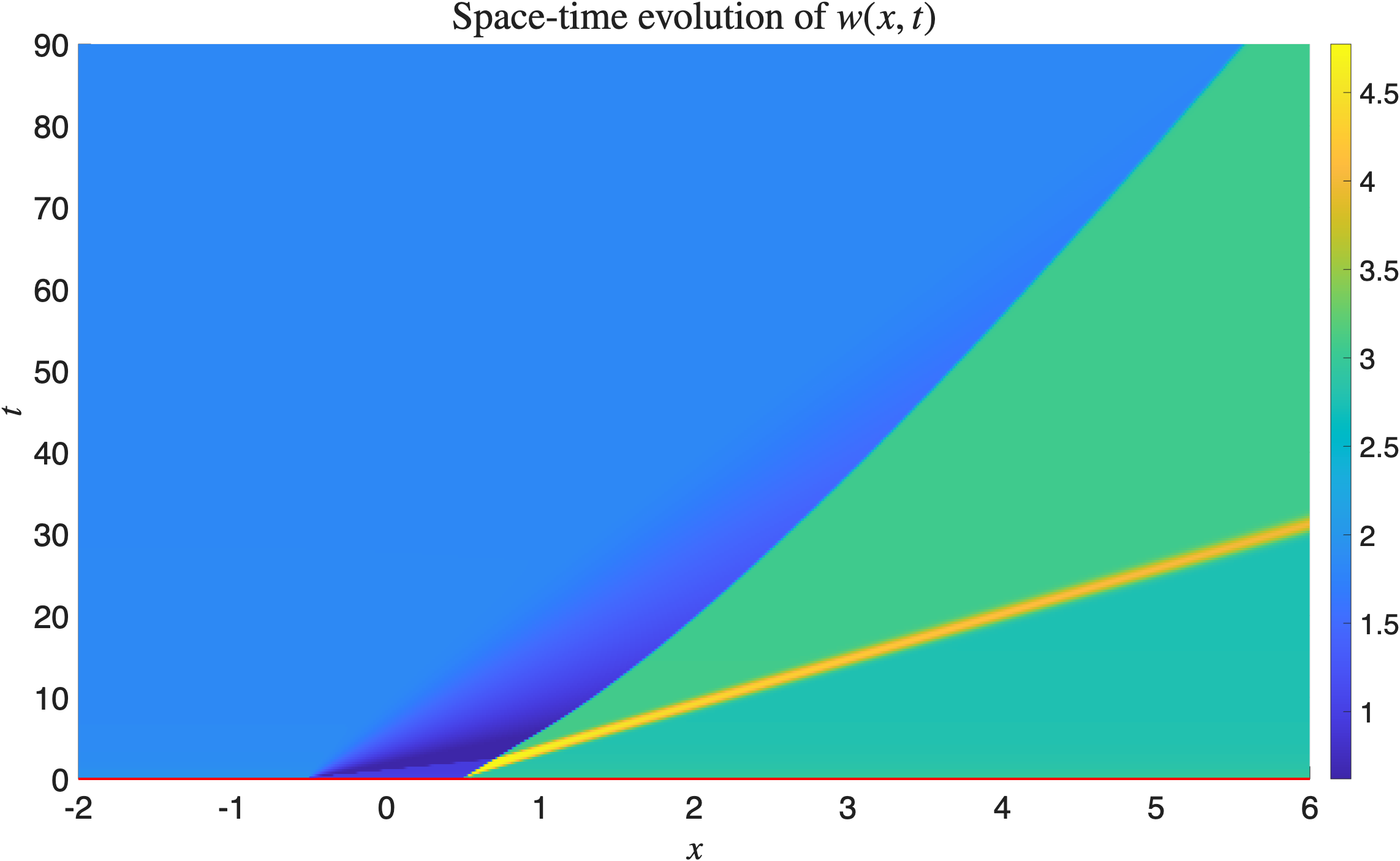}
    \caption{Space-time evolution of $v(x,t)$ (left) and $w(x,t)$ 
    (right) for Case~2, Subcase~2.1 
    ($0 < v_\thicksim < v_- < v_+$).}
    \label{fig:case2Sub21}
\end{figure}

The numerical solution using the Lax-Friedrichs type scheme is shown in Figure~\ref{fig:case2Sub21}. In the left figure, the rarefaction wave $R_1^{(-)}$ is visible as the smooth transition between the states $v = 3$ and $v = 1$ near the origin, while the curved boundary between $t_{*2} \approx 5.6582$ and $t_{*3} \approx 88.3000$ corresponds to the penetration of $S_1^{(*1)}$ through $R_1^{(-)}$, after which the boundary becomes the straight shock $S_1^{(*3)}$. In the right figure, the contact discontinuity $J_2^{(*2)}$ is clearly visible as the sharp orange line separating the intermediate state $W_{*3} = \frac{10}{3}$ from the right state $w_+ = 3$.

{\bf Case 2: $0<v_\thicksim<v_+<v_-$.} (Subcase 2.2)

In this subcase, the initial data and the times and positions of the wave interactions are showed in Table~\ref{tb2}. The local Riemann solution at $x = -\frac{1}{2}$ consists of a 1-rarefaction wave $R_1^{(-)}$ and a 2-contact discontinuity $J_2^{(*1)}$, with intermediate state $(V_{*1}, W_{*1}) = \left(1, \frac{1}{2}\right)$, while the local Riemann solution at $x = \frac{1}{2}$ consists of a 1-shock wave $S_1^{(\sim)}$ and a 2-contact discontinuity $J_2^{(*2)}$, with intermediate state $(V_{*2}, W_{*2}) = (2, 2)$.
The first interaction at $(x_{*1}, t_{*1})$ is between $J_2^{(*1)}$ and $S_1^{(\sim)}$, producing a new shock $S_1^{(*1)}$ with speed $\frac{1}{(1+V_{*1})(1+V_{*3})h(t)} = \frac{1}{6\,h(t)}$ and a new contact discontinuity $J_2^{(*3)}$ with speed $\frac{1}{(1+V_{*3})h(t)} = \frac{1}{3\,h(t)}$, with intermediate state $(V_{*3}, W_{*3}) = (2, 1)$. Since $V_{*3} = V_{*2} = 2$, the contact discontinuities $J_2^{(*3)}$ and $J_2^{(*2)}$ satisfy $\frac{dg_2^*}{dt} = \frac{dg_2}{dt} = \frac{1}{3h(t)}$ for all $t > t_{*1}$, so they are parallel and will never interact.
The second interaction at $(x_{*2}, t_{*2})$ is between the wave front of $R_1^{(-)}$, with speed $\frac{1}{(1+V_{*1})^2 h(t)} = \frac{1}{4\,h(t)}$, and the shock $S_1^{(*1)}$, with speed $\frac{1}{6\,h(t)}$. Since $\frac{1}{4} > \frac{1}{6}$, the wave front of $R_1^{(-)}$ catches $S_1^{(*1)}$ and the shock begins to cross the rarefaction wave $R_1^{(-)}$. However, since $v_+ = 2 < v_- = 4$, the shock $S_1^{(*1)}$ cannot fully cross $R_1^{(-)}$, and the curve 
\begin{equation*}
x + \frac{1}{2} = \frac{1}{(1+V_{*3})^2} \int_0^t \exp \left(0.0880\frac{s}{1+s}\right)ds = \frac{1}{9}\int_0^t \exp \left(0.0880\frac{s}{1+s}\right)ds
\end{equation*}
becomes its asymptote. Therefore, for sufficiently large $t > t_{*2}$, a residual rarefaction wave $R_1$ remains, and no further interactions occur.

\begin{table}[h]
\centering
\begin{tabular}{|c|c|c|} \hline
Initial data & Times of interaction & Space of interaction\\ \hline
$v_- = 4$, $v_\thicksim = 1$, $v_+ = 2$ & $t_{*1} \approx 2.8634$ 
& $x_{*1} = 1.0000$ \\ 
$w_- = 2$, $w_\thicksim = 1$, $w_+ = 3$ & $t_{*2} \approx 11.2058$ 
& $x_{*2} = 2.5000$ \\ \hline
\end{tabular}
\caption{Initial data, times and positions of interactions for 
Case~2, Subcase~2.2 ($0 < v_\thicksim < v_+ < v_-$).}
\label{tb2}
\end{table}

\begin{figure}[h]
    \centering
    \includegraphics[width=0.45\linewidth]{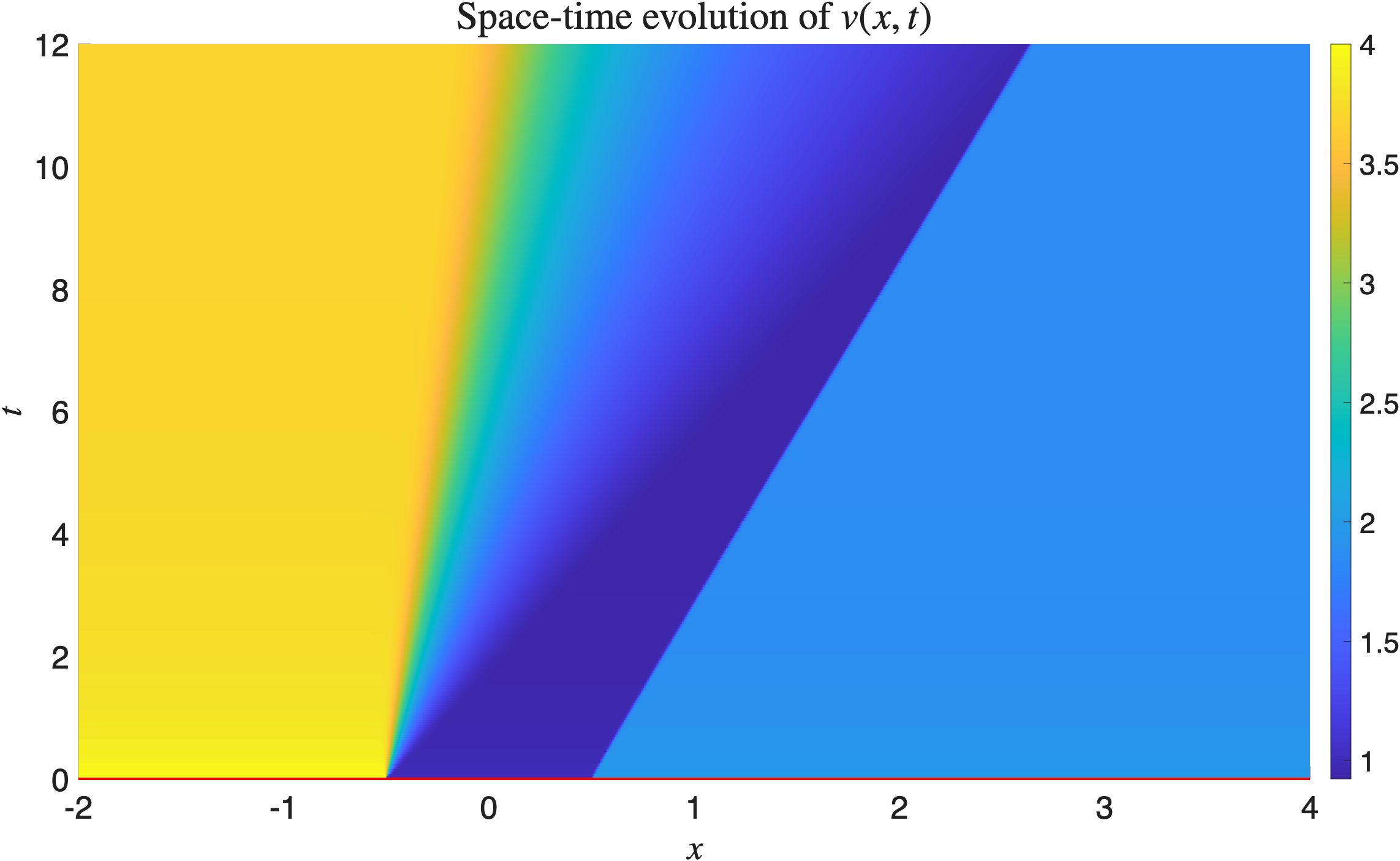}
    \includegraphics[width=0.45\linewidth]{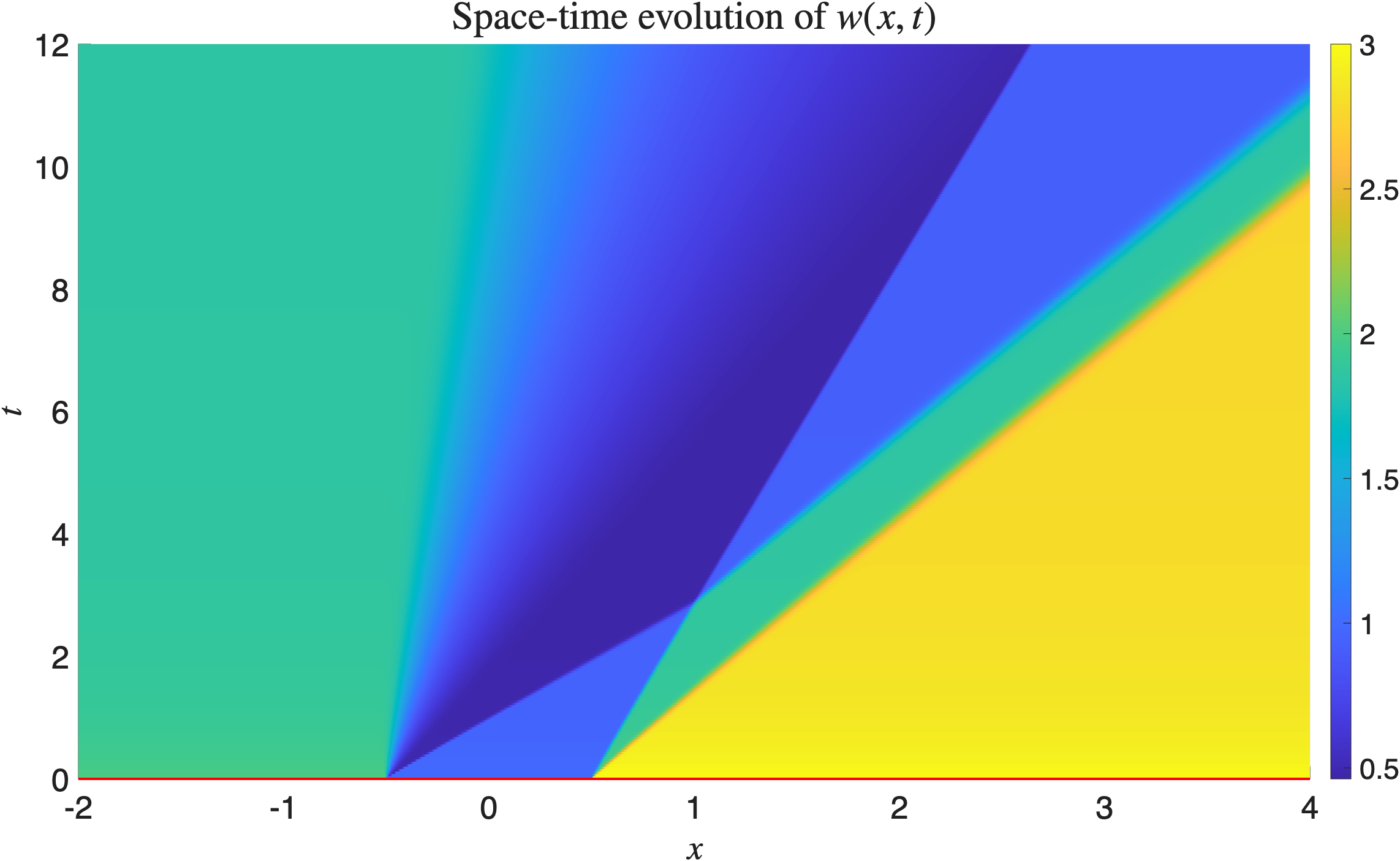}
    \caption{Space-time evolution of $v(x,t)$ (left) and $w(x,t)$ 
    (right) for Case~2, Subcase~2.2 
    ($0 < v_\thicksim < v_+ < v_-$).}
    \label{fig:case2Sub22}
\end{figure}

The numerical solution using the Lax-Friedrichs type scheme is shown in Figure~\ref{fig:case2Sub22}. In the left figure, the smooth transition between $v = 4$ and $v = 2$ persisting for all $t > t_{*2}$ is consistent with the fact that $S_1^{(*1)}$ cannot fully cross $R_1^{(-)}$, and illustrates the presence of the residual rarefaction $R_1$. In the right figure, the two parallel contact discontinuities $J_2^{(*3)}$ and $J_2^{(*2)}$ are clearly visible as two lines with the same slope for $t > t_{*1}$.

{\bf Case 3: $0<v_-<v_+<v_\thicksim$.} (Subcase 3.1: $v_- < v_+$)

In this case, the initial data and the times and positions of the wave interactions are showed in Table~\ref{tb3}. The local Riemann solution at $x = -\frac{1}{2}$ consists of a 1-shock wave $S_1^{(-)}$ and a 2-contact discontinuity $J_2^{(*1)}$, with intermediate state $(V_{*1}, W_{*1}) = (4, 8)$, while the local Riemann solution at $x = \frac{1}{2}$ consists of a 1-rarefaction wave $R_1^{(\sim)}$ and a 2-contact discontinuity $J_2^{(*2)}$, with intermediate state $(V_{*2}, W_{*2}) = \left(3, \frac{9}{4}\right)$.
The first interaction at $(x_{*1}, t_{*1})$ is between $J_2^{(*1)}$ and the tail of $R_1^{(\sim)}$, after which $J_2^{(*1)}$ begins to cross $R_1^{(\sim)}$. The second interaction at $(x_{*2}, t_{*2})$ occurs when $J_2^{(*1)}$ exits the front of $R_1^{(\sim)}$, producing a new rarefaction wave $R_1^{(*1)}$ and a new contact discontinuity $J_2^{(*3)}$ with intermediate state $(V_{*3}, W_{*3}) = (3, 6)$. Since $V_{*3} = V_{*2} = 3$, the contact discontinuities $J_2^{(*3)}$ and $J_2^{(*2)}$ satisfy $\frac{dg_2^*}{dt} = \frac{dg_2}{dt} = \frac{1}{4h(t)}$ for all $t > t_{*2}$, so they are parallel and will never interact.
The third interaction at $(x_{*3}, t_{*3})$ is between $S_1^{(-)}$ and the tail of $R_1^{(*1)}$, after which $S_1^{(-)}$ begins to cross $R_1^{(*1)}$. Since $v_- = 1 < v_+ = 3$, the shock $S_1^{(-)}$ fully crosses $R_1^{(*1)}$, producing at the fourth interaction $(x_{*4}, t_{*4})$ a new local Riemann problem with left state $(v_-, w_-) = (1, 2)$ and right state $(V_{*3}, W_{*3}) = (3, 6)$. Since $\frac{w_-}{v_-} = \frac{2}{1} = 2 = \frac{6}{3} = \frac{W_{*3}}{V_{*3}}$, the contact discontinuity is trivial and the solution consists of the single shock $S_1^{(*3)}$ with speed $\frac{1}{(1+v_-)(1+V_{*3})h(t)} = \frac{1}{8\,h(t)}$, after which no further interactions occur.

\begin{table}[h]
\centering
\begin{tabular}{|c|c|c|} \hline
Initial data & Times of interaction & Space of interaction\\ \hline
$v_- = 1$, $v_\thicksim = 4$, $v_+ = 3$ & $t_{*1} \approx 5.8901$ 
& $x_{*1} = 0.7500$ \\ 
$w_- = 2$, $w_\thicksim = 3$, $w_+ = 1$ & $t_{*2} \approx 6.6883$ 
& $x_{*2} \approx 0.9444$ \\
& $t_{*3} \approx 15.5058$ & $x_{*3} \approx 1.1667$ \\
& $t_{*4}$ very large & \\ \hline
\end{tabular}
\caption{Initial data, times and positions of interactions for Case~3, Subcase~3.1 ($0 < v_- < v_+ < v_\thicksim$).}
\label{tb3}
\end{table}

\begin{figure}[h]
    \centering
    \includegraphics[width=0.45\linewidth]{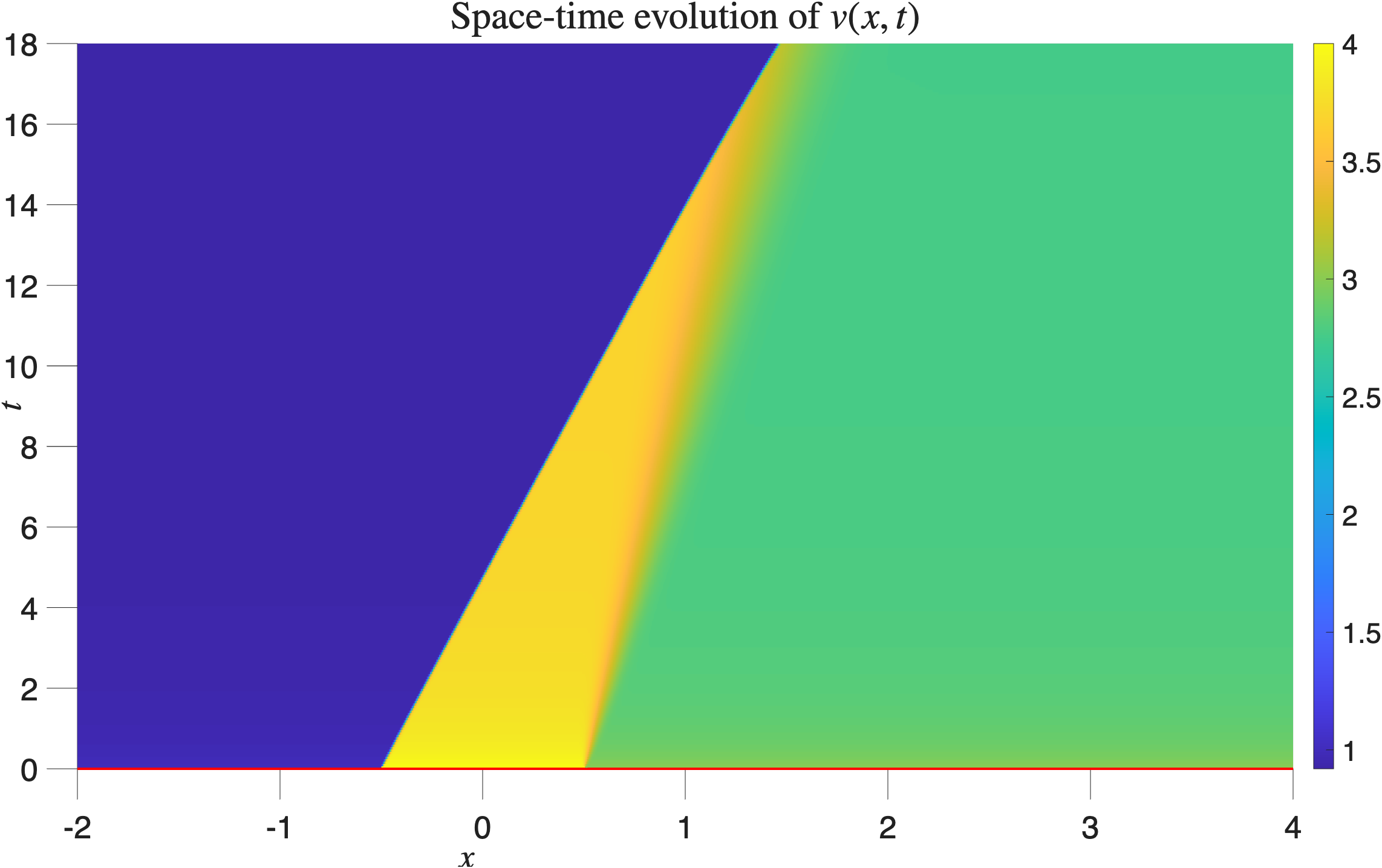}
    \includegraphics[width=0.45\linewidth]{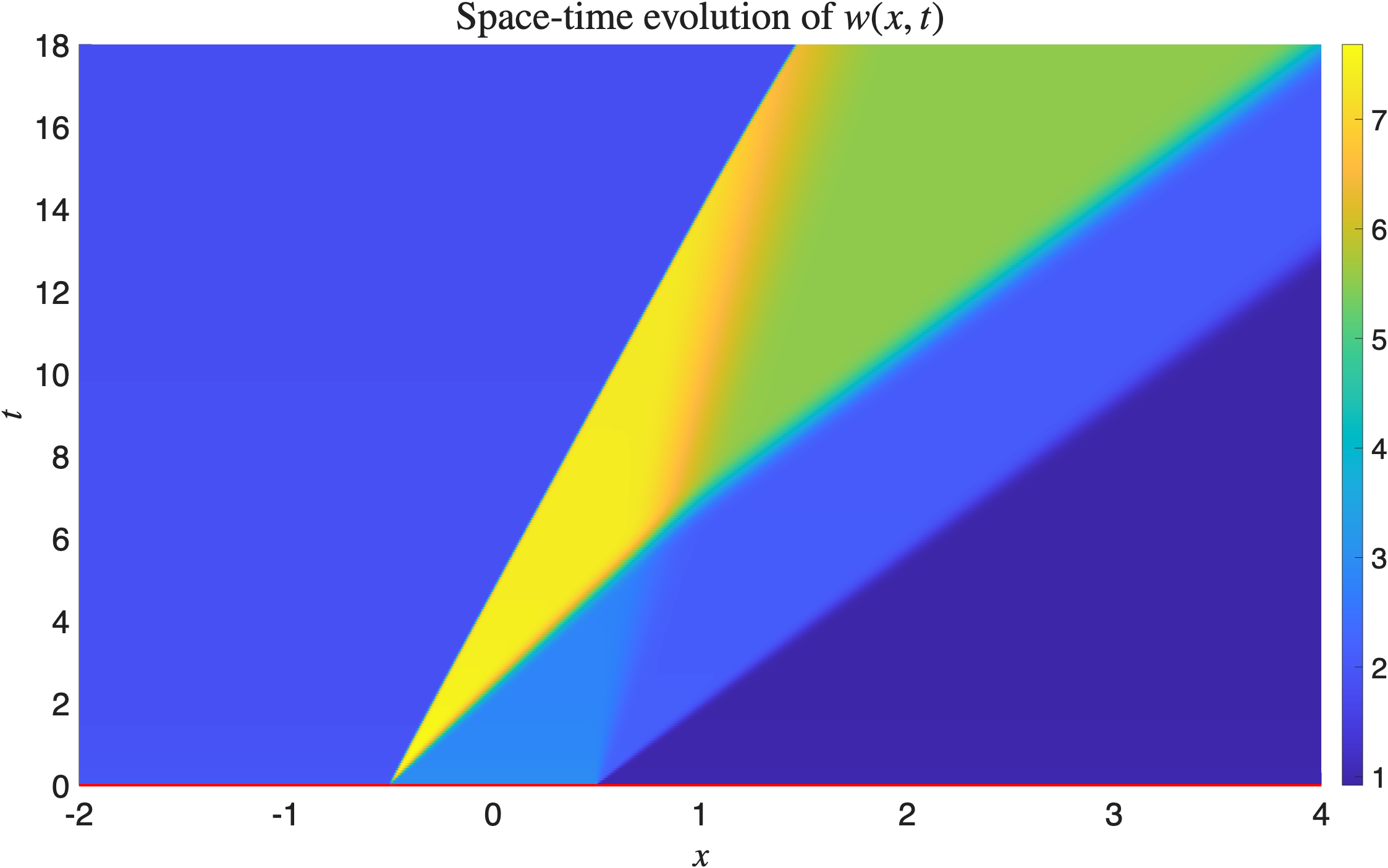}
    \caption{Space-time evolution of $v(x,t)$ (left) and $w(x,t)$ (right) for Case~3, Subcase~3.1 ($0 < v_- < v_+ < v_\thicksim$), with $T = 18$.}
    \label{fig:case3Sub31}
\end{figure}

The numerical solution using the Lax-Friedrichs type scheme is shown in Figure~\ref{fig:case3Sub31}. In the left figure, the narrow yellow strip near $x\approx 0$ for $t_{*1} \leq t \leq t_{*2}$ corresponds to the intermediate state $V_{*1}=4$ during the crossing of $J_2^{(*1)}$ through $R_1^{(\sim)}$, while the smooth left boundary of the green region for $t > t_{*3} \approx 15.5058$ corresponds to the crossing of $S_1^{(-)}$ through $R_1^{(*1)}$, which is still ongoing at $T=18$ since $t_{*4}$ is very large. In the right figure, the intermediate states $W_{*1}=8$ (yellow), $W_{*3}=6$ (green) and $W_{*2}=\frac{9}{4}$ (cyan) are clearly distinguishable, and the two parallel contact discontinuities $J_2^{(*3)}$ and $J_2^{(*2)}$ are visible with the same slope for $t > t_{*2}$.

{\bf Case 3: $0<v_+<v_-<v_\thicksim$.} (Subcase 3.2: $v_+ < v_-$)

In this subcase, the initial data and the times and positions of the wave interactions are showed in Table~\ref{tb4}. The local Riemann solution at $x = -\frac{1}{2}$ consists of a 1-shock wave $S_1^{(-)}$ and a 2-contact discontinuity $J_2^{(*1)}$, with intermediate state $(V_{*1}, W_{*1}) = (4, 4)$, while the local Riemann solution at $x = \frac{1}{2}$ consists of a 1-rarefaction wave $R_1^{(\sim)}$ and a 2-contact discontinuity $J_2^{(*2)}$, with intermediate state $(V_{*2}, W_{*2}) = \left(1, \frac{3}{4}\right)$.
The first interaction at $(x_{*1}, t_{*1})$ is between $J_2^{(*1)}$ and the tail of $R_1^{(\sim)}$, after which $J_2^{(*1)}$ begins to cross $R_1^{(\sim)}$. The second interaction at $(x_{*2}, t_{*2})$ occurs when $J_2^{(*1)}$ exits the front of $R_1^{(\sim)}$, producing a new rarefaction wave $R_1^{(*1)}$ and a new contact discontinuity $J_2^{(*3)}$ with intermediate state $(V_{*3}, W_{*3}) = (1, 1)$. Since $V_{*3} = V_{*2} = 1$, the contact discontinuities $J_2^{(*3)}$ and $J_2^{(*2)}$ satisfy $\frac{dg_2^*}{dt} = \frac{dg_2}{dt} = \frac{1}{2h(t)}$ for all $t > t_{*2}$, so they are parallel and will never interact.
The third interaction at $(x_{*3}, t_{*3})$ is between $S_1^{(-)}$ and the tail of $R_1^{(*1)}$, after which $S_1^{(-)}$ begins to cross $R_1^{(*1)}$. However, since $v_+ = 1 < v_- = 2$, the shock $S_1^{(-)}$ cannot fully cross $R_1^{(*1)}$, and the 1-characteristic curve $\frac{1}{(1+v_-)^2 h(t)} = \frac{1}{9h(t)}$ inside $R_1^{(*1)}$ becomes its asymptote. Therefore, for sufficiently large $t > t_{*3}$, a residual rarefaction wave $R_1$ remains, and no further interactions occur. 

\begin{table}[h]
\centering
\begin{tabular}{|c|c|c|} \hline
Initial data & Times of interaction & Space of interaction\\ \hline
$v_- = 2$, $v_\thicksim = 4$, $v_+ = 1$ & $t_{*1} \approx 5.8901$ 
& $x_{*1} = 0.7500$ \\ 
$w_- = 2$, $w_\thicksim = 3$, $w_+ = 1$ & $t_{*2} \approx 14.8920$ 
& $x_{*2} = 4.5000$ \\
& $t_{*3} \approx 34.6518$ & $x_{*3} = 2.0000$ \\ \hline
\end{tabular}
\caption{Initial data, times and positions of interactions for Case~3, Subcase~3.2 ($0 < v_+ < v_- < v_\thicksim$).}
\label{tb4}
\end{table}

\begin{figure}[h]
    \centering
    \includegraphics[width=0.45\linewidth]{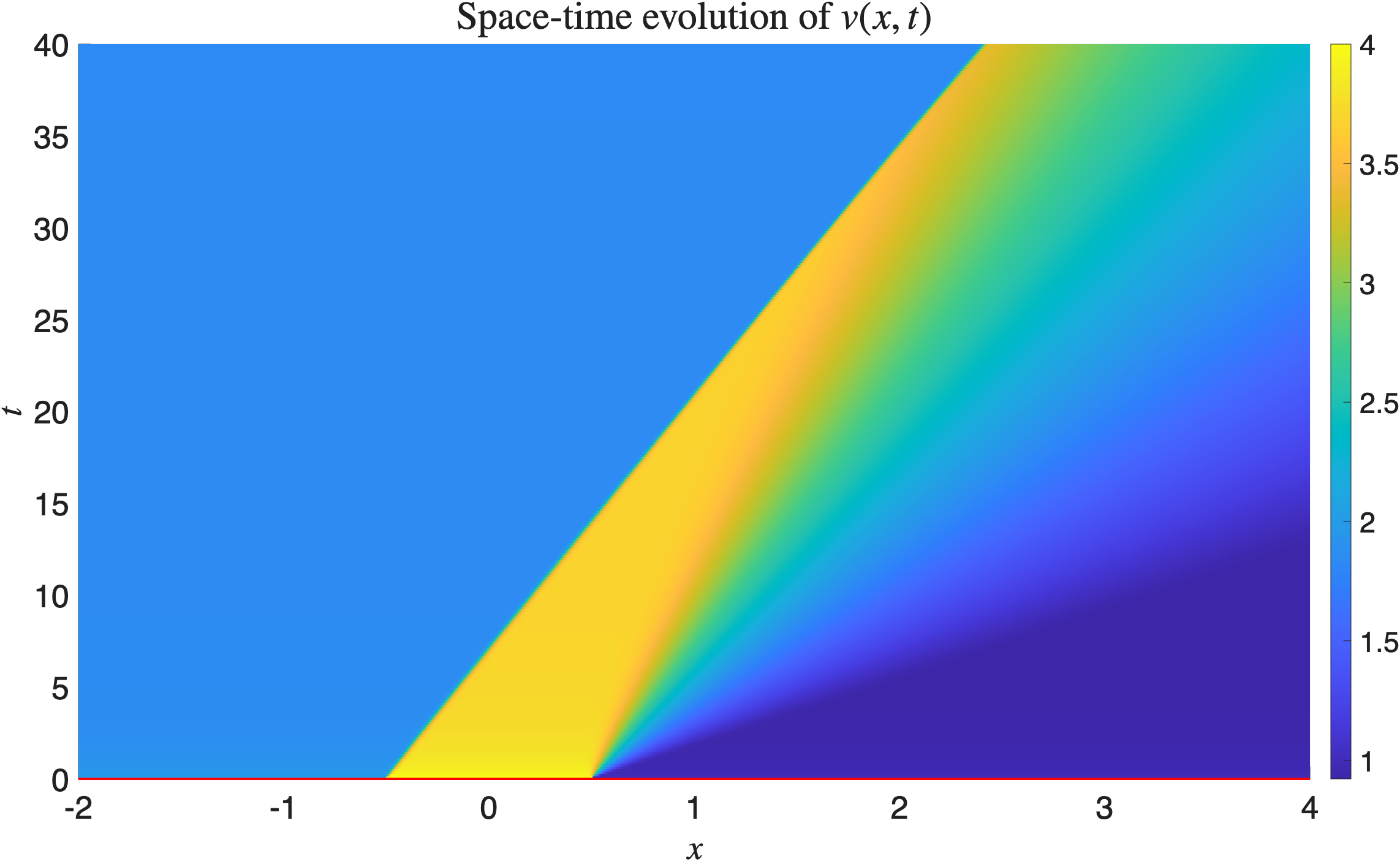}
    \includegraphics[width=0.45\linewidth]{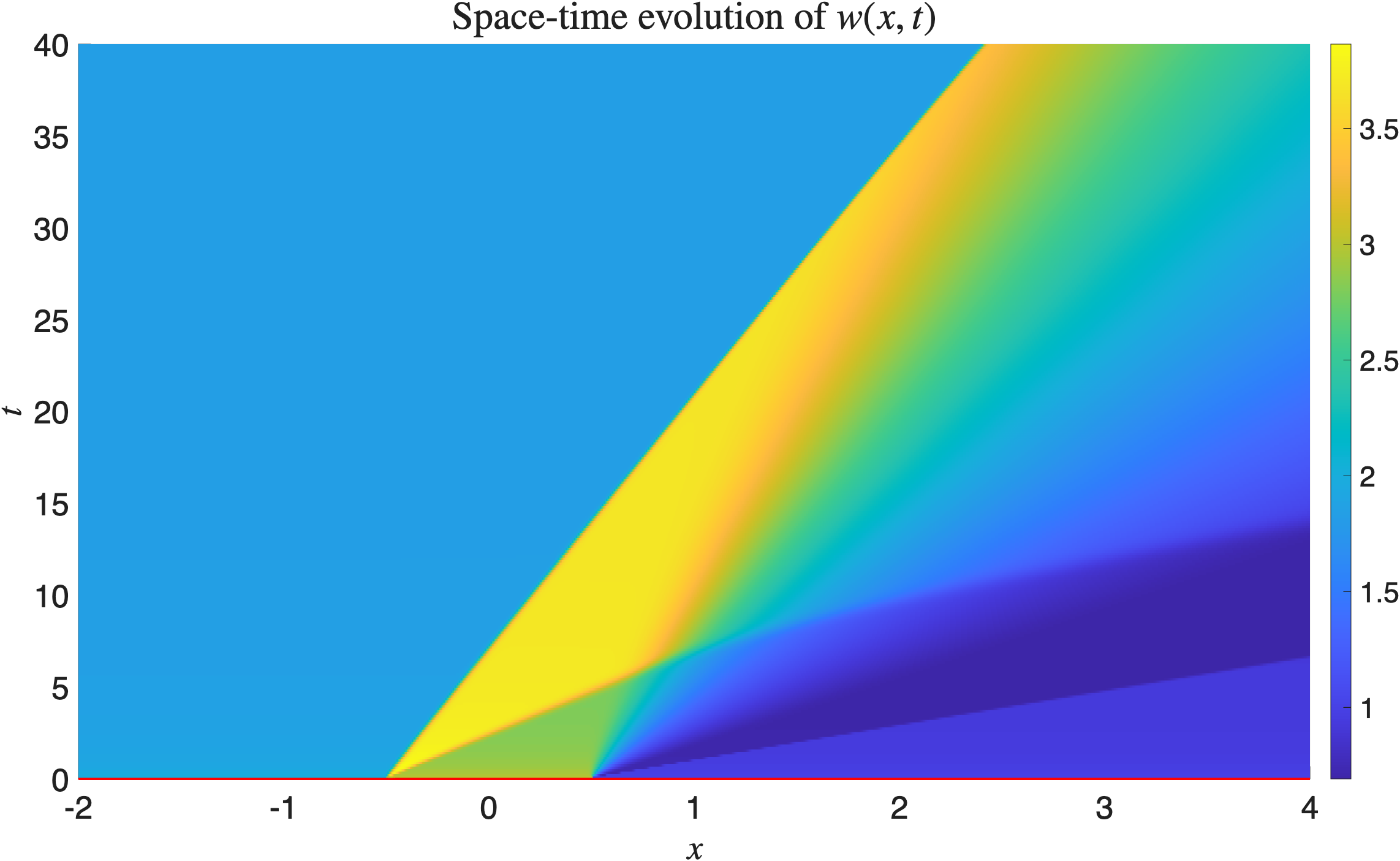}
    \caption{Space-time evolution of $v(x,t)$ (left) and $w(x,t)$ 
    (right) for Case~3, Subcase~3.2 
    ($0 < v_+ < v_- < v_\thicksim$), with $T = 40$.}
    \label{fig:case3Sub32}
\end{figure}

The numerical solution using the Lax-Friedrichs type scheme is shown in Figure~\ref{fig:case3Sub32}. In the left figure, the yellow/orange band corresponds to the intermediate state $V_{*1}=4$ during the crossing of $J_2^{(*1)}$ through $R_1^{(\sim)}$ between $t_{*1}$ and $t_{*2}$. The smooth transition persisting for all $t > t_{*3} \approx 34.6518$ is consistent with the fact that $S_1^{(-)}$ cannot fully cross $R_1^{(*1)}$, and illustrates the presence of the residual rarefaction $R_1$, in agreement with the limiting solution $R_1 + J_2$.
In the right figure, since $(V_{*3}, W_{*3}) = (V_{*2}, W_{*2}) = (1,1)$, the two parallel contact discontinuities $J_2^{(*3)}$ and $J_2^{(*2)}$ carry the same state and appear merged into a single boundary.

{\bf Case 4: $0<v_+<v_\thicksim<v_-$.}

For this case, Table~\ref{tb5} contains the initial data and the times and positions of the wave interactions.
Both local Riemann solutions consist of rarefaction waves and contact discontinuities. The local Riemann solution at $x = -\frac{1}{2}$ consists of a 1-rarefaction wave $R_1^{(-)}$ and a 2-contact discontinuity $J_2^{(*1)}$, with intermediate state $(V_{*1}, W_{*1}) = (2, 1)$, while the local Riemann solution at $x = \frac{1}{2}$ consists of a 1-rarefaction wave $R_1^{(\sim)}$ and a 2-contact discontinuity $J_2^{(*2)}$, with intermediate state $(V_{*2}, W_{*2}) = \left(1, \frac{1}{2}\right)$.
The first interaction at $(x_{*1}, t_{*1})$ is between $J_2^{(*1)}$ and the tail of $R_1^{(\sim)}$, after which $J_2^{(*1)}$ begins to cross the rarefaction wave $R_1^{(\sim)}$. The second interaction at $(x_{*2}, t_{*2})$ occurs when $J_2^{(*1)}$ exits the front of $R_1^{(\sim)}$, producing a new rarefaction wave $R_1^{(*3)}$ and a new contact discontinuity $J_2^{(*3)}$ with intermediate state $(V_{*3}, W_{*3}) = \left(1, \frac{1}{2}\right)$. Since $V_{*3} = V_{*2} = 1$, the contact discontinuities $J_2^{(*3)}$ and $J_2^{(*2)}$ satisfy $\frac{dg_2^*}{dt} = \frac{dg_2}{dt} = \frac{1}{2h(t)}$ for all $t > t_{*2}$, so they are parallel and will never interact. Furthermore, the head of $R_1^{(-)}$ and the tail of $R_1^{(*3)}$ propagate with the same speed $\frac{1}{(1+V_{*1})^2 h(t)} = \frac{1}{9h(t)}$, so the two rarefaction waves $R_1^{(-)}$ and $R_1^{(*3)}$ never interact. Therefore, no further interactions occur after $t_{*2}$. 

\begin{table}[h]
\centering
\begin{tabular}{|c|c|c|} \hline
Initial data & Times of interaction & Space of interaction\\ \hline
$v_- = 4$, $v_\thicksim = 2$, $v_+ = 1$ & $t_{*1} \approx 4.2640$ 
& $x_{*1} = 1.0000$ \\ 
$w_- = 2$, $w_\thicksim = 1$, $w_+ = 3$ & $t_{*2} \approx 7.5112$ 
& $x_{*2} = 2.5000$ \\ \hline
\end{tabular}
\caption{Initial data, times and positions of interactions for Case~4 ($0 < v_+ < v_\thicksim < v_-$).}
\label{tb5}
\end{table}

\begin{figure}[h]
    \centering
    \includegraphics[width=0.45\linewidth]{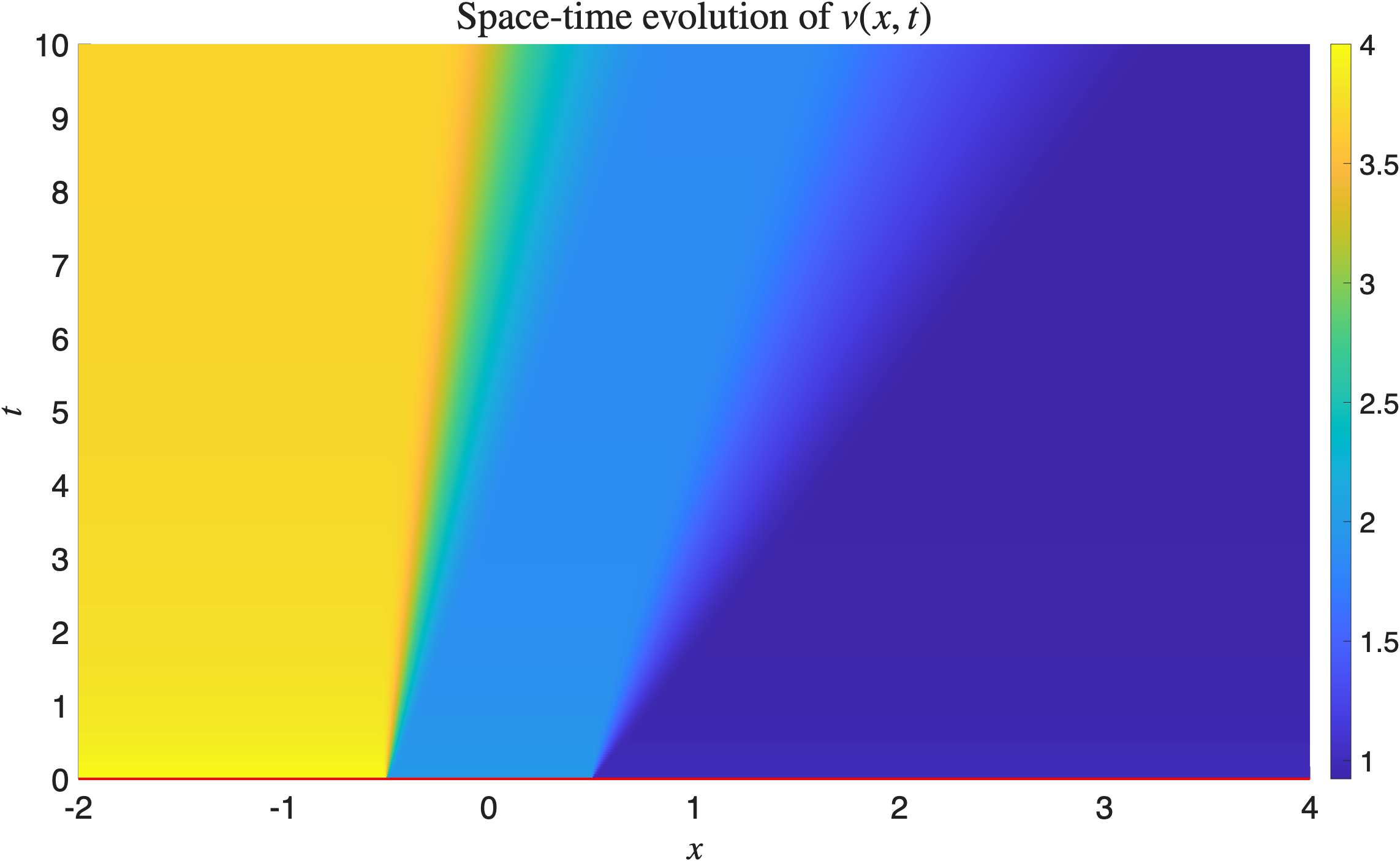}
    \includegraphics[width=0.45\linewidth]{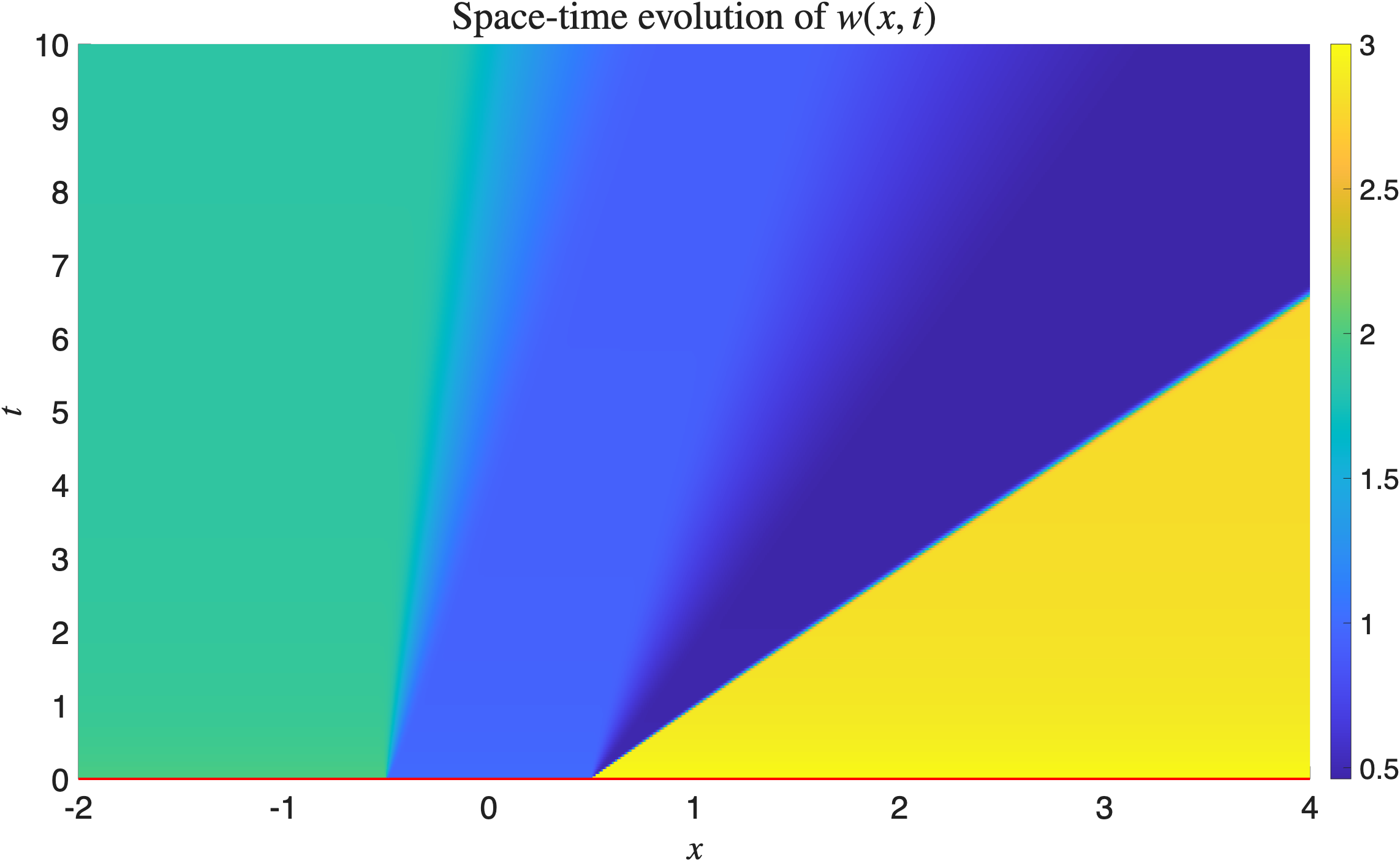}
    \caption{Space-time evolution of $v(x,t)$ (left) and $w(x,t)$ (right) for Case~4 ($0 < v_+ < v_\thicksim < v_-$), with $T = 10$.}
    \label{fig:case4}
\end{figure}

The numerical solution using the Lax-Friedrichs type scheme is shown in Figure~\ref{fig:case4}. In the left figure, the two smooth transitions between $v=4$ and $v=2$, and between $v=2$ and $v=1$, correspond respectively to the rarefaction waves $R_1^{(-)}$ and $R_1^{(*3)}$, which is consistent with the fact that their head and tail propagate at the same speed $\frac{1}{9h(t)}$ and they never interact. In the right figure, since $(V_{*3},W_{*3})=(V_{*2},W_{*2})=\left(1,\frac{1}{2}\right)$, the contact discontinuities $J_2^{(*3)}$ and $J_2^{(*2)}$ carry the same state and appear as a single sharp line separating the region $w\approx \frac{1}{2}$ from the right state $w_+=3$.

{\bf Case 5: $v_\thicksim=0$, $v_-<v_+$.} (Subcase 5.1)

For this subcase, Table~\ref{tb6} contains the initial data and the times and positions of the wave interactions.
Since $v_\thicksim = 0$, the local Riemann solution at $x = -\frac{1}{2}$ consists of a composite wave $R_1^{(-)}J_2^{(*1)}$, in which the wave front of $R_1^{(-)}$ and $J_2^{(*1)}$ are supported on the same curve $x + \frac{1}{2} = \int_0^t \exp \left(0.0880 \frac{s}{1+s}\right)ds$, with speed $\frac{1}{h(t)}$. The local Riemann solution at $x = \frac{1}{2}$ is a delta shock wave $\delta S$, with speed $\frac{1}{(1+v_+)h(t)} = \frac{1}{4h(t)}$.
The first interaction at $(x_{*1}, t_{*1})$ occurs when the composite wave $R_1^{(-)}J_2^{(*1)}$ meets the delta shock wave $\delta S$. Since $v_- \neq 0$, the delta shock wave $\delta S$ splits into a 1-shock wave $S_1$ and a delta contact discontinuity $\delta J_2$, with intermediate state $(V_{*2}, W_{*2}) = (3, 6)$ and strength $\alpha(t_{*1}) = h(t_{*1})$. The delta contact discontinuity $\delta J_2$ propagates from $(x_{*1}, t_{*1})$ with speed $\frac{1}{(1+v_+)h(t)} = \frac{1}{4h(t)}$, while $S_1$ begins to cross the rarefaction wave $R_1^{(-)}$.
Since $v_- = 1 < v_+ = 3$, the shock $S_1$ fully crosses $R_1^{(-)}$, and the second interaction at $(x_{*2}, t_{*2})$ produces a new local Riemann problem with left state $(v_-, w_-) = (1, 2)$ and right state $(V_{*2}, W_{*2}) = (3, 6)$. Since $\frac{w_-}{v_-} = \frac{2}{1} = 2 = \frac{6}{3} = \frac{W_{*2}}{V_{*2}}$, the contact discontinuity is trivial and the solution consists of the single shock $S_1^{(*2)}$ with speed $\frac{1}{(1+v_-)(1+v_+)h(t)} = \frac{1}{8h(t)}$. Since $\frac{1}{8h(t)} < \frac{1}{4h(t)}$, the shock $S_1^{(*2)}$ is always slower than $\delta J_2$, so no further interactions occur. 

\begin{table}[h]
\centering
\begin{tabular}{|c|c|c|} \hline
Initial data & Times of interaction & Space of interaction\\ \hline
$v_- = 1$, $v_\thicksim = 0$, $v_+ = 3$ & $t_{*1} \approx 1.2919$ 
& $x_{*1} \approx 0.8333$ \\ 
$w_- = 2$, $w_\thicksim = 1$, $w_+ = 1$ & $t_{*2} \approx 11.2058$ 
& $x_{*2} = 2.5000$ \\ \hline
\end{tabular}
\caption{Initial data, times and positions of interactions for Case~5, Subcase~5.1 ($v_\thicksim = 0$, $v_- < v_+$).}
\label{tb6}
\end{table}

\begin{figure}[h]
    \centering
    \includegraphics[width=0.45\linewidth]{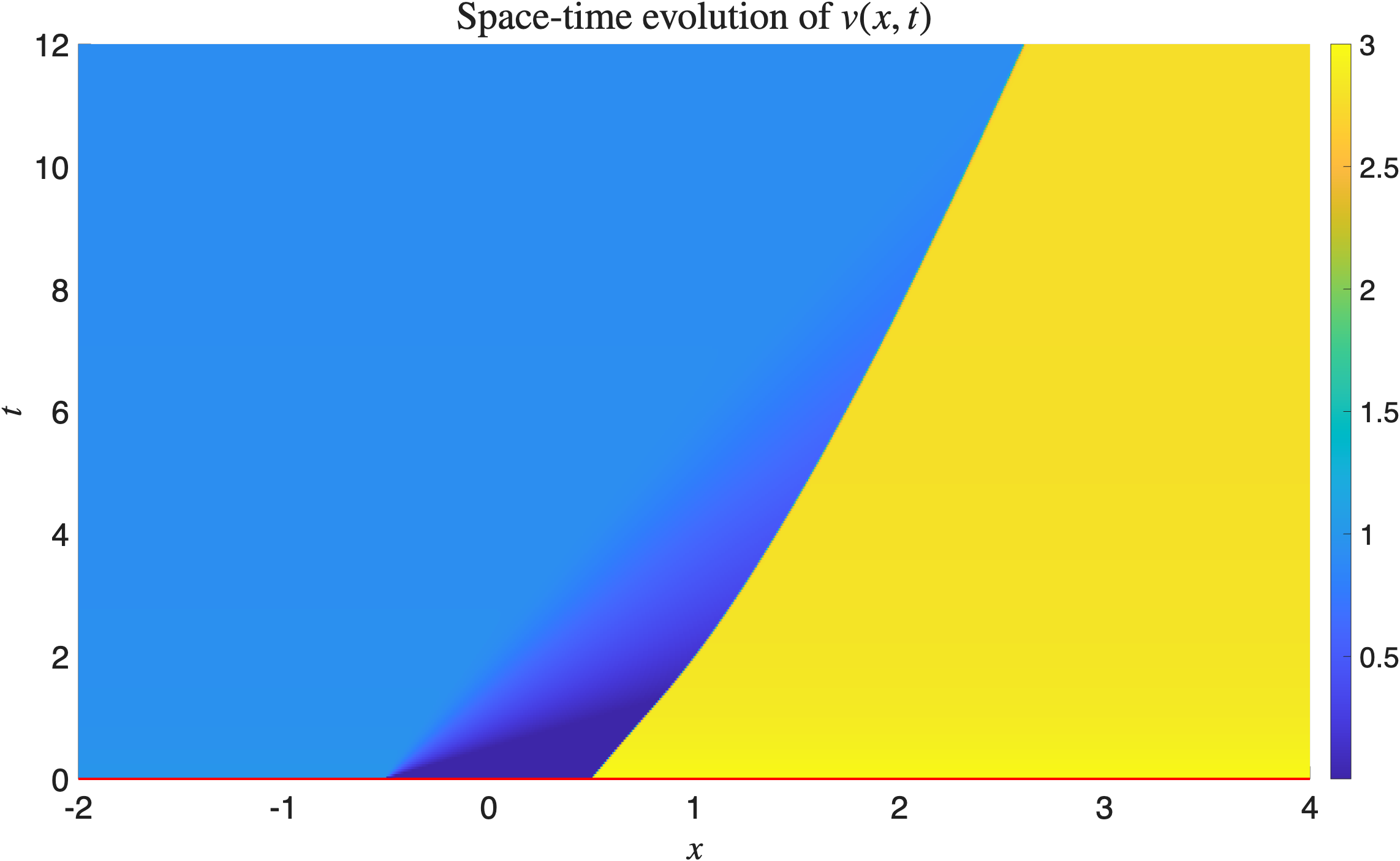}
    \includegraphics[width=0.45\linewidth]{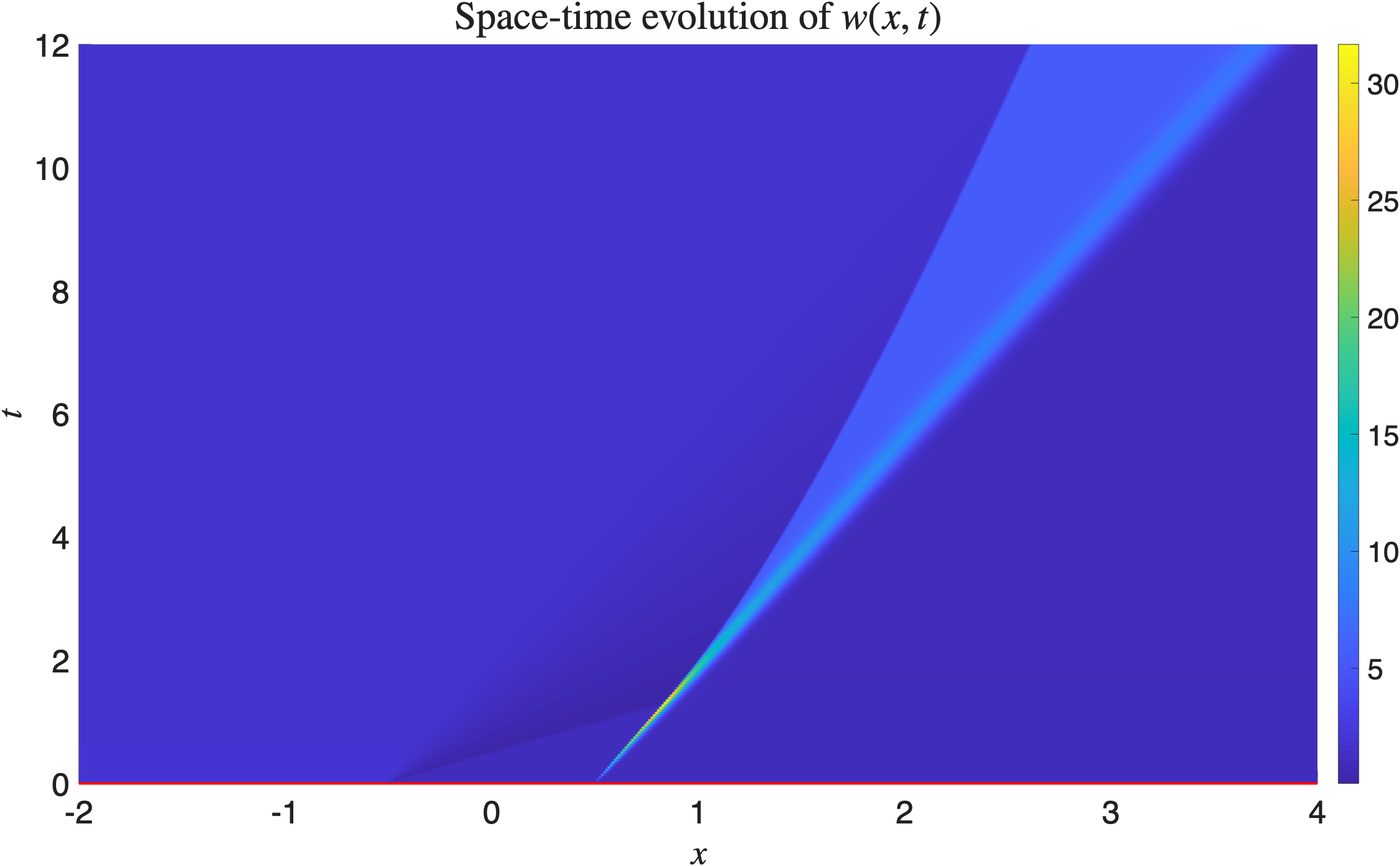}
    \caption{Space-time evolution of $v(x,t)$ (left) and $w(x,t)$ (right) for Case~5, Subcase~5.1 ($v_\thicksim = 0$, $v_- < v_+$).}
    \label{fig:case5Sub51}
\end{figure}

The numerical solution using the Lax-Friedrichs type scheme is shown in Figure~\ref{fig:case5Sub51}. In the left figure, the dark blue spot near $(x_{*1}, t_{*1}) \approx (0.8333, 1.2919)$ corresponds to the $v_\thicksim = 0$ state from the delta shock $\delta S$, while the curved boundary between $t_{*1}$ and $t_{*2} \approx 11.2058$ corresponds to the crossing of $S_1$ through the rarefaction $R_1^{(-)}$. The sharp straight boundary for $t > t_{*2}$ is the shock $S_1^{(*2)}$. In the right figure, the highly concentrated peak at $(x_{*1}, t_{*1})$ corresponds to the large strength of the delta contact discontinuity $\delta J_2$ at the moment of interaction, which then propagates as the cyan diagonal band with speed $\frac{1}{4h(t)}$ for $t > t_{*1}$.

{\bf Case 5: $v_\thicksim=0$, $v_->v_+$.} (Subcase 5.2)

In this subcase, Table~\ref{tb7} contains the initial data and the times and positions of the wave interactions.
The local Riemann solution at $x = -\frac{1}{2}$ consists of a composite wave $R_1^{(-)}J_2^{(*1)}$ with front speed $\frac{1}{h(t)}$, and the local Riemann solution at $x = \frac{1}{2}$ is a delta shock wave $\delta S$ with speed $\frac{1}{(1+v_+)h(t)} = \frac{1}{2h(t)}$.
The first interaction at $(x_{*1}, t_{*1})$ occurs when the composite wave $R_1^{(-)}J_2^{(*1)}$ meets the delta shock wave $\delta S$. Since $v_- \neq 0$, the delta shock wave $\delta S$ splits into a 1-shock wave $S_1$ and a delta contact discontinuity $\delta J_2$, with intermediate state $(V_{*2}, W_{*2}) = \left(1, \frac{2}{3}\right)$ and strength $\alpha(t_{*1}) = h(t_{*1})$. The delta contact discontinuity $\delta J_2$ propagates from $(x_{*1}, t_{*1})$ with speed $\frac{1}{(1+v_+)h(t)} = \frac{1}{2h(t)}$, while $S_1$ begins to cross the rarefaction wave $R_1^{(-)}$. However, since $v_- = 3 > v_+ = 1$, the shock $S_1$ cannot fully cross $R_1^{(-)}$, and the curve $x + \frac{1}{2} = \frac{1}{(1+v_+)^2}\int_0^t \exp \left(0.0880\frac{s}{1+s}\right)ds = \frac{1}{4}\int_0^t \exp \left(0.0880\frac{s}{1+s}\right)ds$ becomes its asymptote. Therefore, for sufficiently large $t > t_{*1}$, a residual rarefaction wave $R_1^{(-)}$ remains, and no further interactions occur. 

\begin{table}[h]
\centering
\begin{tabular}{|c|c|c|} \hline
Initial data & Times of interaction & Space of interaction\\ \hline
$v_- = 3$, $v_\thicksim = 0$, $v_+ = 1$ & $t_{*1} \approx 1.9234$ 
& $x_{*1} = 1.5000$ \\ 
$w_- = 2$, $w_\thicksim = 1$, $w_+ = 3$ & & \\ \hline
\end{tabular}
\caption{Initial data and time and position of interaction for Case~5, Subcase~5.2 ($v_\thicksim = 0$, $v_- > v_+$).}
\label{tb7}
\end{table}

\begin{figure}[h]
    \centering
    \includegraphics[width=0.45\linewidth]{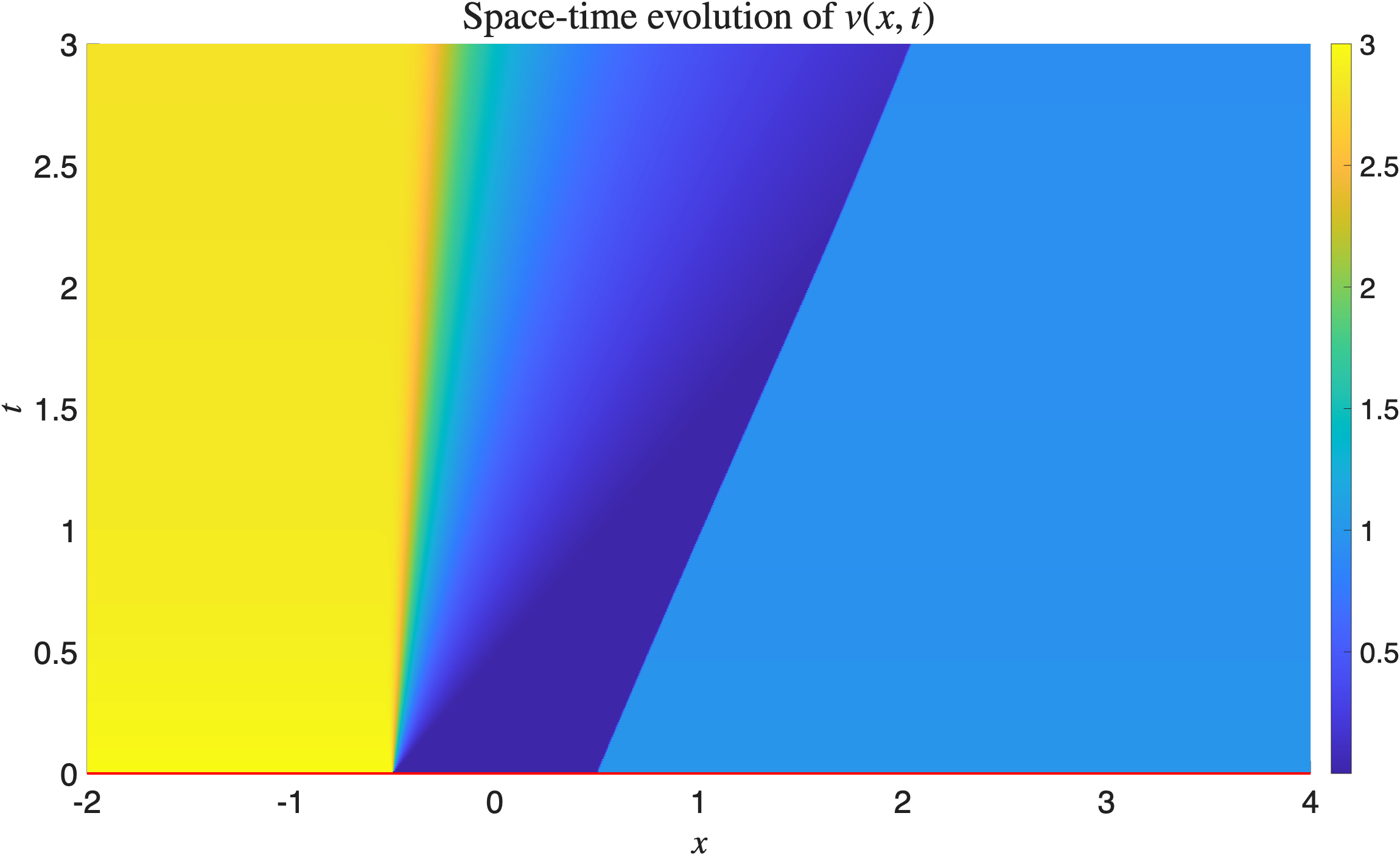}
    \includegraphics[width=0.45\linewidth]{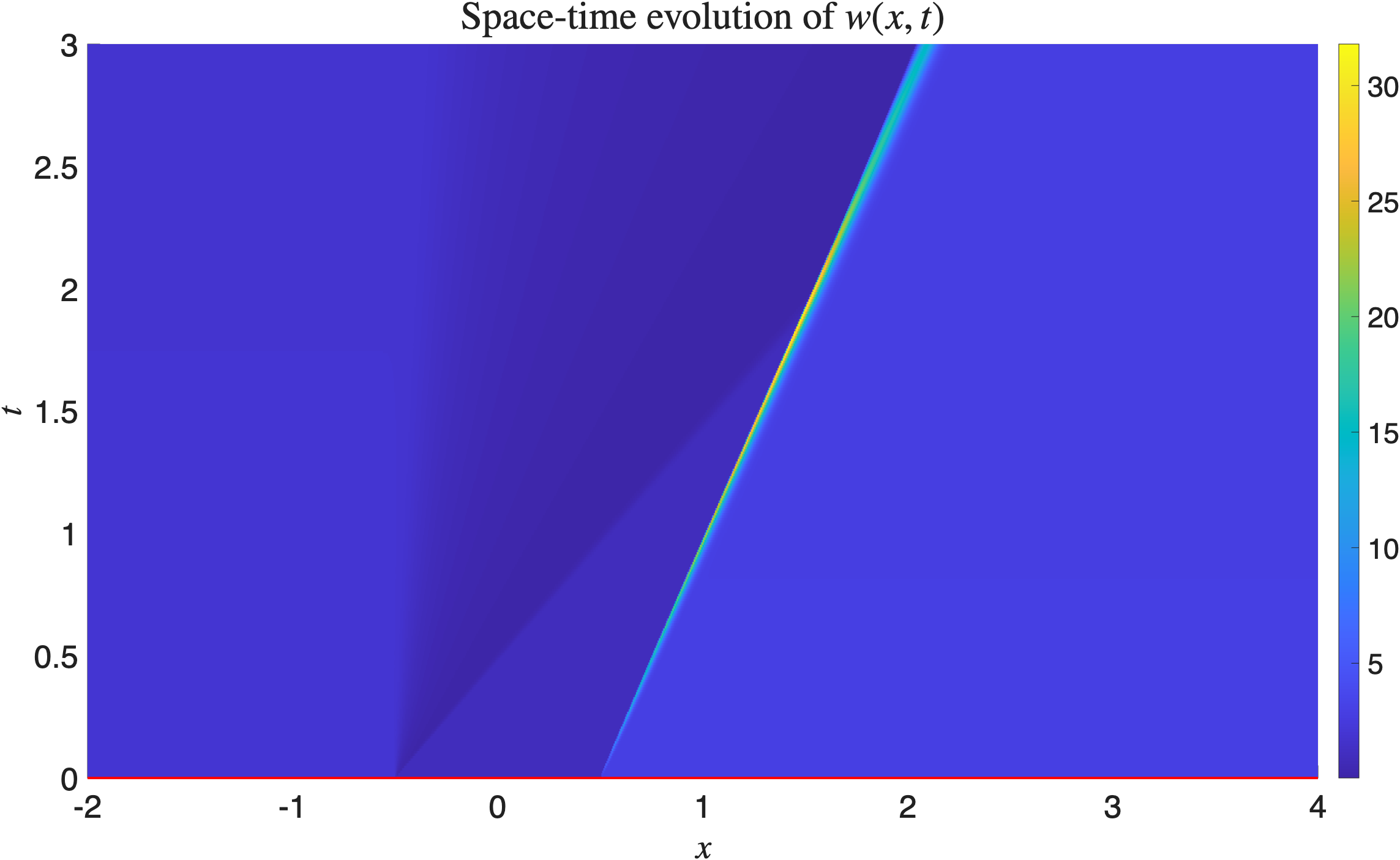}
    \caption{Space-time evolution of $v(x,t)$ (left) and $w(x,t)$ (right) for Case~5, Subcase~5.2 ($v_\thicksim = 0$, $v_- > v_+$), with $T = 3$.}
    \label{fig:case5Sub52}
\end{figure}

The numerical solution using the Lax-Friedrichs type scheme is shown in Figure~\ref{fig:case5Sub52}. In the left figure, the smooth transition from $v=3$ to $v\approx 0$ near the origin corresponds to the rarefaction $R_1^{(-)}$ and the $v_\thicksim=0$ state of the delta shock $\delta S$. The smooth left boundary persisting for all $t > t_{*1} \approx 1.9234$ is consistent with the fact that $S_1$ cannot fully cross $R_1^{(-)}$ since $v_+ = 1 < v_- = 3$, and illustrates the presence of the residual rarefaction $R_1^{(-)}$, in agreement with the limiting solution $R_1 + J_2$. In the right figure, the bright sharp line corresponds to the delta contact discontinuity $\delta J_2$ propagating from $(x_{*1}, t_{*1})$ with speed $\frac{1}{2h(t)}$.

{\bf Case 6: $v_-=0$, $v_\thicksim<v_+$.}

In this subcase, Table~\ref{tb8} contains the initial data and the times and positions of the wave interactions.
Since $v_- = 0$, the local Riemann solution at $x = -\frac{1}{2}$ is a delta shock wave $\delta S$ with speed $\frac{1}{(1+v_\thicksim)h(t)} = \frac{1}{2h(t)}$ and strength $\alpha(t) = h(t)\frac{w_- v_\thicksim}{1+v_\thicksim}\int_0^t \exp \left(0.0880\frac{s}{1+s}\right)ds$. The local Riemann solution at $x = \frac{1}{2}$ consists of a 1-shock wave $S_1^{(\sim)}$ and a 2-contact discontinuity $J_2^{(*1)}$, with intermediate state $(V_{*1}, W_{*1}) = (3, 3)$, and speeds $\frac{1}{8h(t)}$ and $\frac{1}{4h(t)}$, respectively.
Since the speed of $\delta S$ is $\frac{1}{2h(t)} > \frac{1}{8h(t)}$, the delta shock wave $\delta S$ catches the shock $S_1^{(\sim)}$ at $(x_{*1}, t_{*1})$. Since $v_- = 0$, the interaction produces a new delta shock wave $\delta S_1$ with speed $\frac{1}{(1+v_+)h(t)} = \frac{1}{4h(t)}$ and strength
\begin{equation*}
\alpha_1(t) = \frac{4}{3}\,h(t) + \frac{w_- v_+}{1+v_+}\,h(t) \int_{t_{*1}}^t \exp \left(0.0880\frac{s}{1+s}\right)ds.
\end{equation*}
Since the speed of $\delta S_1$ coincides with the speed of $J_2^{(*1)}$, namely $\frac{1}{(1+V_{*1})h(t)} = \frac{1}{4h(t)}$, they are parallel and will never interact. Therefore, no further interactions occur.

\begin{table}[h]
\centering
\begin{tabular}{|c|c|c|} \hline
Initial data & Times of interaction & Space of interaction\\ \hline
$v_- = 0$, $v_\thicksim = 1$, $v_+ = 3$ & $t_{*1} \approx 2.5508$ 
& $x_{*1} \approx 0.8333$ \\ 
$w_- = 2$, $w_\thicksim = 1$, $w_+ = 1$ & & \\ \hline
\end{tabular}
\caption{Initial data and time and position of interaction for Case~6 ($v_- = 0$, $v_\thicksim < v_+$).}
\label{tb8}
\end{table}

\begin{figure}[h]
    \centering
    \includegraphics[width=0.45\linewidth]{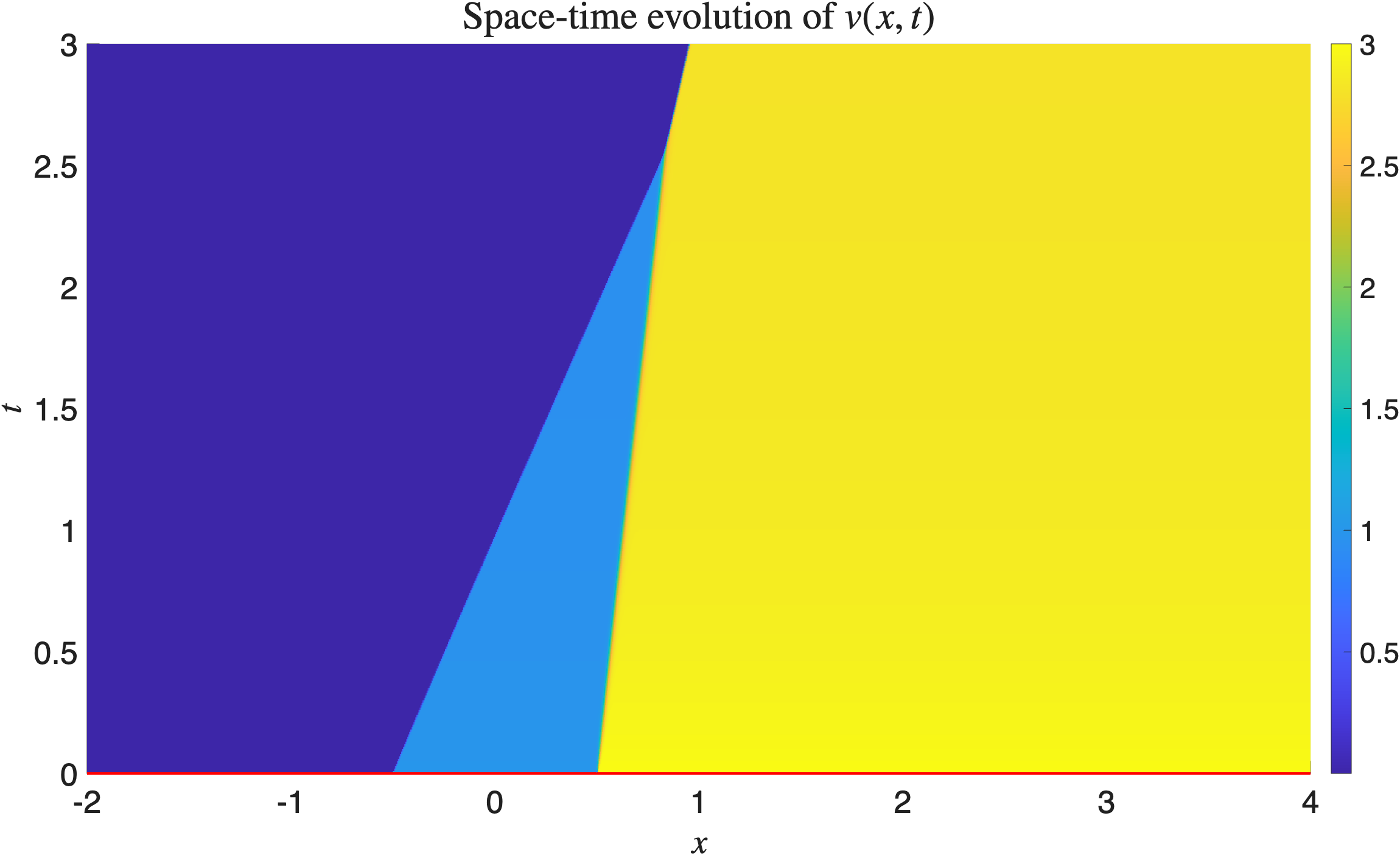}
    \includegraphics[width=0.45\linewidth]{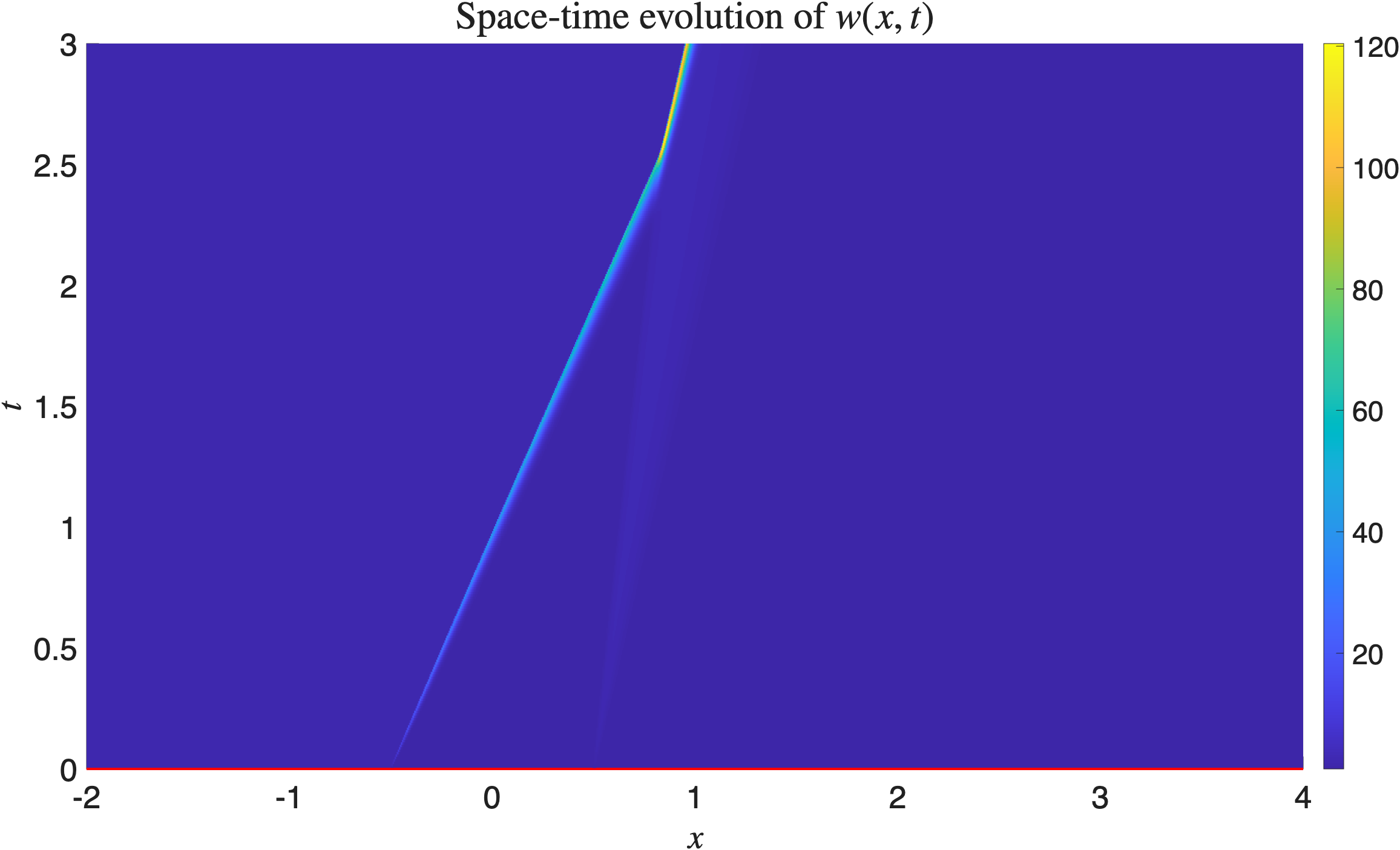}
    \caption{Space-time evolution of $v(x,t)$ (left) and $w(x,t)$ (right) for Case~6 ($v_- = 0$, $v_\thicksim < v_+$), with $T = 3$.}
    \label{fig:case6}
\end{figure}

The numerical solution using the Lax-Friedrichs type scheme is shown in Figure~\ref{fig:case6}. In the left figure, the dark blue region corresponds to the left state $v_-=0$, the cyan region to the intermediate state $V_{*1}=v_\thicksim=1$, and the yellow region to the right state $v_+=3$. The interaction at $(x_{*1},t_{*1})\approx(0.8333, 2.5508)$ is clearly visible as the change in slope of the left boundary, where the delta shock $\delta S$ meets $S_1^{(\sim)}$ and produces the new delta shock $\delta S_1$ propagating with speed $\frac{1}{4h(t)}$, parallel to $J_2^{(*1)}$. In the right figure, the bright sharp line corresponds to the delta shock carrying the concentrated $w$ mass, whose strength grows over time according to \eqref{sol_case6}.

{\bf Case 7: $v_-=0$, $v_+<v_\thicksim$.}

In this subcase, Table~\ref{tb9} contains the initial data and the times and positions of the wave interactions.
Since $v_- = 0$, the local Riemann solution at $x = -\frac{1}{2}$ is a delta shock wave $\delta S$ with speed $\frac{1}{(1+v_\thicksim)h(t)} = \frac{1}{4h(t)}$. The local Riemann solution at $x = \frac{1}{2}$ consists of a 1-rarefaction wave $R_1^{(\sim)}$ and a 2-contact discontinuity $J_2^{(*1)}$, with intermediate state $(V_{*1}, W_{*1}) = \left(1, \frac{1}{3}\right)$ and speeds $\frac{1}{4h(t)}$ (front of $R_1^{(\sim)}$) and $\frac{1}{2h(t)}$ ($J_2^{(*1)}$), respectively.
Since the speed of $\delta S$ equals the speed of the front of $R_1^{(\sim)}$, namely $\frac{1}{4h(t)}$, the delta shock wave $\delta S$ catches the tail of $R_1^{(\sim)}$ at $(x_{*1}, t_{*1})$ and begins to cross it as a new delta shock wave $\delta S_1$. The crossing curve $\delta S_1$ is given by
\begin{equation*}
x = \frac{1}{2} + \left(\sqrt{\int_0^t \exp \left(0.0880 \frac{s}{1+s}\right)ds} - \sqrt{3}\right)^2,
\end{equation*}
with increasing speed. The second interaction at $(x_{*2}, t_{*2})$ occurs when $\delta S_1$ exits the front of $R_1^{(\sim)}$, producing a new delta shock wave $\delta S_2$ with speed $\frac{1}{(1+v_+)h(t)} = \frac{1}{2h(t)}$ and strength
\begin{equation*}
\alpha_2(t) = \frac{w_- v_+}{1+v_+}\,h(t)\int_{t_{*2}}^t \exp \left(0.0880\frac{s}{1+s}\right)ds + \frac{\alpha_1(t_{*2})}{h(t_{*2})}\,h(t),
\end{equation*}
together with a 2-contact discontinuity $J_2^{(*1)}$ with speed $\frac{1}{(1+V_{*1})h(t)} = \frac{1}{2h(t)}$. Since $\delta S_2$ and $J_2^{(*1)}$ propagate with the same speed, they are parallel and will never interact. Therefore, no further interactions occur. 

\begin{table}[h]
\centering
\begin{tabular}{|c|c|c|} \hline
Initial data & Times of interaction & Space of interaction\\ \hline
$v_- = 0$, $v_\thicksim = 3$, $v_+ = 1$ & $t_{*1} \approx 5.0391$ 
& $x_{*1} \approx 0.8333$ \\ 
$w_- = 2$, $w_\thicksim = 1$, $w_+ = 3$ & $t_{*2} \approx 11.2058$ 
& $x_{*2} = 3.5000$ \\ \hline
\end{tabular}
\caption{Initial data, times and positions of interactions for 
Case~7 ($v_- = 0$, $v_+ < v_\thicksim$).}
\label{tb9}
\end{table}

\begin{figure}[h]
    \centering
    \includegraphics[width=0.45\linewidth]{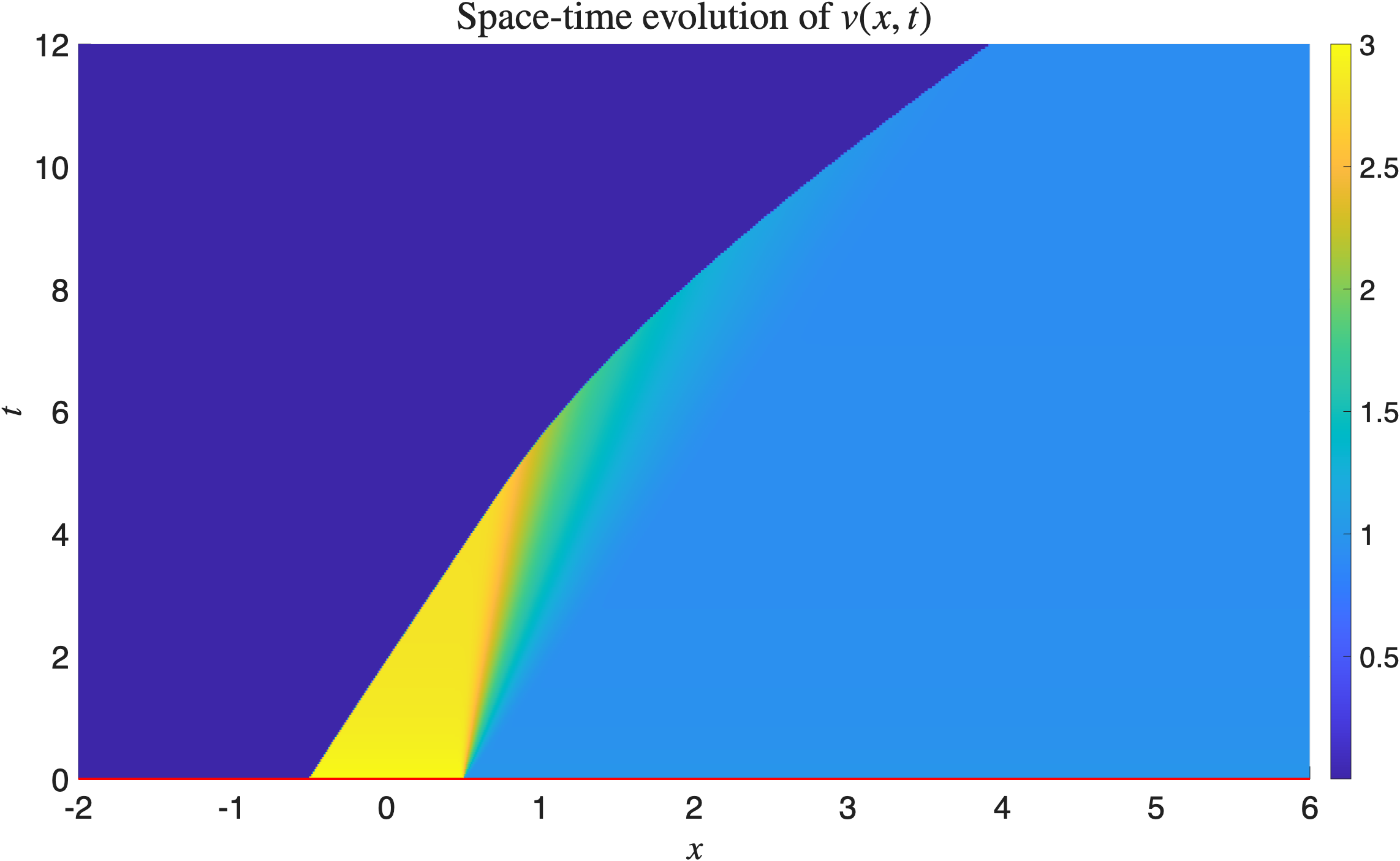}
    \includegraphics[width=0.45\linewidth]{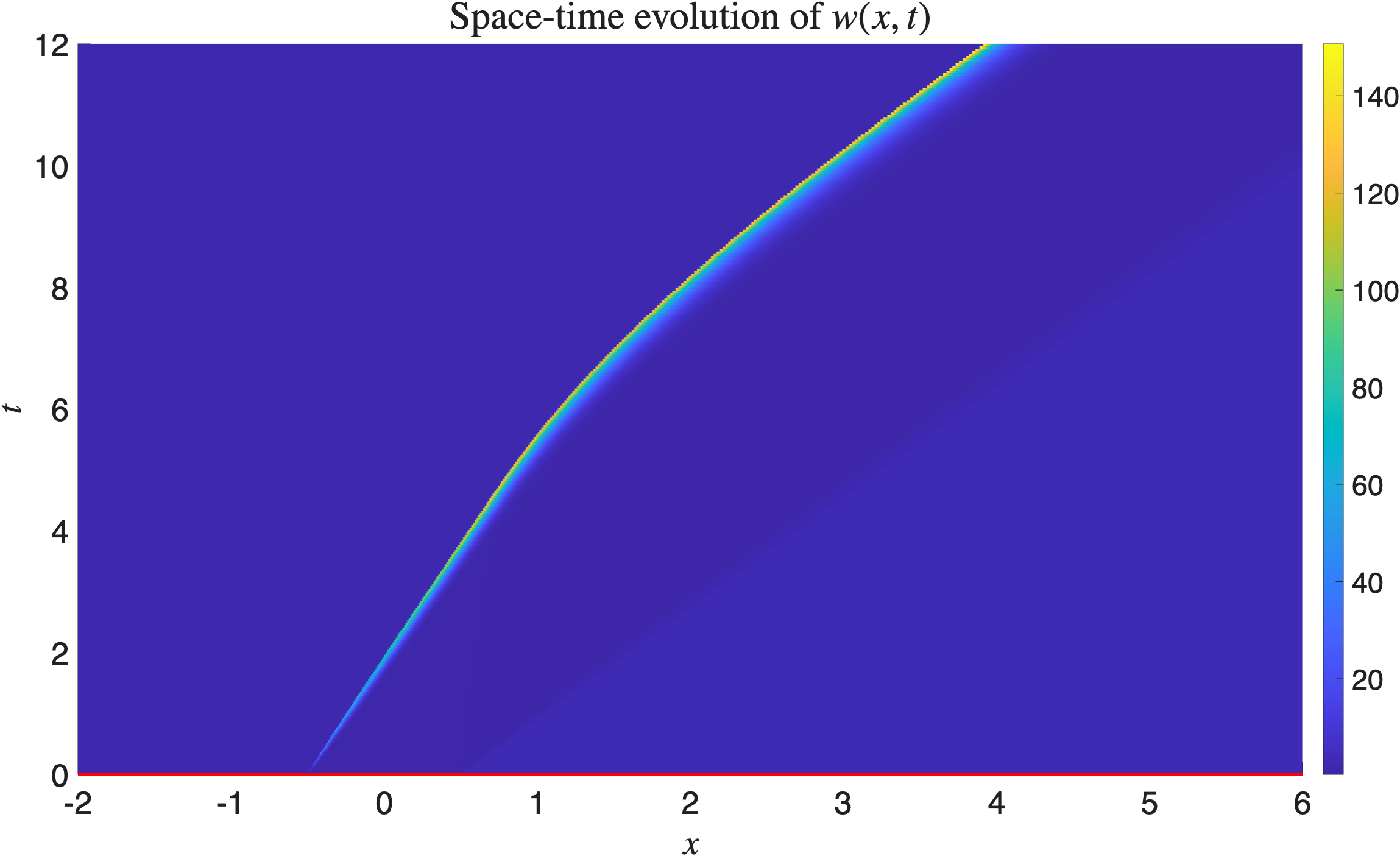}
    \caption{Space-time evolution of $v(x,t)$ (left) and $w(x,t)$ 
    (right) for Case~7 ($v_- = 0$, $v_+ < v_\thicksim$), 
    with $T = 12$.}
    \label{fig:case7}
\end{figure}

The numerical solution using the Lax-Friedrichs type scheme is shown in Figure~\ref{fig:case7}. In the left figure, the dark blue region corresponds to the left state $v_-=0$, and the yellow peak near the origin corresponds to the $v_\thicksim=3$ state of $R_1^{(\sim)}$ being crossed by the delta shock. The curved boundary between $t_{*1} \approx 5.0391$ and $t_{*2} \approx 11.2058$ corresponds to the delta shock $\delta S_1$ crossing $R_1^{(\sim)}$ along the parabolic curve $x = \frac{1}{2} + \left(\sqrt{\int_0^t \exp \left(0.0880 \frac{s}{1+s}\right)ds} - \sqrt{3}\right)^2$, after which the boundary becomes the straight line $\delta S_2$ with speed $\frac{1}{2h(t)}$. In the right figure, the bright curved line clearly shows the parabolic trajectory of the delta shock $\delta S_1$ during the crossing of $R_1^{(\sim)}$, with growing strength $\alpha_1(t)$, followed by the straight trajectory of $\delta S_2$ for $t > t_{*2}$.

\subsection{Profile of numerical solution}

In this section, we show some profiles of the numerical solution. 
More specifically, we analyse the cases where $0 < v_- < v_\thicksim < v_+$ (Case~1) and $v_- = 0$, $v_+ < v_\thicksim$ (Case~7), which represent respectively a classical wave interaction converging to $S_1 + J_2$, and a nonclassical interaction involving a delta shock wave converging to $\delta S$.\\

{\bf Case 1: $0<v_-<v_\thicksim<v_+$.}

\begin{figure}[h]
    \centering
\begin{tabular}{cc}
\includegraphics[width=0.45\linewidth]{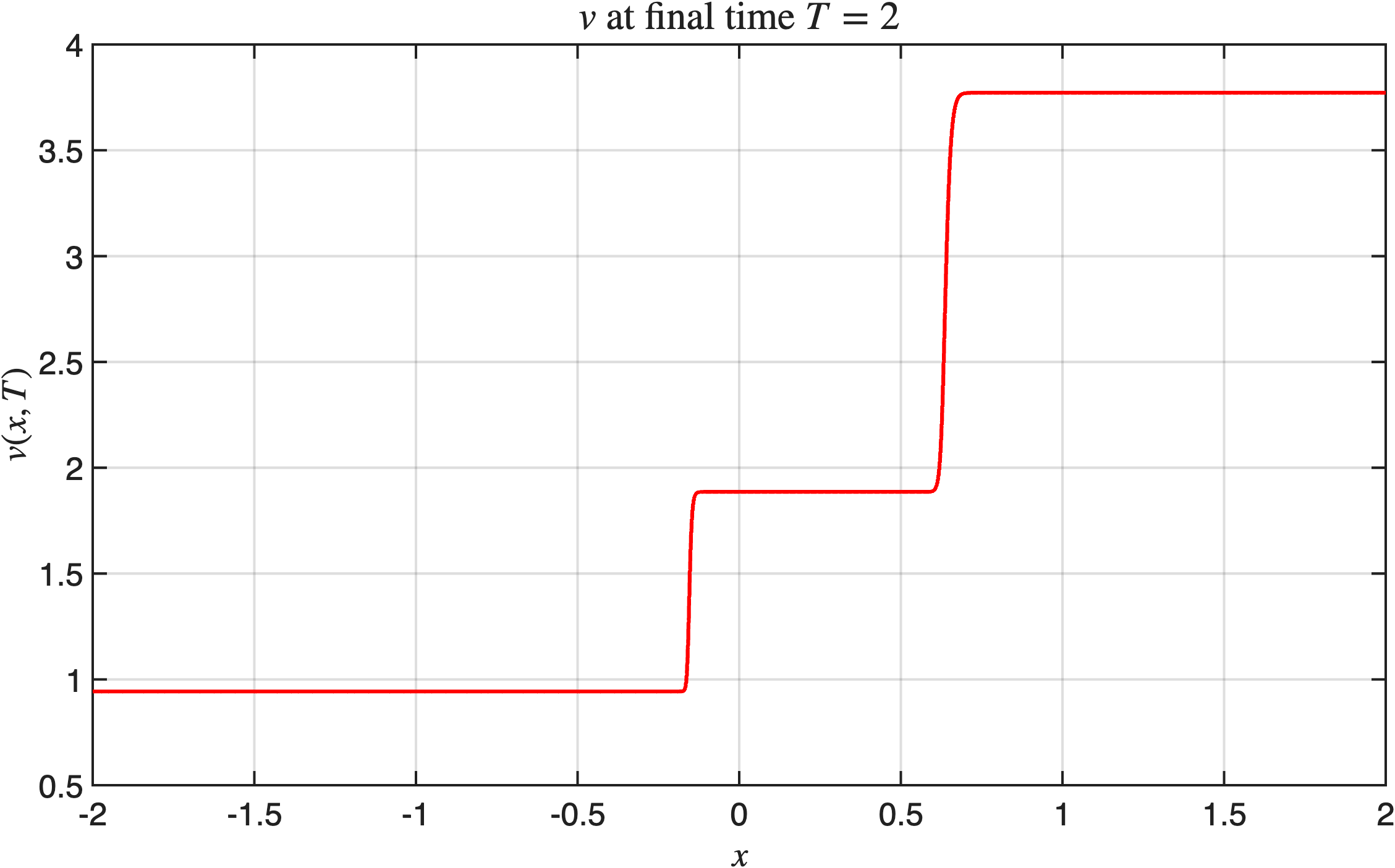} & 
\includegraphics[width=0.45\linewidth]{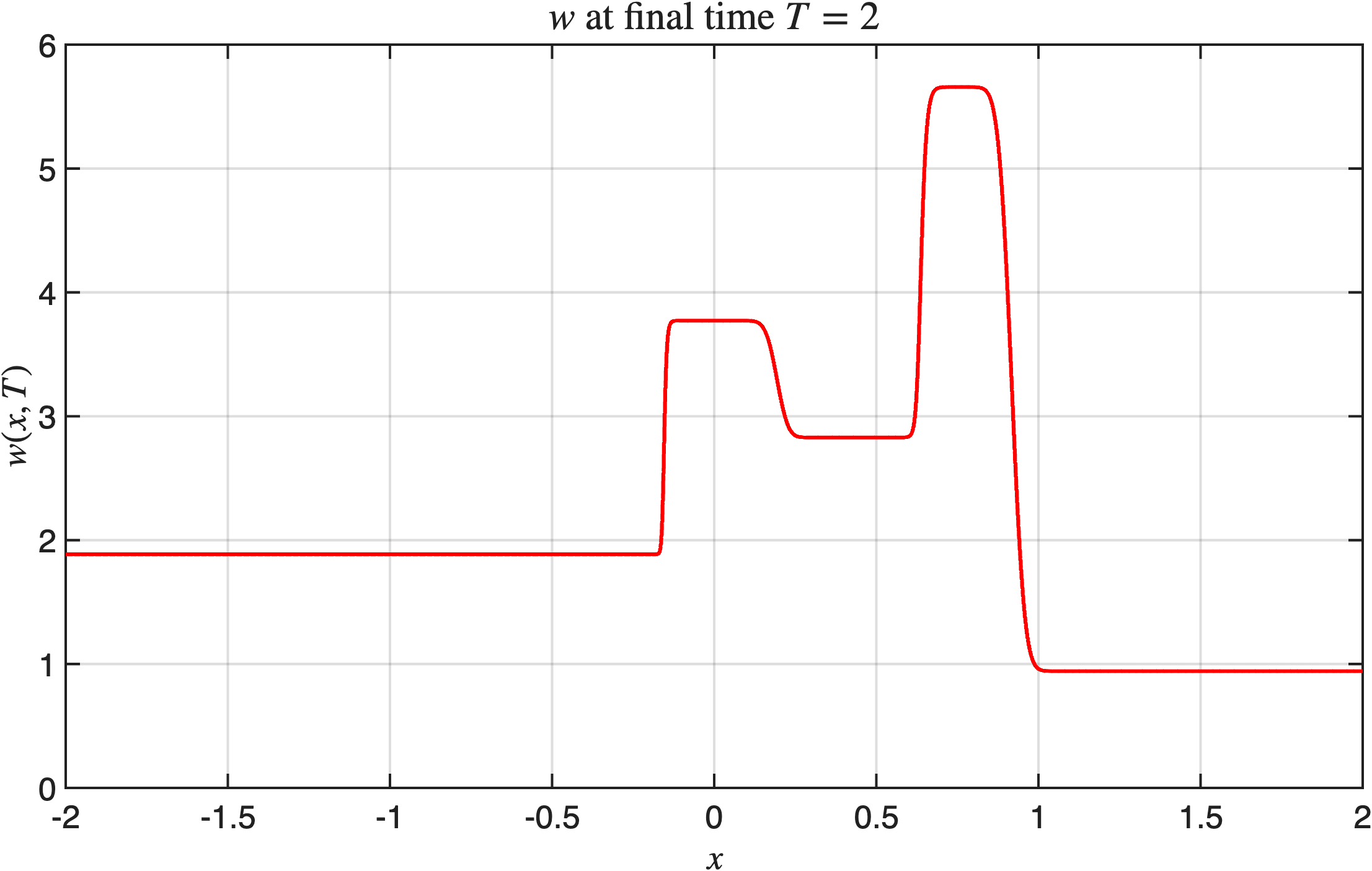} \\
(a) & (b) \\
\includegraphics[width=0.45\linewidth]{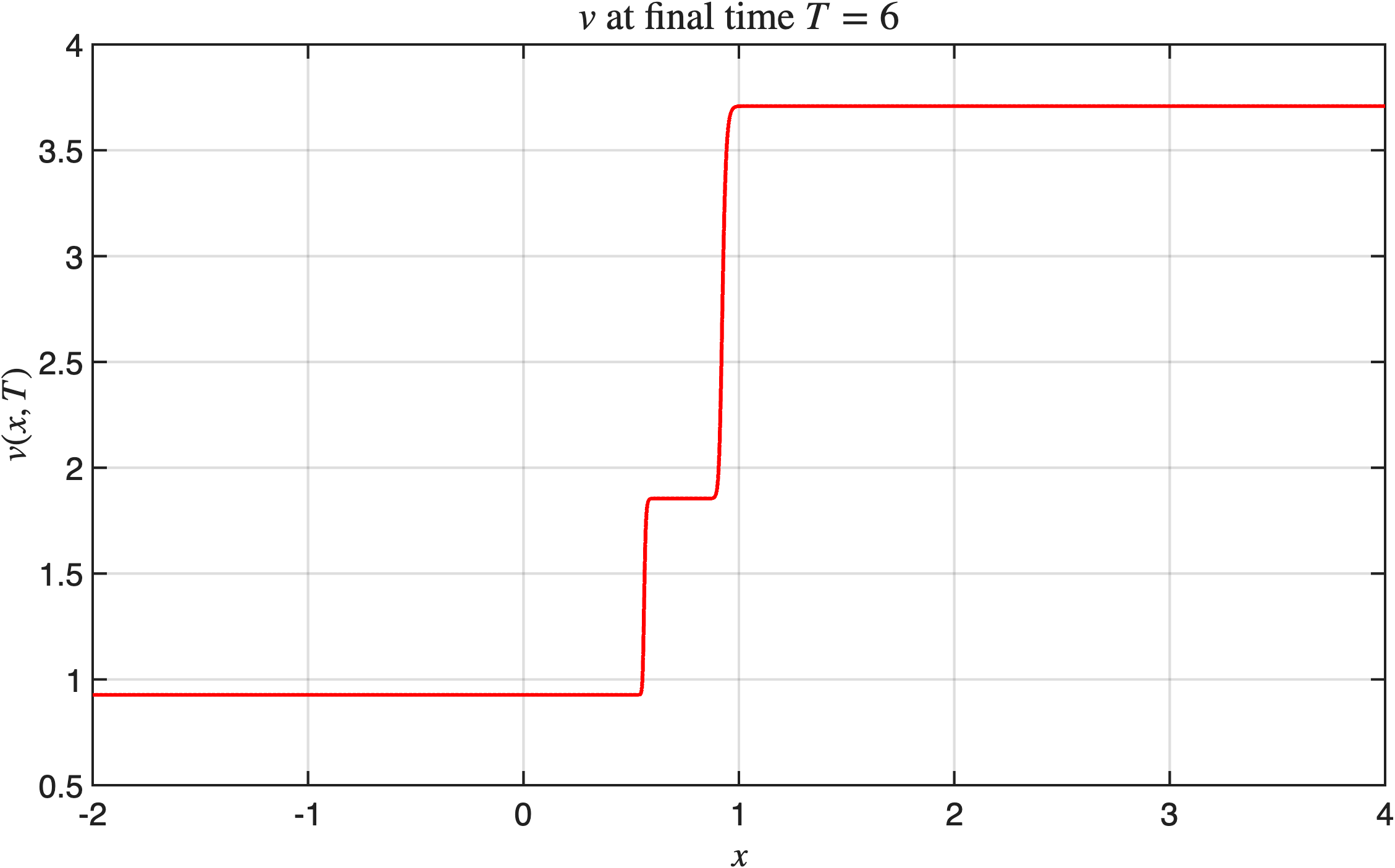} & 
\includegraphics[width=0.45\linewidth]{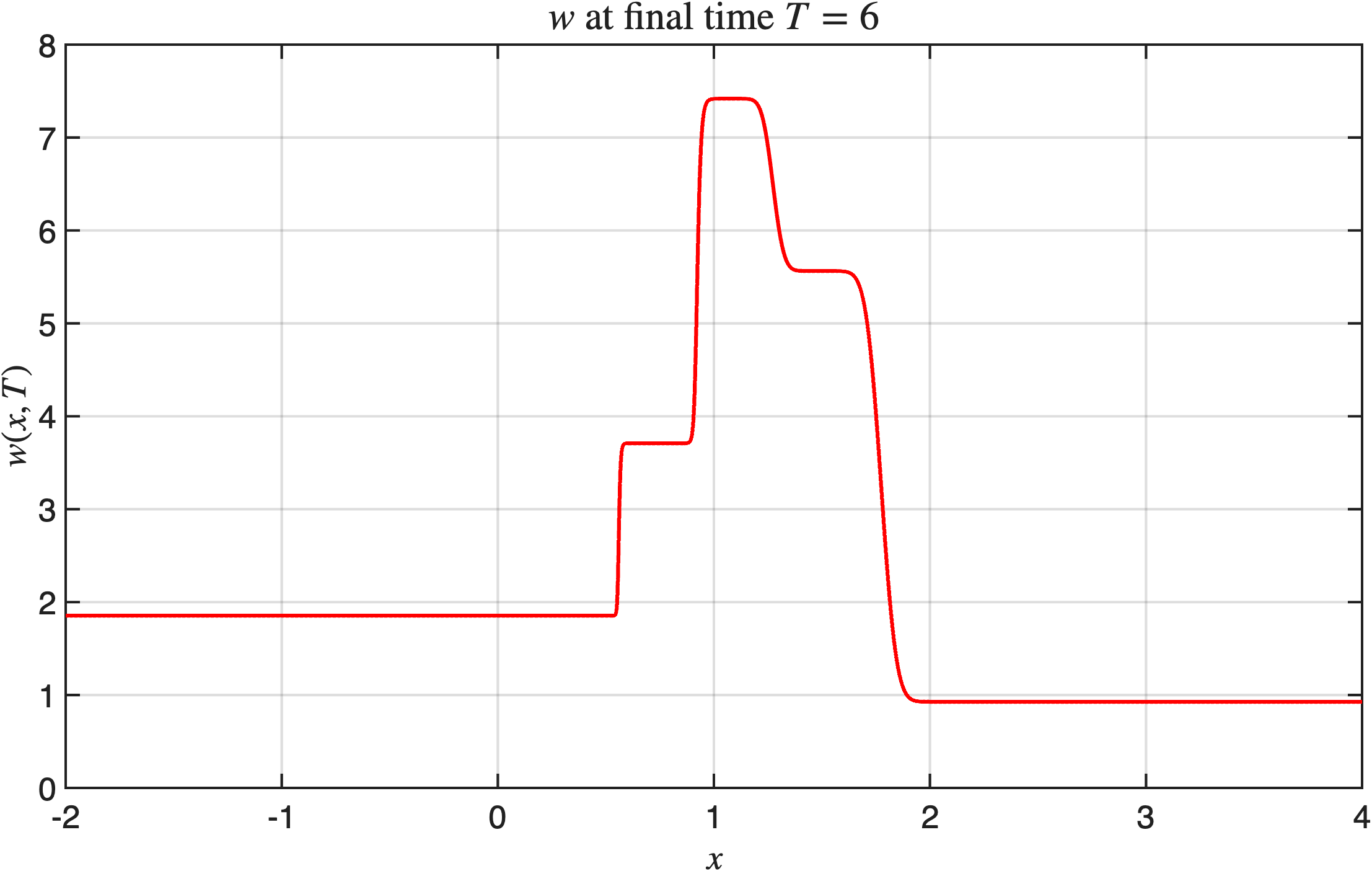} \\
(c) & (d) \\
\includegraphics[width=0.45\linewidth]{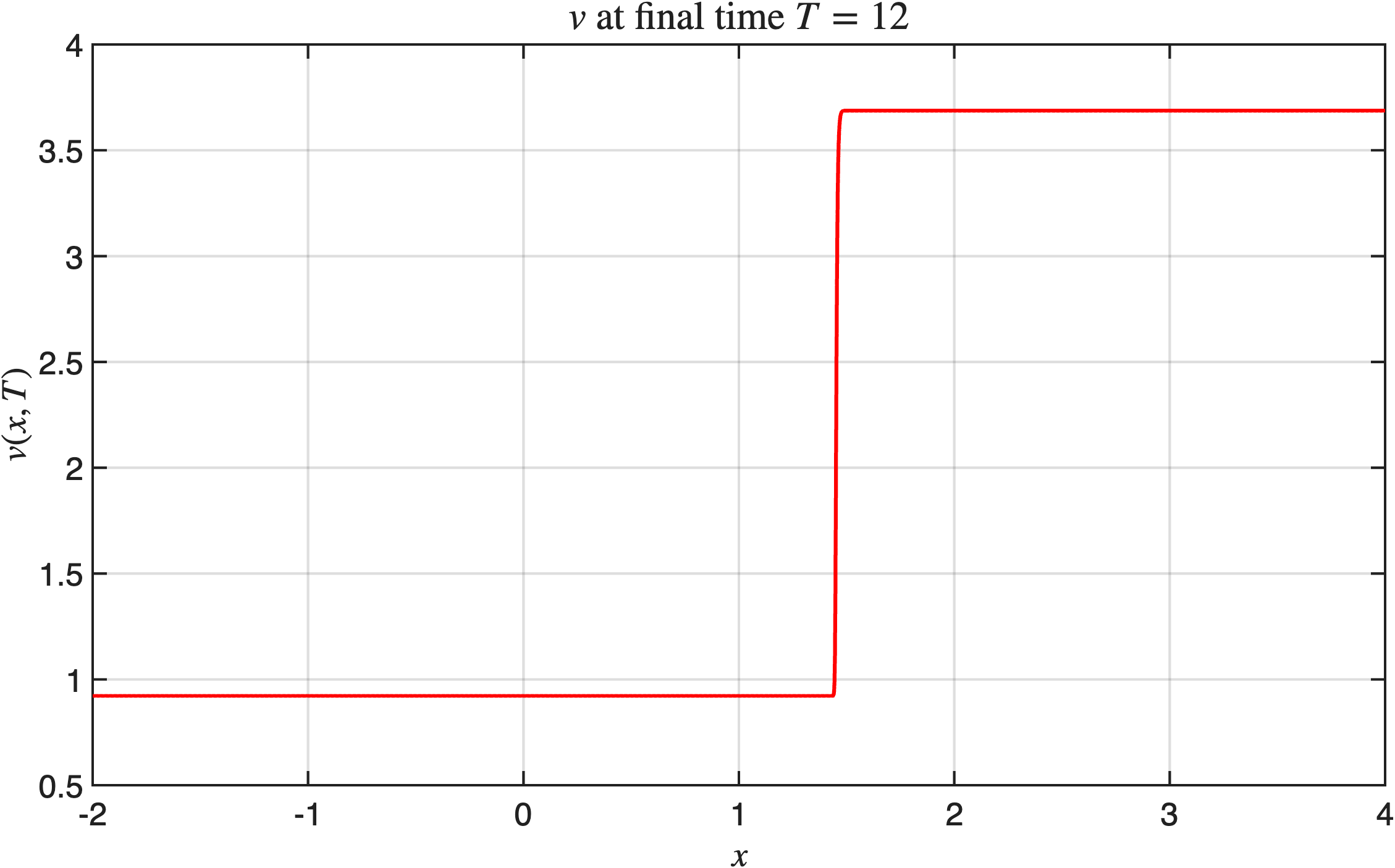} & 
\includegraphics[width=0.45\linewidth]{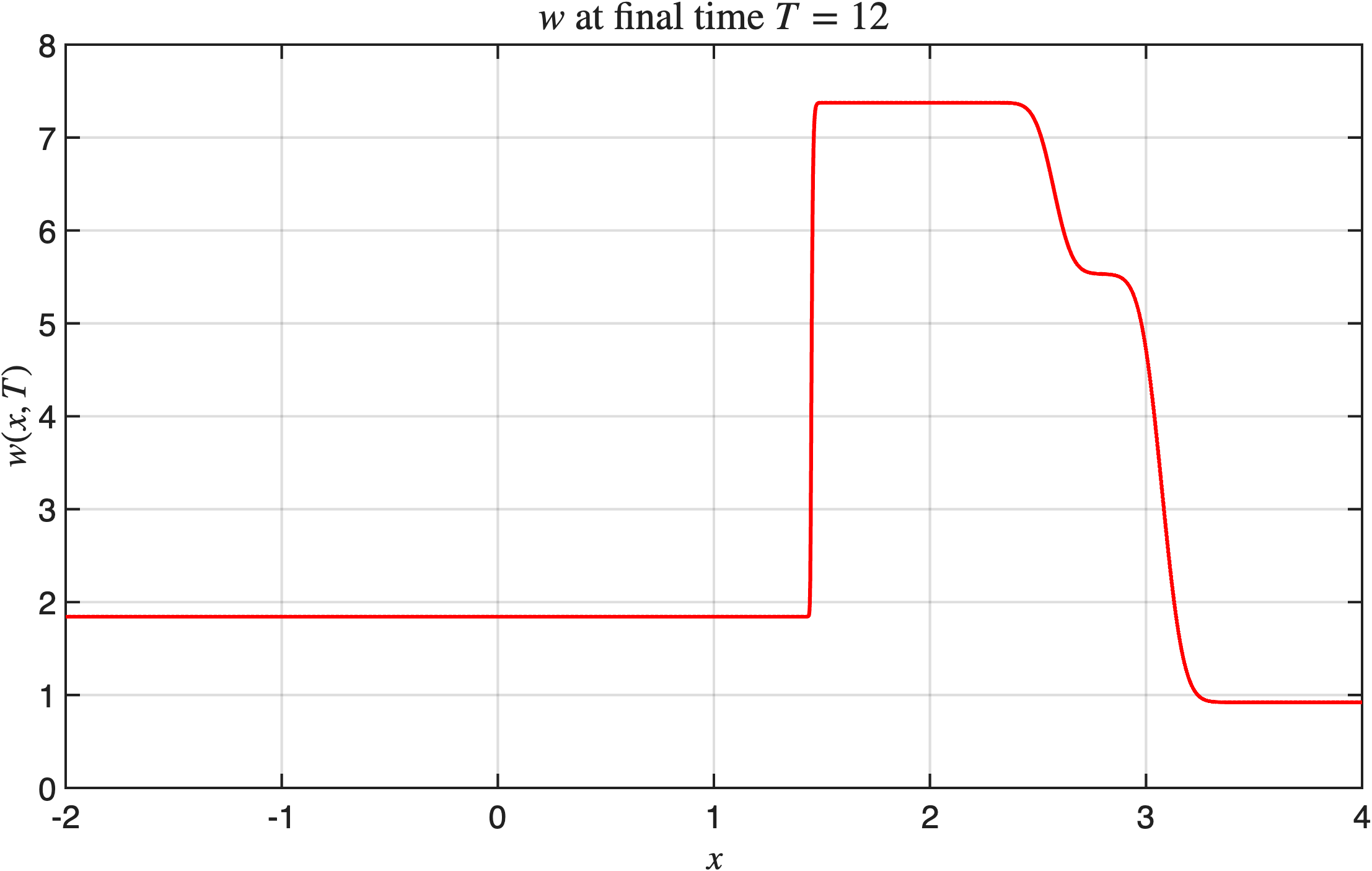} \\
(e) & (f)
\end{tabular}
    \caption{Profiles of $v(x,T)$ (left) and $w(x,T)$ (right) for Case~1 ($0 < v_- < v_\thicksim < v_+$) at times $T=2$ (a)--(b), $T=6$ (c)--(d), and $T=12$ (e)--(f).}
    \label{fig:ProfileCase1}
\end{figure}

Figures~\ref{fig:ProfileCase1}(a)--(b) show the profiles of $v(x,T)$ and $w(x,T)$ at $T=2$, which corresponds to the time interval $0 \leq t < t_{*1} \approx 3.5647$. The five-region structure of the solution is clearly visible: in (a), the solution takes the constant values $v_-h(2)=0.9430$, $V_{*1}h(2)=1.8860$ and $V_{*2}h(2)=3.7721$ in three distinct regions, separated by the shock $S_1^{(-)}$, the contact discontinuity $J_2^{(*1)}$ and the shock $S_1^{(\sim)}$. In (b), the solution takes the constant values $w_-h(2)=1.8860$, $W_{*1}h(2)=3.7721$, $w_\thicksim h(2)=2.8291$ and $W_{*2}h(2)=5.6581$ in four distinct regions.
Figures~\ref{fig:ProfileCase1}(c)--(d) show the profiles at $T=6$, which lies in the interval $t_{*1} \leq t < t_{*2} \approx 9.3599$. After the first interaction, the contact discontinuity $J_2^{(*1)}$ and the shock $S_1^{(\sim)}$ have merged, producing the new shock $S_1^{(*1)}$ and contact discontinuity $J_2^{(*3)}$. In (c), the intermediate region $V_{*1}h(6)=1.8547$ is now very narrow since $S_1^{(-)}$ and $S_1^{(*1)}$ are approaching each other. In (d), the intermediate state $W_{*3}h(6)=7.4188$ is visible as the peak between $S_1^{(*1)}$ and $J_2^{(*3)}$, while $W_{*2}h(6)=5.5641$ occupies the region between $J_2^{(*3)}$ and $J_2^{(*2)}$.
Figures~\ref{fig:ProfileCase1}(e)--(f) show the profiles at $T=12$, which lies in the interval $t \geq t_{*2} \approx 9.3599$. After the second interaction, the shocks $S_1^{(-)}$ and $S_1^{(*1)}$ have merged into the single shock $S_1^{(*3)}$. In (e), only one sharp jump remains, at $x \approx 1.4$, separating the states $v_-h(12)=0.9220$ and $V_{*3}h(12)=3.6879$. In (f), the two parallel contact discontinuities $J_2^{(*3)}$ and $J_2^{(*2)}$ are visible as the smooth transition between the states $W_{*3}h(12)=7.3758$ and $W_{*2}h(12)=5.5319$, and the right state $w_+h(12)=0.9220$.

\begin{figure}[h]
    \centering
    \includegraphics[width=0.45\linewidth]{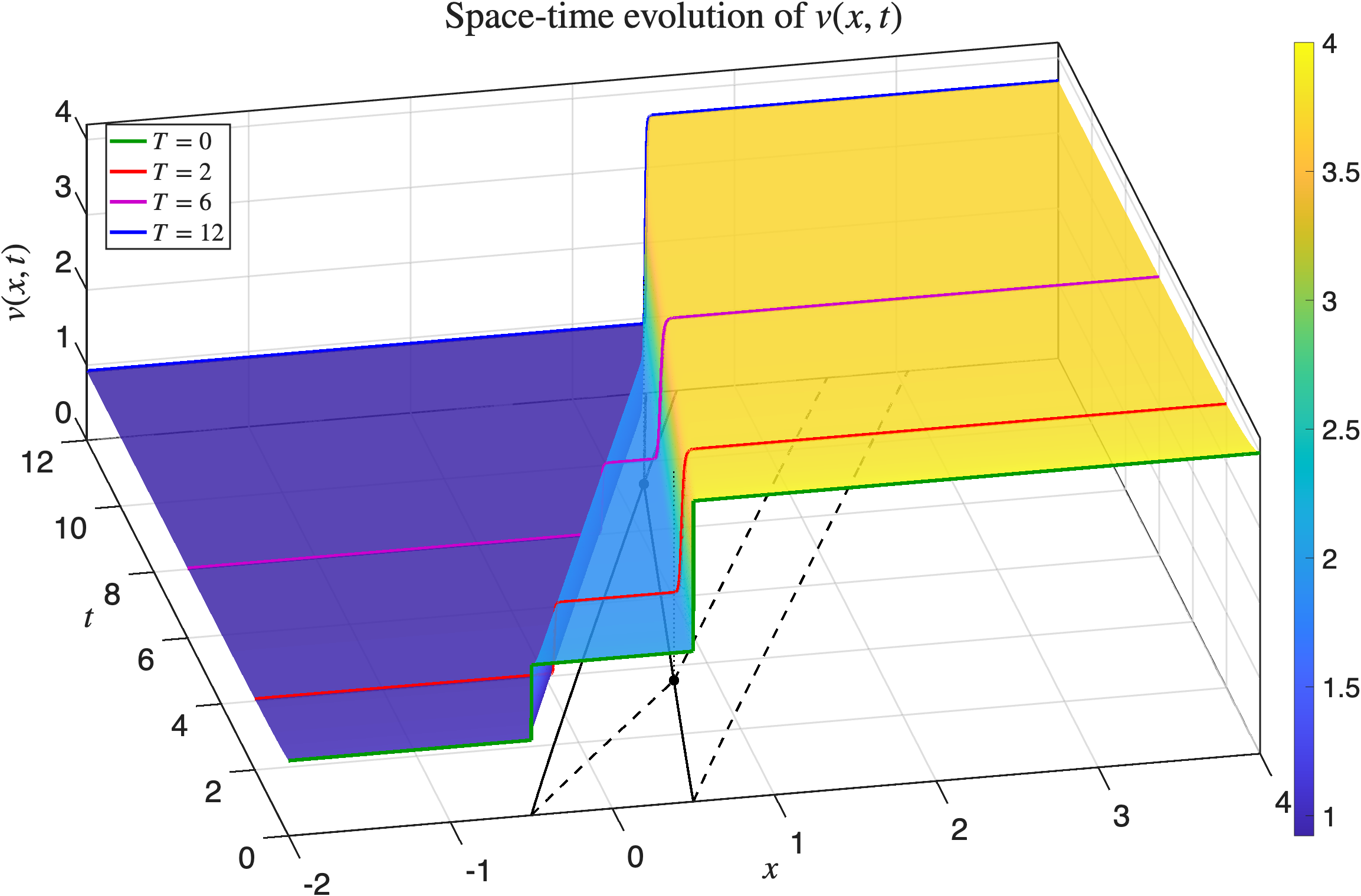}
    \includegraphics[width=0.45\linewidth]{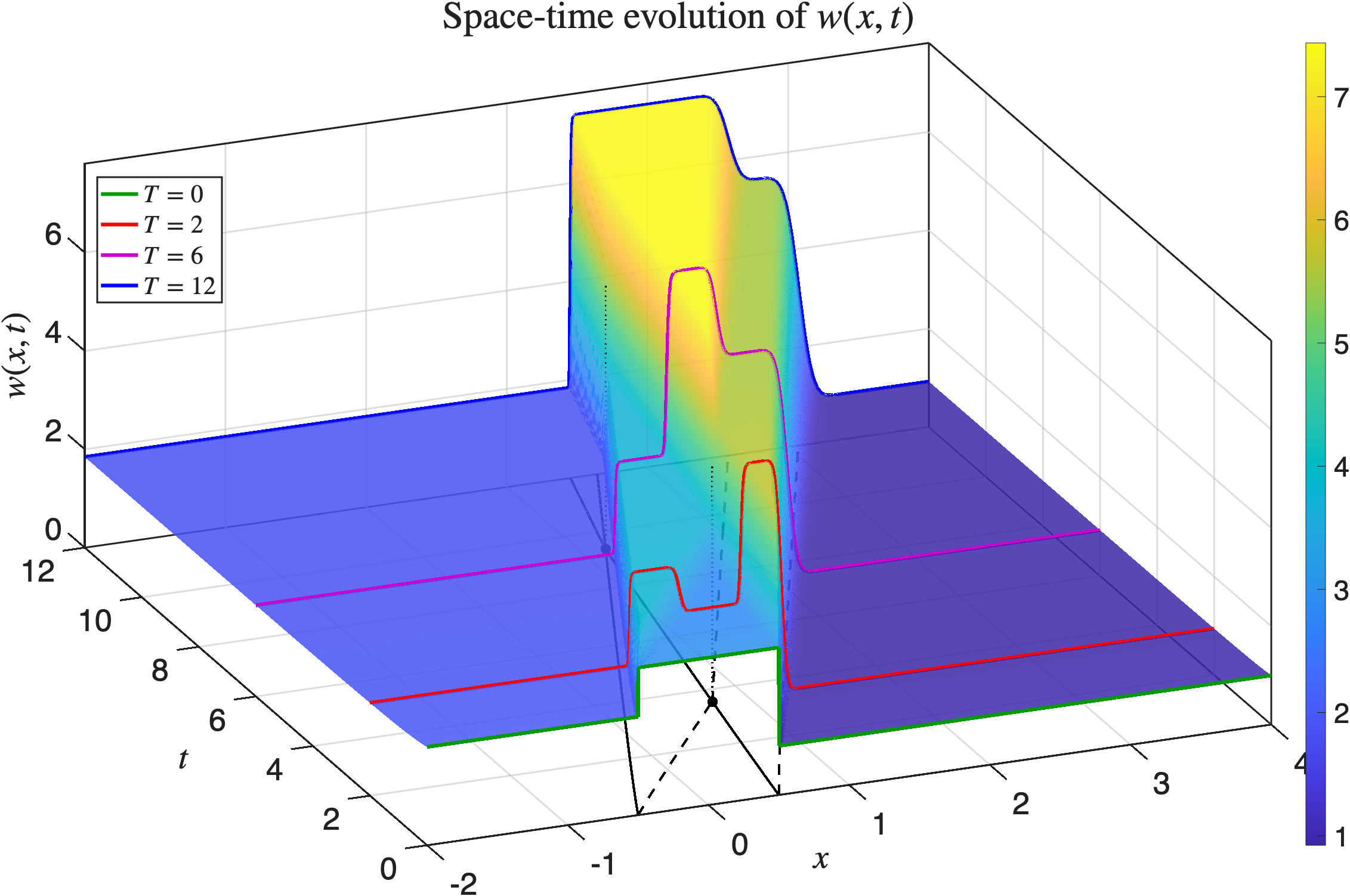}
    \caption{Space-time evolution of $v(x,t)$ (left) and 
    $w(x,t)$ (right) for Case~1 ($0 < v_- < v_\thicksim < v_+$), complementing Figure~\ref{fig:ProfileCase1}. The solid colored lines show the profiles at $T=0$, $T=2$, $T=6$, and $T=12$, and the black curves on the $x$-$t$ plane show the shock curves (solid) and contact discontinuities (dashed), with the interaction points marked as black dots.}
    \label{fig:3DCase1}
\end{figure}

Figure~\ref{fig:3DCase1} shows the space-time evolution of $v(x,t)$ and $w(x,t)$ for Case~I, complementing the profiles shown in Figure~\ref{fig:ProfileCase1}. The solid colored lines correspond to the profiles at times $T=0$, $T=2$, $T=6$ and $T=12$ already analyzed in Figure~\ref{fig:ProfileCase1}(a)--(f), while the wave curves $g_1(t;1,2)$, $g_1^*(t;2,4)$, $g_1^{**}(t;1,4)$ (shocks, solid) and $g_2(t;2)$, $g_2^*(t;4)$, $g_2(t;4)$ (contacts, dashed) are projected onto the $x$-$t$ plane, allowing the reader to identify precisely where each profile intersects the wave fronts. 
The interaction points $(x_{*1},t_{*1})\approx(0.7500, 3.5647)$ and $(x_{*2},t_{*2})\approx(1.1667, 9.3599)$ are marked with black dots on the $x$-$t$ plane.

{\bf Case 7: $v_-=0$, $v_+<v_\thicksim$.}

\begin{figure}[h]
    \centering
\begin{tabular}{cc}
\includegraphics[width=0.45\linewidth]{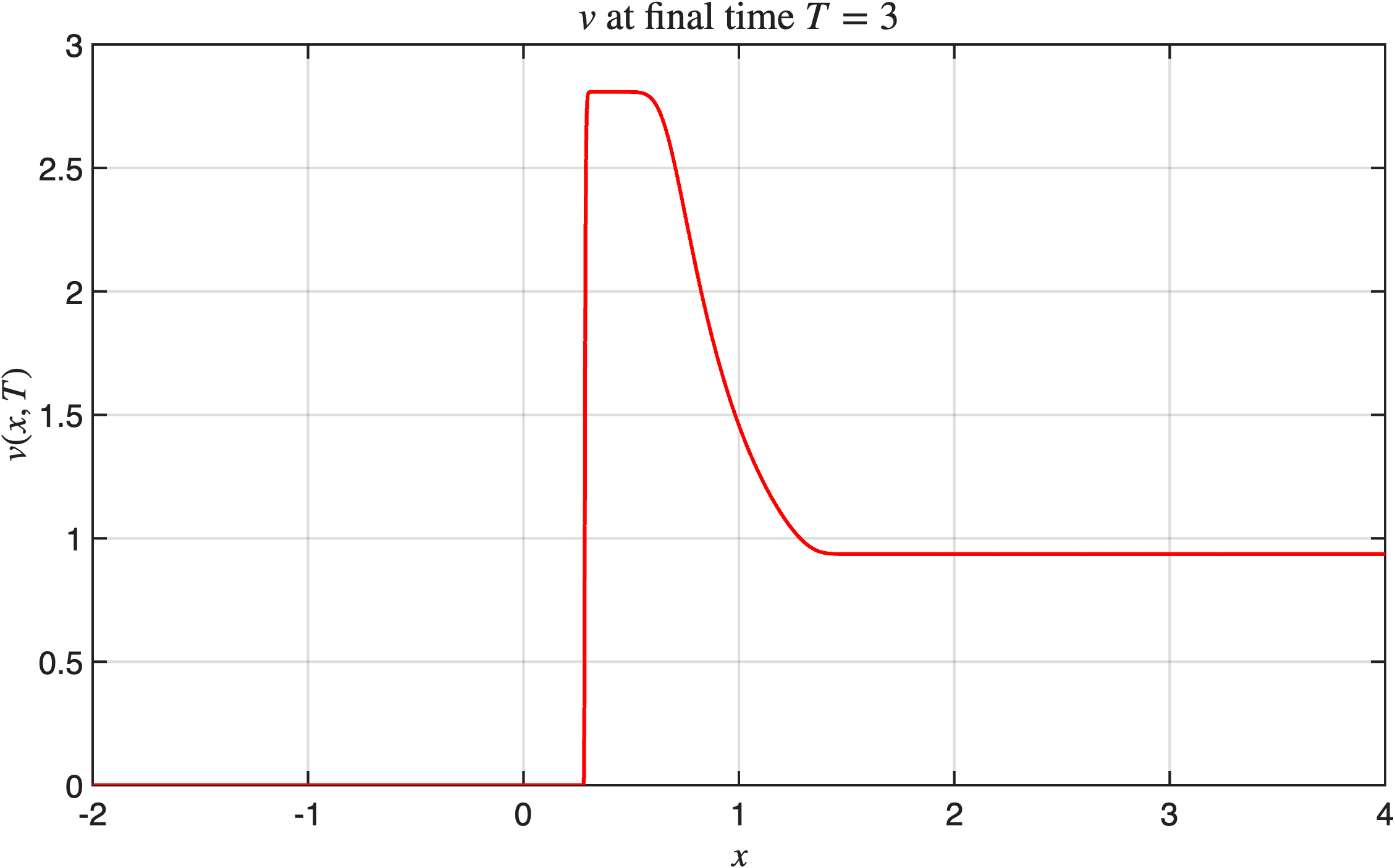} & 
\includegraphics[width=0.45\linewidth]{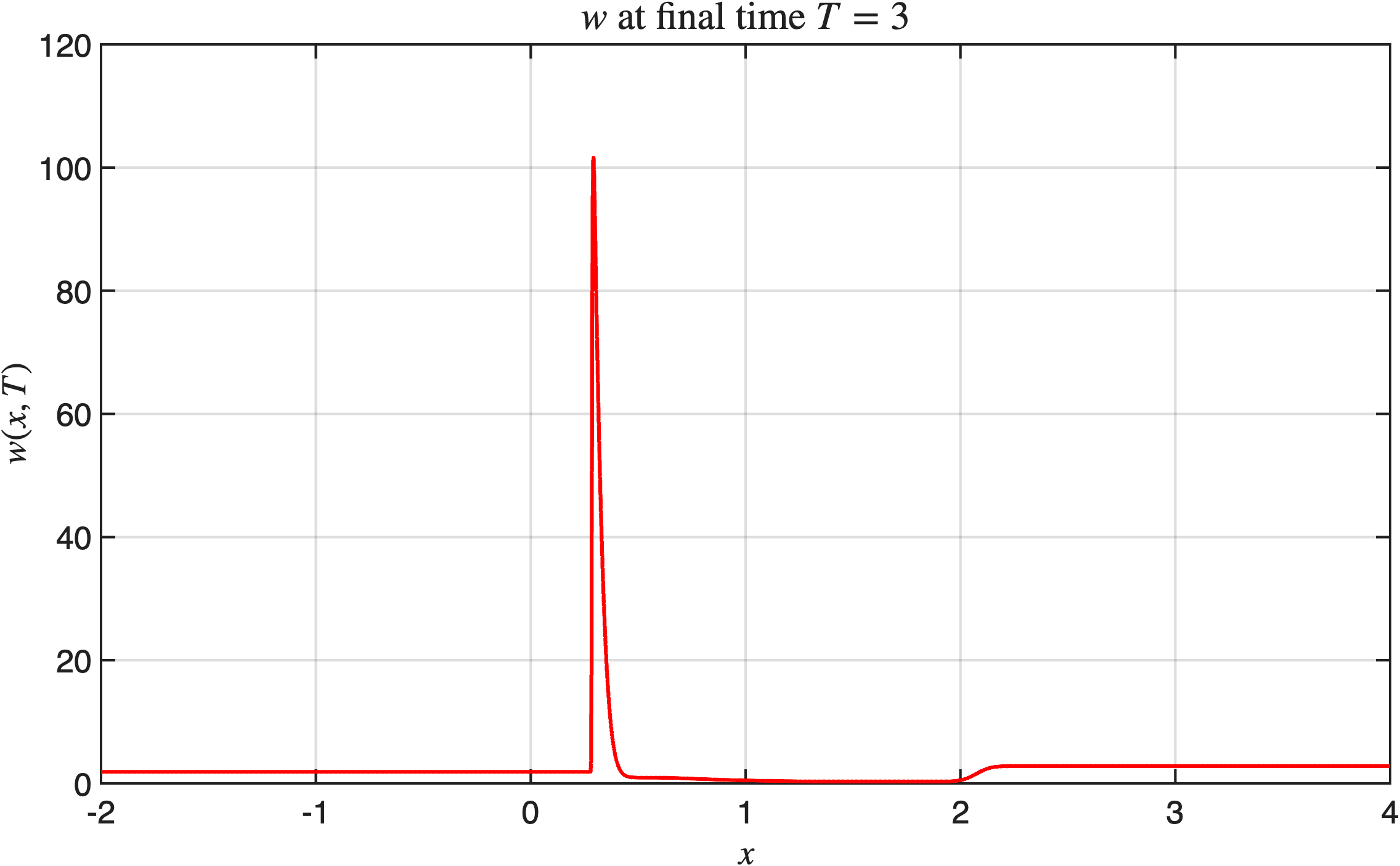} \\
(a) & (b) \\
\includegraphics[width=0.45\linewidth]{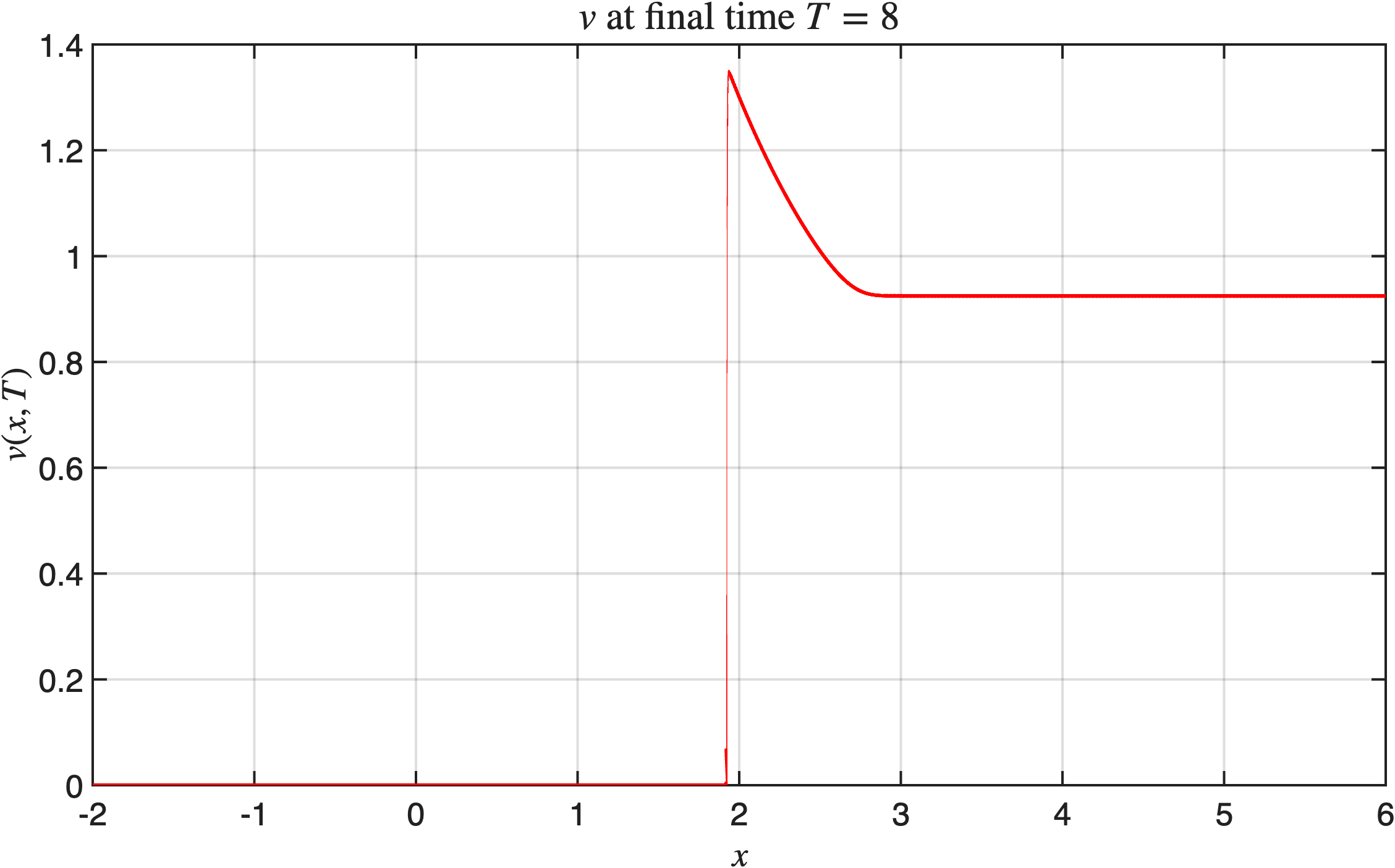} & 
\includegraphics[width=0.45\linewidth]{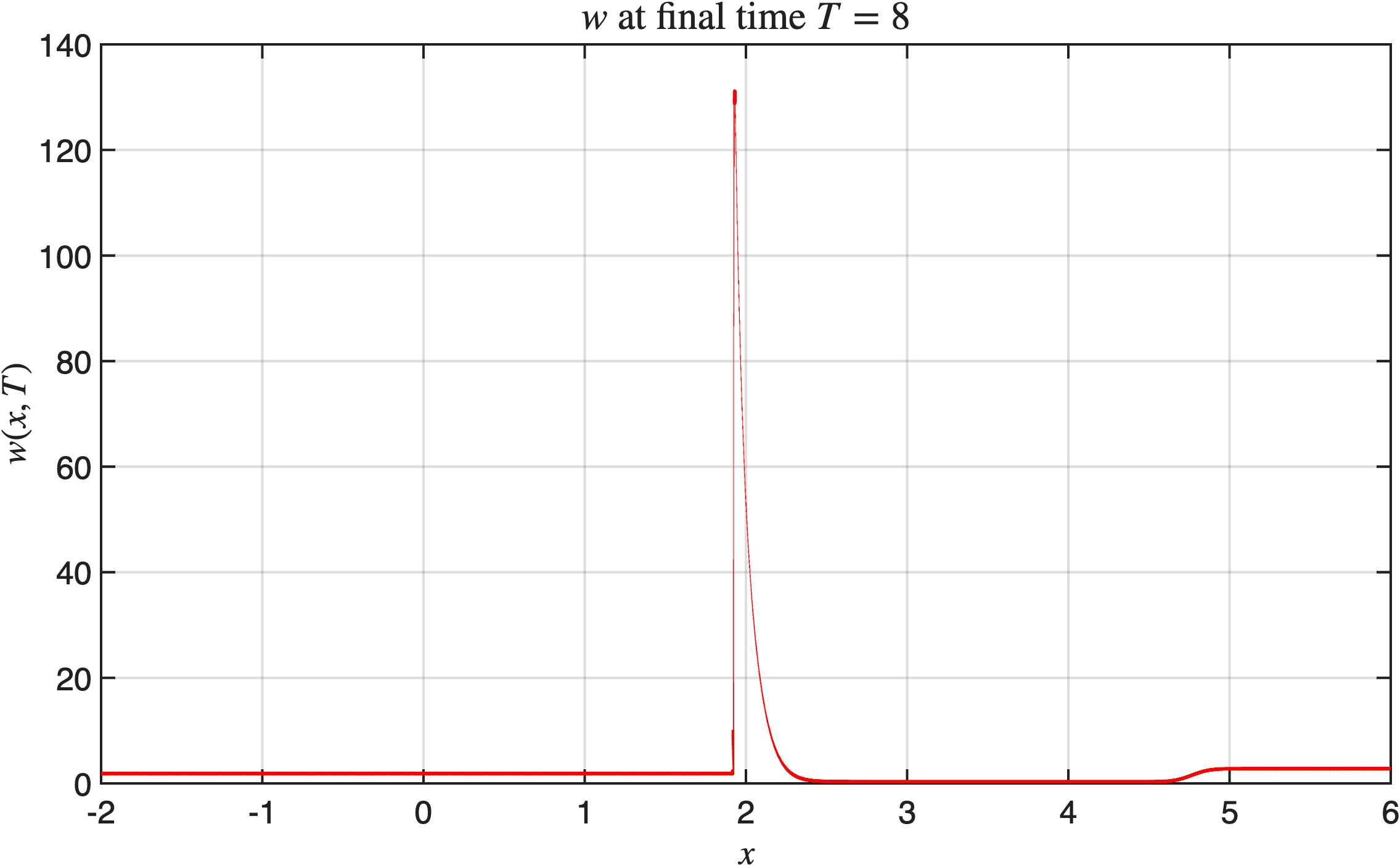} \\
(c) & (d) \\
\includegraphics[width=0.45\linewidth]{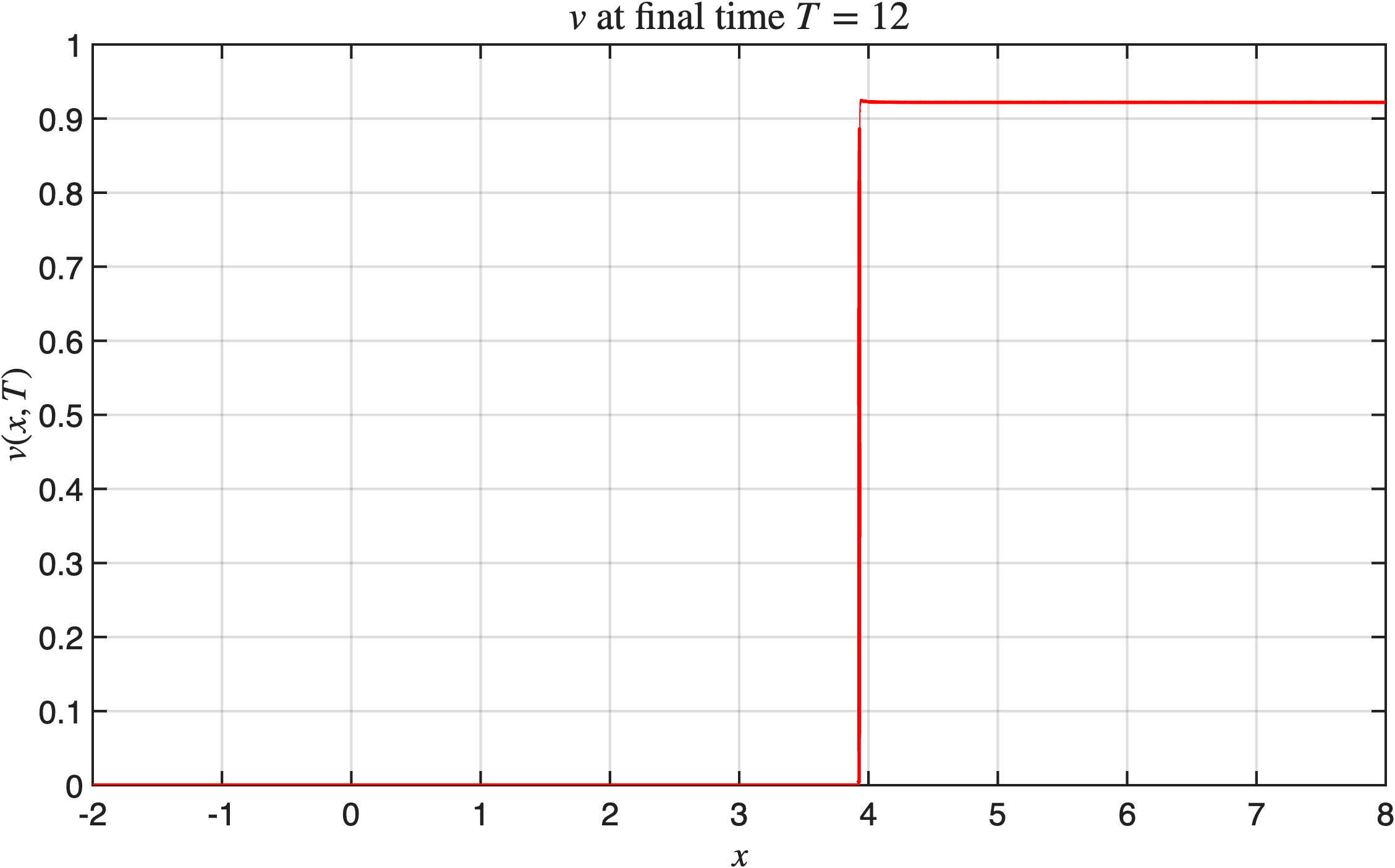} & 
\includegraphics[width=0.45\linewidth]{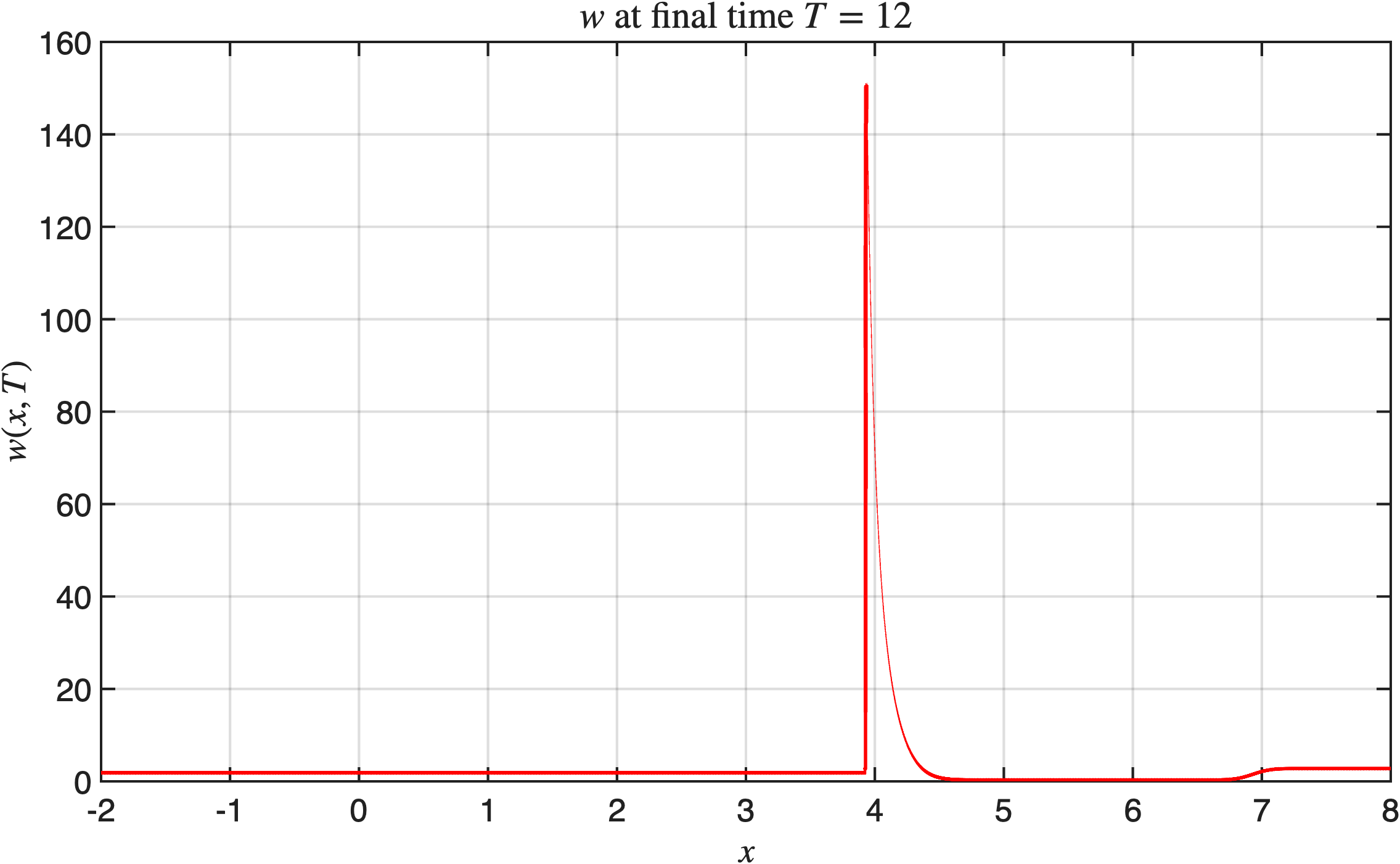} \\
(e) & (f)
\end{tabular}
    \caption{Profiles of $v(x,T)$ (left) and $w(x,T)$ (right) for Case~7 ($v_- = 0$, $v_+ < v_\thicksim$) at times $T=3$ (a)--(b), $T=8$ (c)--(d), and $T=12$ (e)--(f).}
    \label{fig:ProfileCase7}
\end{figure}

Figures~\ref{fig:ProfileCase7}(a)--(b) show the profiles of $v(x,T)$ and $w(x,T)$ at $T=3$, which corresponds to the time interval $0 \leq t < t_{*1} \approx 5.0391$. In (a), the left state $v_-h(3)=0$ is separated from the state $v_\thicksim h(3)=2.8084$ by the delta shock $\delta S$ at $x\approx 0.2865$, followed by the smooth transition of the rarefaction $R_1^{(\sim)}$ between $x\approx 0.6966$ and $x\approx 1.2865$, and the right state $v_+h(3)=0.9361$. 
In (b), the delta shock $\delta S$ appears as a very sharp spike at $x\approx 0.2865$ with strength $\alpha(3)\approx 4.4174$, while the small bump near $x\approx 2.0729$ corresponds to the contact discontinuity $J_2^{(*1)}$.
Figures~\ref{fig:ProfileCase7}(c)--(d) show the profiles at $T=8$, which lies in the interval $t_{*1} \leq t < t_{*2} \approx 11.2058$. After the first interaction, the delta shock $\delta S$ begins to cross the rarefaction $R_1^{(\sim)}$ as $\delta S_1$. In (c), the sharp jump at $x\approx 1.9121$ corresponds to $\delta S_1$ crossing $R_1^{(\sim)}$, with the residual rarefaction visible as the smooth transition from $V h(8)\approx 1.3479$ to $v_+h(8)=0.9248$ on the right. In (d), the spike of $\delta S_1$ at $x\approx 1.9121$ has grown to strength $\alpha_1(8)\approx 12.2328$.
Figures~\ref{fig:ProfileCase7}(e)--(f) show the profiles at $T=12$, which lies in the interval $t \geq t_{*2} \approx 11.2058$. After the second interaction, $\delta S_1$ has fully crossed $R_1^{(\sim)}$ producing the new delta shock $\delta S_2$. In (e), the smooth transition has completely disappeared and only one sharp jump remains at $x\approx 3.9306$, separating the left state $v_-h(12)=0$ from the right state $v_+h(12)=0.9220$. In (f), the spike of $\delta S_2$ at $x\approx 3.9306$ has grown further to strength $\alpha_2(12)\approx 16.7870$, while the small bump near $x\approx 6.9306$ corresponds to the contact discontinuity $J_2^{(*1)}$.

\begin{figure}[h]
    \centering
    \includegraphics[width=0.45\linewidth]{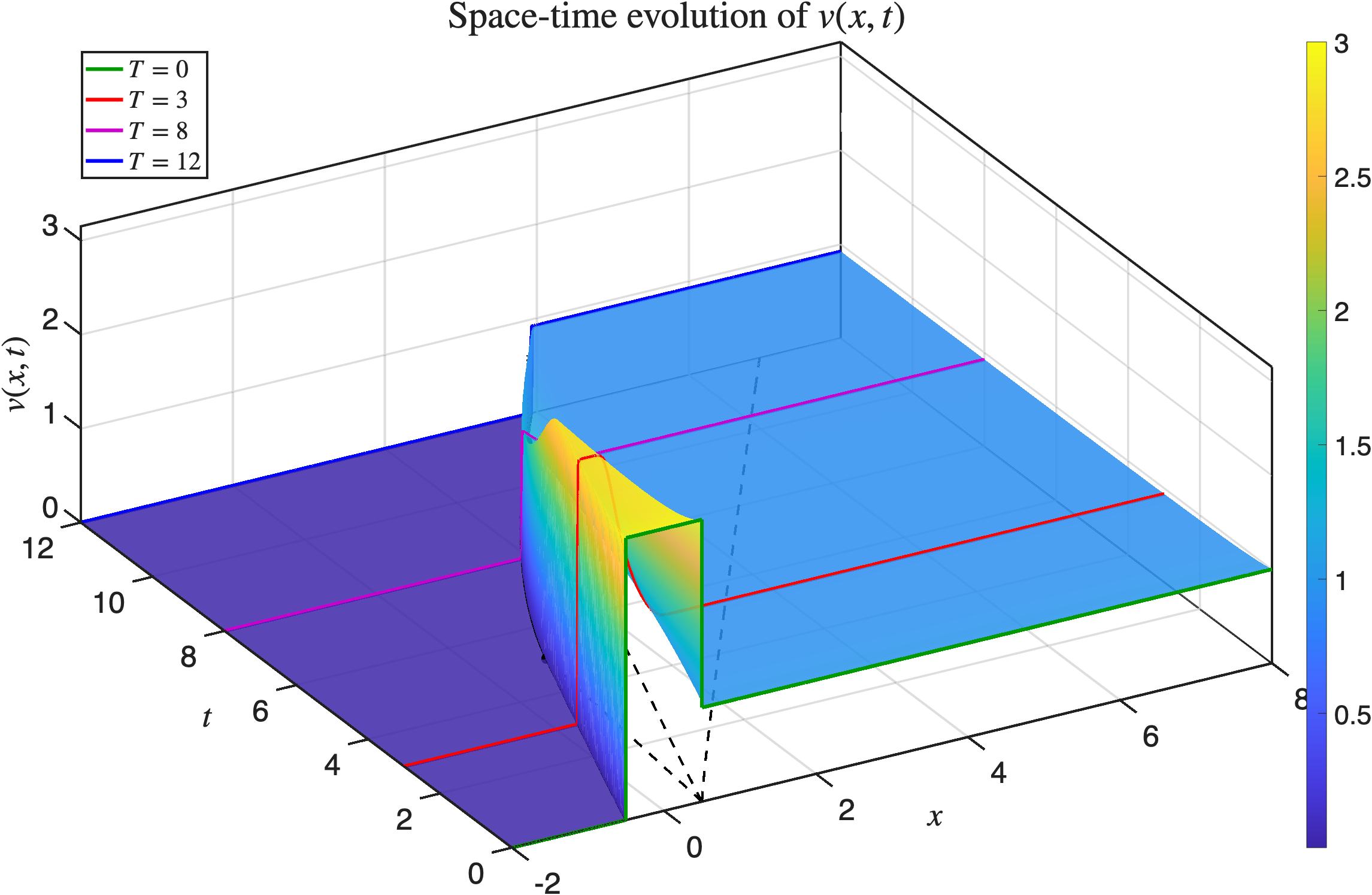}
    \includegraphics[width=0.45\linewidth]{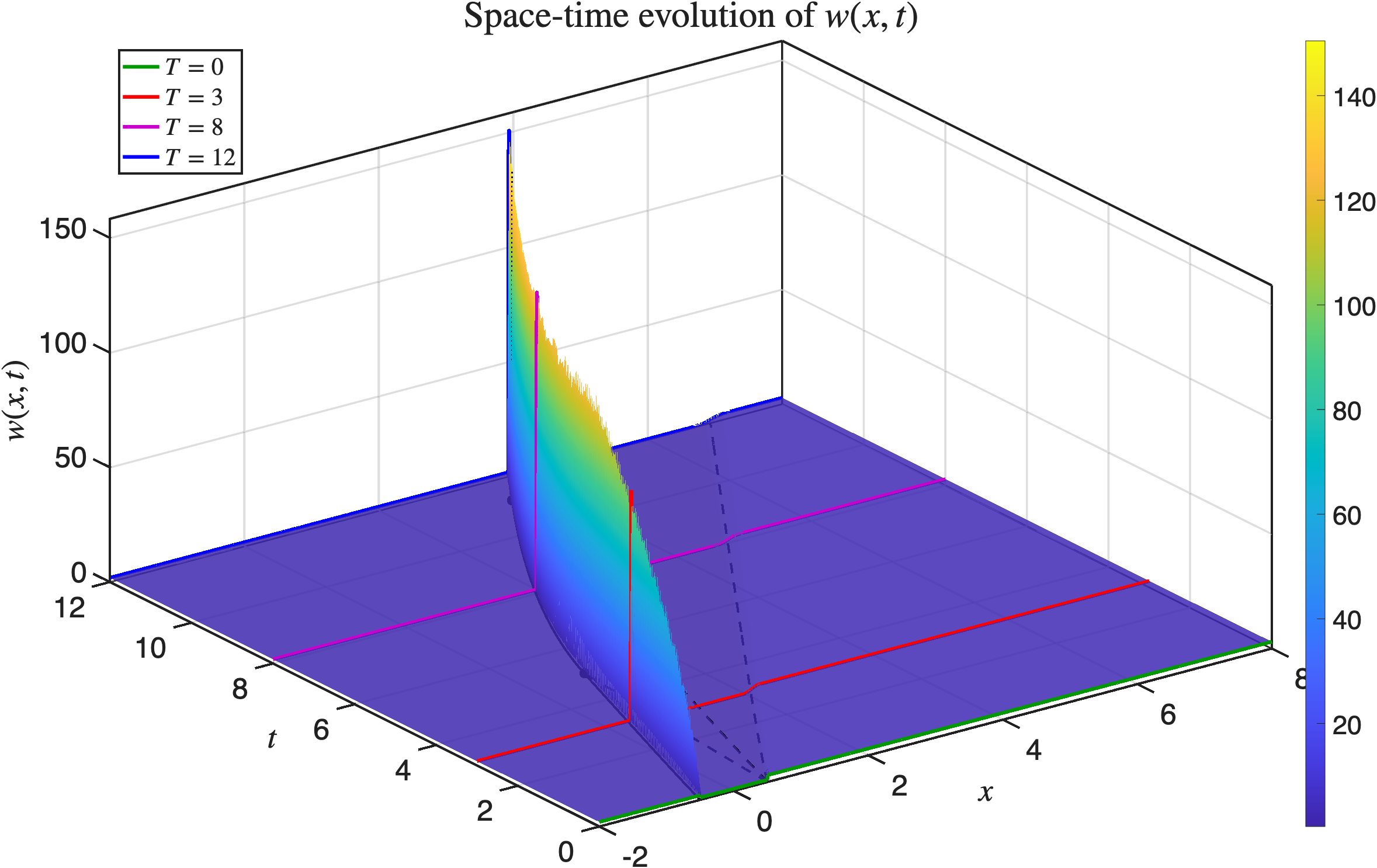}
    \caption{Space-time evolution of $v(x,t)$ (left) and $w(x,t)$ (right) for Case~7 ($v_- = 0$, $v_+ < v_\thicksim$), complementing Figure~\ref{fig:ProfileCase7}. 
    The solid colored lines show the profiles at $T=0$, $T=3$, $T=8$, and $T=12$, and the black curves on the $x$-$t$ plane show the delta shock trajectories (solid) and the rarefaction boundaries and contact discontinuity (dashed), with the interaction points marked as black dots.}
    \label{fig:3DCase7}
\end{figure}

Figure~\ref{fig:3DCase7} shows the space-time evolution of $v(x,t)$ and $w(x,t)$ for Case~7, complementing the profiles shown in Figure~\ref{fig:ProfileCase7}. The solid colored lines correspond to the profiles at times $T=0$, $T=3$, $T=8$, and $T=12$ already analyzed in Figure~\ref{fig:ProfileCase7}(a)--(f). 
On the $x$-$t$ plane, the solid black curve corresponds to the delta shock $\delta S$ (segment 1), the parabolic curve $x = \frac{1}{2} + \left(\sqrt{\int_0^t \exp\!\left(0.0880 \frac{s}{1+s}\right)ds} - \sqrt{3}\right)^2$ corresponds to $\delta S_1$ crossing $R_1^{(\sim)}$ (segment 2), and the straight line corresponds to $\delta S_2$ (segment 3), while the dashed curves correspond to the boundaries of $R_1^{(\sim)}$ and the contact discontinuity $J_2^{(*1)}$. The interaction points $(x_{*1},t_{*1})\approx(0.8333, 5.0391)$ and $(x_{*2},t_{*2})\approx(3.5000, 11.2058)$ are marked with black dots. In the right figure, the growing spike of $w(x,t)$ clearly illustrates the increasing strength of the delta shock over time.

\subsection{Numerical evidence of asymptotic stability}

In this section, we provide numerical evidence of the convergence of the perturbed Riemann solution to the Riemann solution as $\epsilon \to 0$, as established in Theorem~\ref{main_thm}. 
Specifically, we consider Case~1 ($0 < v_- < v_\thicksim < v_+$, with initial data $(v_-,w_-)=(1,2)$, $(v_\thicksim,w_\thicksim)=(2,3)$, $(v_+,w_+)=(4,1)$) and Case~7 ($v_-=0$, $v_+<v_\thicksim$, with initial data $(v_-,w_-)=(0,2)$, $(v_\thicksim,w_\thicksim)=(3,1)$, $(v_+,w_+)=(1,3)$). For each case, we fix the final time $T=12$ and vary the perturbation parameter 
\begin{equation*}
\epsilon \in \left\{0.50000,\, 0.20000,\, 0.10000,\, 0.01000,\, 0.00100,\, 0.00010,\, 0.00001\right\}.
\end{equation*}
For each value of $\epsilon$, we compute the exact wave curves of the perturbed Riemann solution and display them in the $x$-$t$ plane, together with the limiting Riemann solution (black thick line). As $\epsilon \to 0$, the wave curves of the perturbed solution converge to those of the Riemann solution. For Case~1, the limiting solution is $S_1 + J_2$, consisting of a 1-shock wave followed by a 2-contact discontinuity, while for Case~7, the limiting solution is a single delta shock wave $\delta S$.\\

{\bf Case 1: $0<v_-<v_\thicksim<v_+$.}

\begin{figure}[h]
    \centering
    \includegraphics[width=0.45\linewidth]{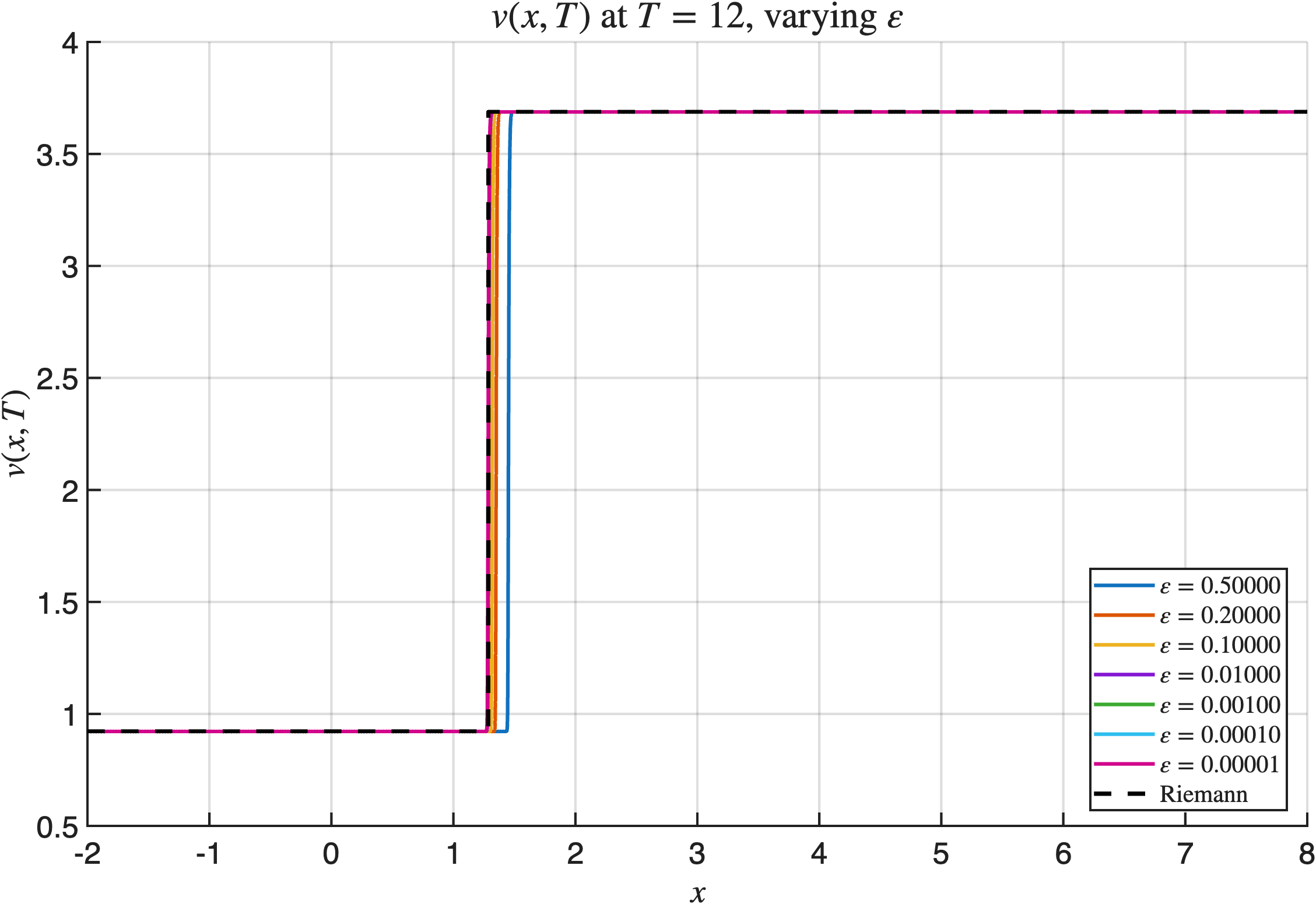}
    \includegraphics[width=0.45\linewidth]{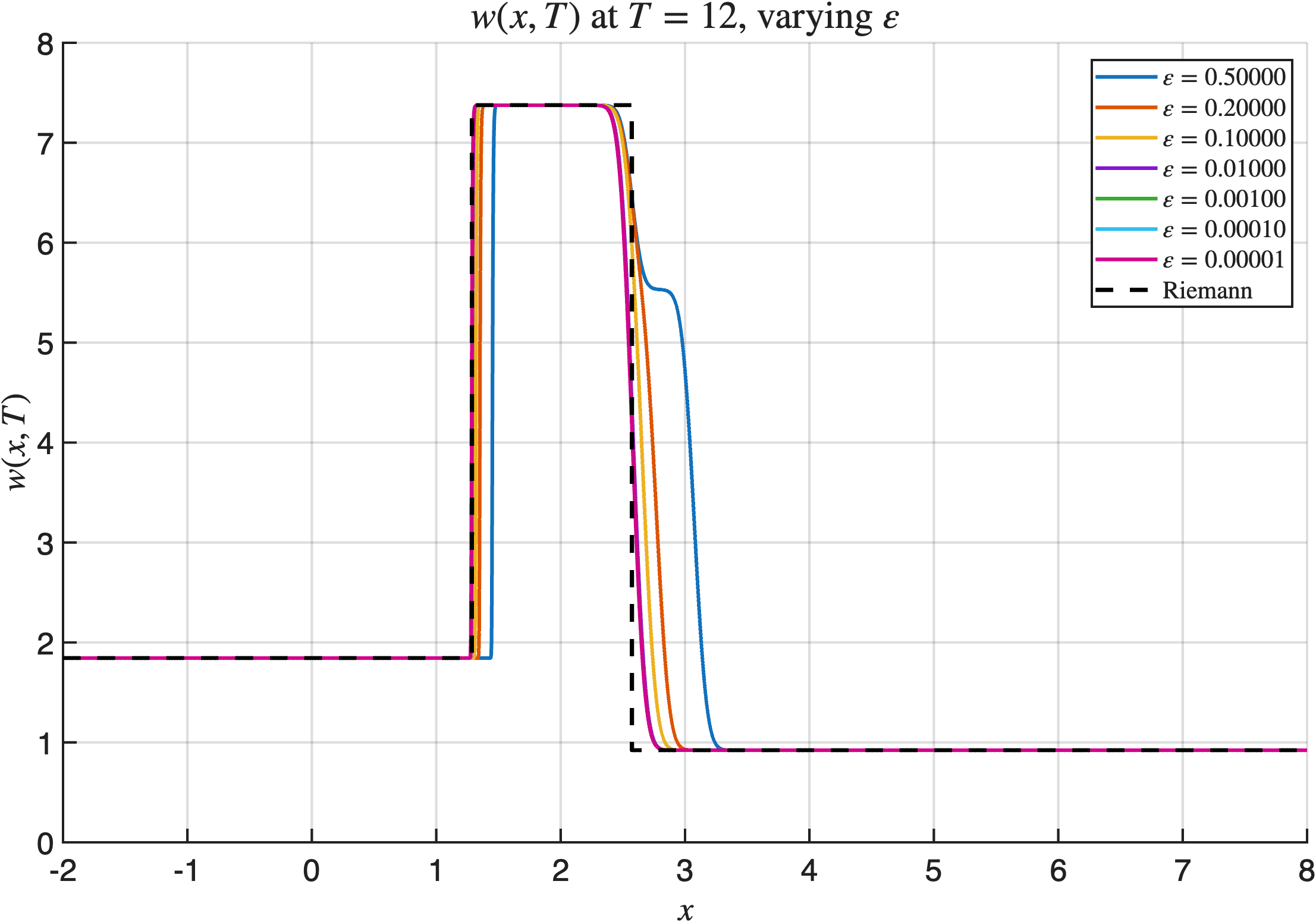}
    \caption{Profiles of $v(x,T)$ (left) and $w(x,T)$ (right) at $T=12$ for Case~1 ($0<v_-<v_\thicksim<v_+$), varying $\epsilon$. The black dashed line is the Riemann solution $S_1+J_2$.}
    \label{fig:AsymCase1}
\end{figure}

Figures~\ref{fig:AsymCase1} show the profiles of $v(x,T)$ and $w(x,T)$ at $T=12$ for varying $\epsilon$, together with the Riemann solution (black dashed line). In the left figure, the shock $S_1^{(*3)}$ shifts toward the Riemann shock $S_1$ as $\epsilon$ decreases, and for $\epsilon \leq 0.01000$ the profiles are virtually indistinguishable from the Riemann solution $S_1 + J_2$. In the right figure, the convergence is more visible: for $\epsilon = 0.50000$, the two parallel contact discontinuities $J_2^{(*3)}$ and $J_2^{(*2)}$ are clearly separated, producing a constant region at $W_{*3}h(12) \approx 7.3758$ between them, while as $\epsilon \to 0$ they merge into the single contact discontinuity $J_2$ of the Riemann solution, and this intermediate constant region shrinks to zero width.

\begin{figure}[h]
    \centering
    \includegraphics[width=0.6\linewidth]{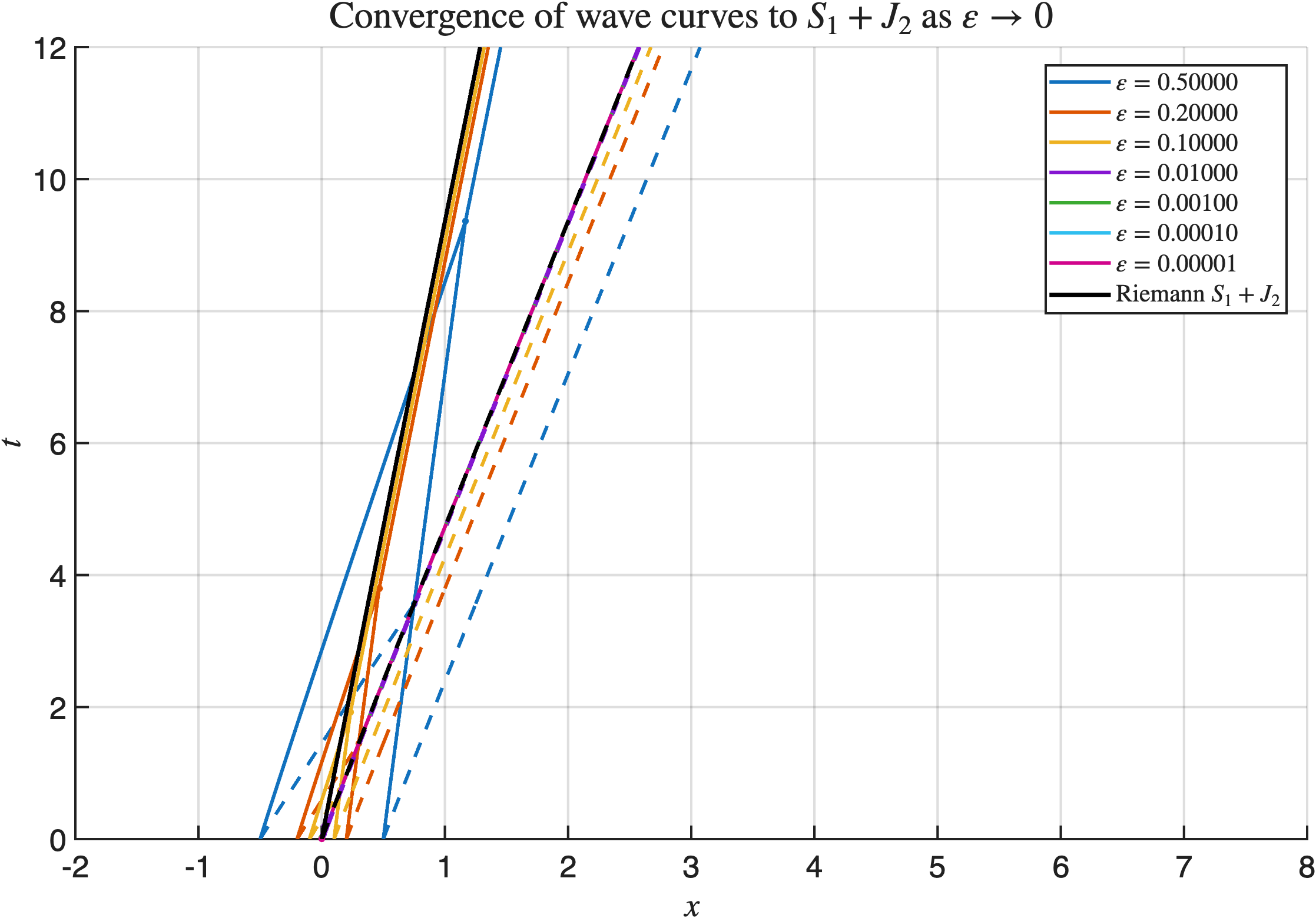}
    \caption{Convergence of the wave curves of the perturbed Riemann solution to $S_1 + J_2$ as $\epsilon \to 0$, for Case~1 ($0 < v_- < v_\thicksim < v_+$). Solid lines: shock curves. Dashed lines: contact discontinuities. Black thick lines: Riemann solution $S_1 + J_2$.}
    \label{fig:AsymCurCase1}
\end{figure}

Figure~\ref{fig:AsymCurCase1} shows the convergence of the wave curves of the perturbed Riemann solution to those of the Riemann solution $S_1 + J_2$ as $\epsilon \to 0$. The solid lines correspond to the shock curves $S_1^{(-)}$, $S_1^{(\sim)}$, $S_1^{(*1)}$ and $S_1^{(*3)}$, while the dashed lines correspond to the contact discontinuities $J_2^{(*1)}$, $J_2^{(*2)}$ and $J_2^{(*3)}$, for each value of $\epsilon$. The black thick solid and dashed lines are respectively the shock $S_1$ and the contact discontinuity $J_2$ of the Riemann solution. As $\epsilon$ decreases, the interaction points (visible as kinks in the solid curves) move toward the origin, and all wave curves converge to the Riemann solution curves. For $\epsilon \leq 0.01000$, the wave curves are virtually indistinguishable from those of the Riemann solution.\\

{\bf Case 7: $v_-=0$, $v_+<v_\thicksim$.}

\begin{figure}[h]
    \centering
    \includegraphics[width=0.45\linewidth]{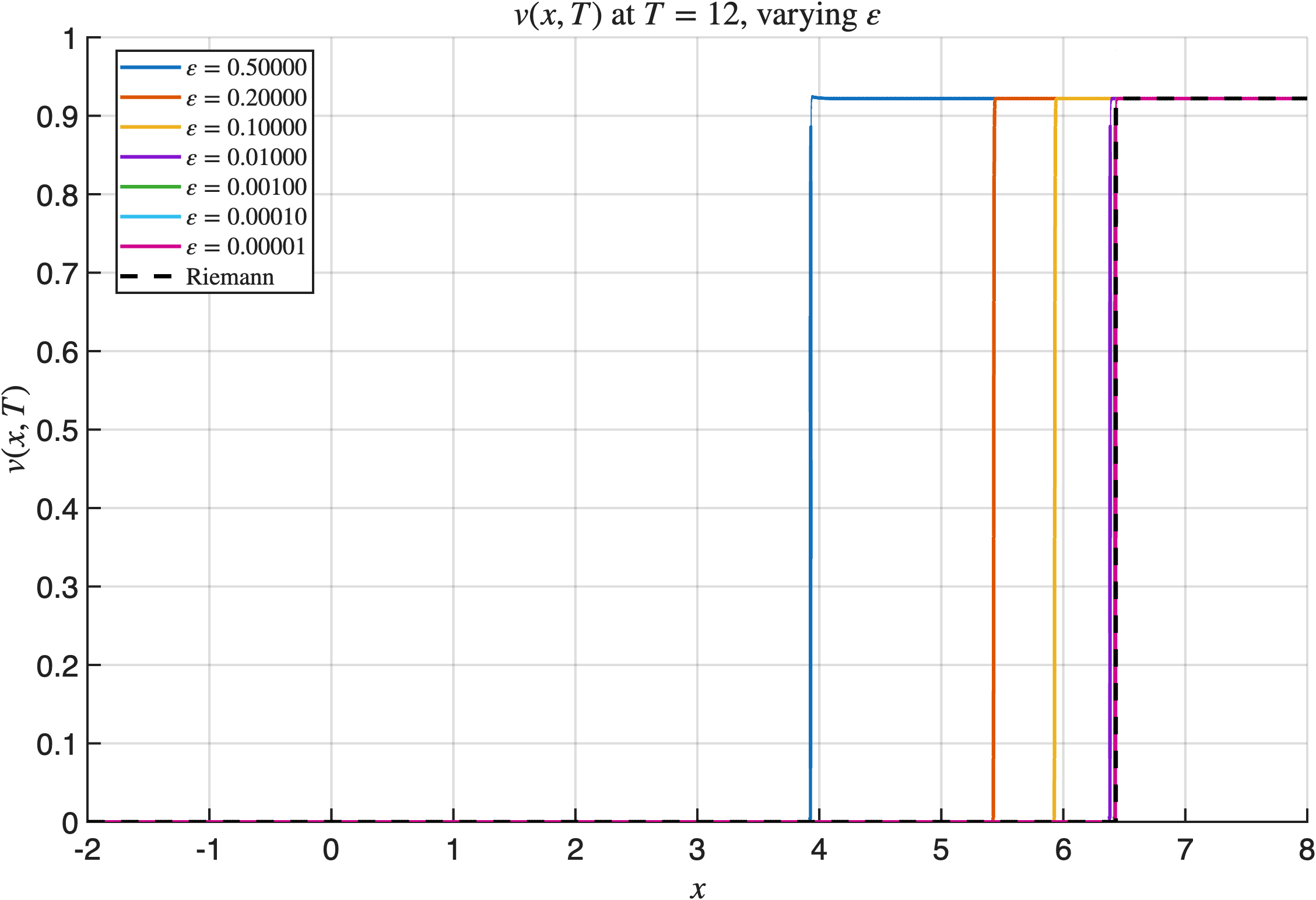}
    \includegraphics[width=0.45\linewidth]{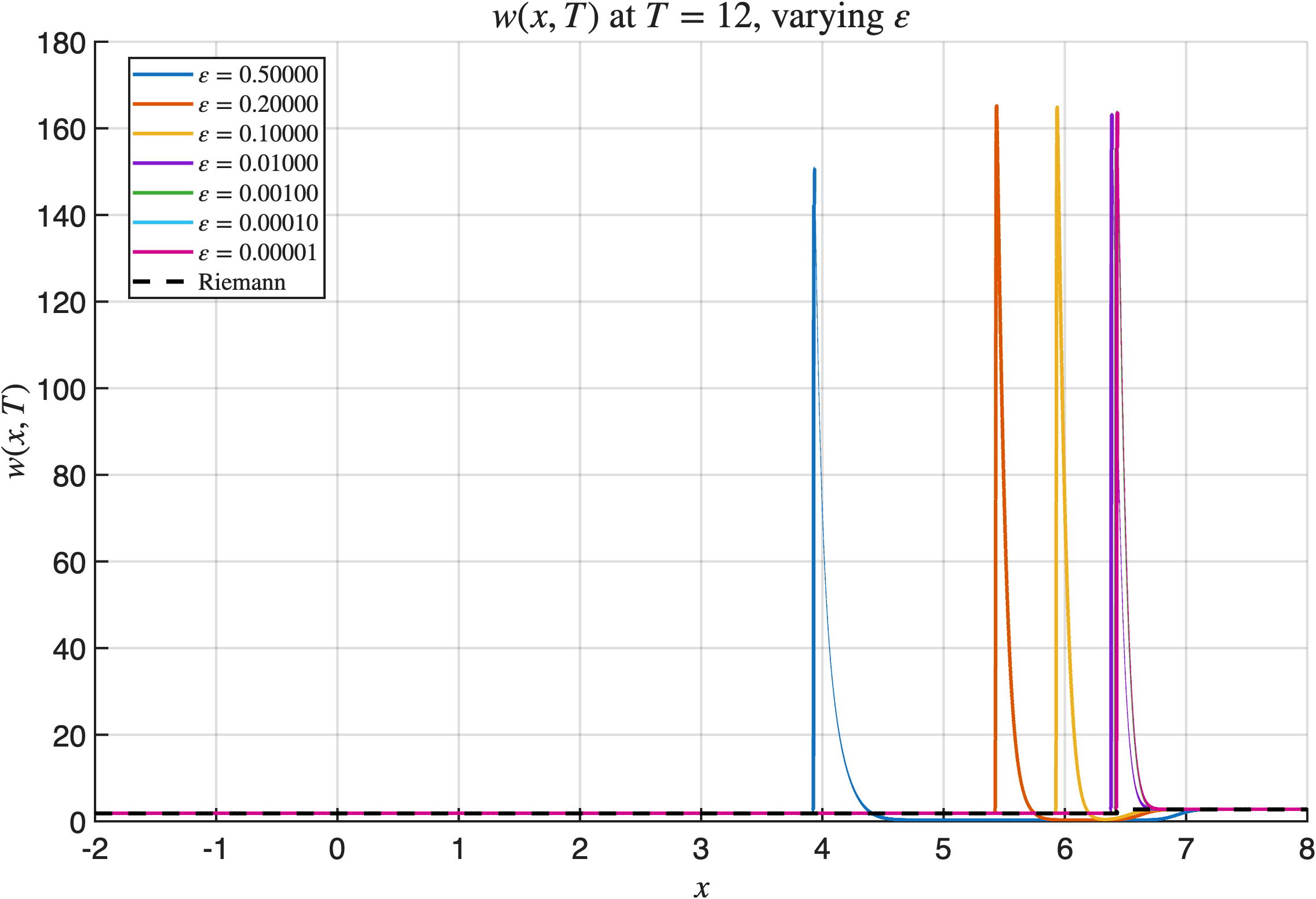}
    \caption{Profiles of $v(x,T)$ (left) and $w(x,T)$ (right) at $T=12$ for Case~7 ($v_-=0$, $v_+<v_\thicksim$), varying $\epsilon$. The black dashed line is the Riemann solution $\delta S$.}
    \label{fig:AsymCase7}
\end{figure}

Figures~\ref{fig:AsymCase7} show the profiles of $v(x,T)$ and $w(x,T)$ at $T=12$ for varying $\epsilon$, together with the Riemann solution (black dashed line). In the left figure, the position of the delta shock $\delta S_2$ moves to the right as $\epsilon$ decreases, converging to the Riemann delta shock position $x = \frac{1}{2}\Lambda(12) \approx 6.4306$. For $\epsilon = 0.50000$, the jump occurs at $x \approx 3.9306$, while for $\epsilon \leq 0.01000$ it is virtually at the Riemann position. In the right figure, the delta shock appears as a sharp spike whose position converges to the Riemann position as $\epsilon \to 0$. Moreover, as $\epsilon$ decreases, the spike becomes narrower and taller, concentrating toward a Dirac delta measure supported on the Riemann delta shock curve, which is the expected limiting behavior in the sense of distributions.

\begin{figure}[h]
    \centering
    \includegraphics[width=0.6\linewidth]{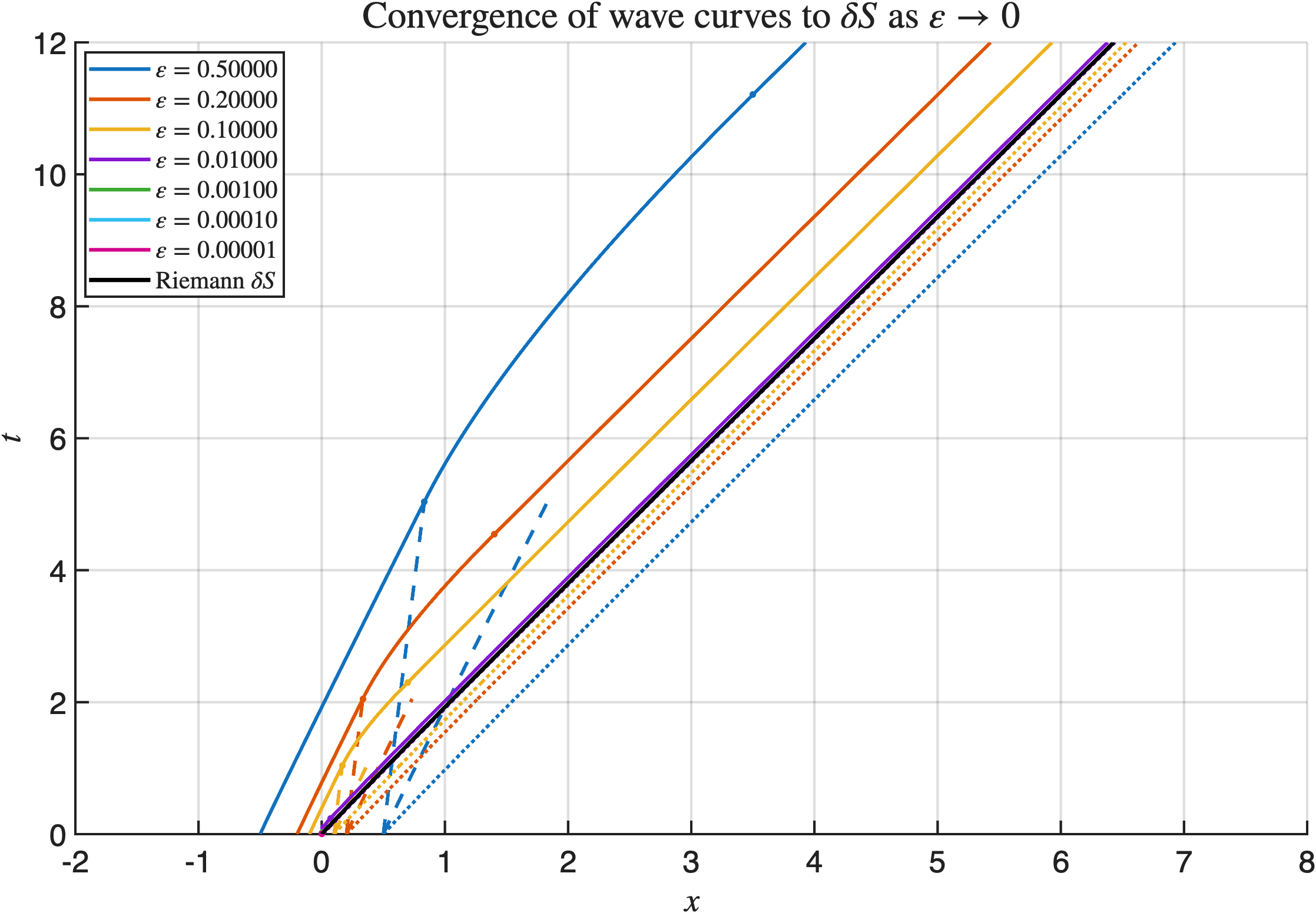}
    \caption{Convergence of the wave curves of the perturbed Riemann solution to $\delta S$ as $\epsilon \to 0$, for Case~7 ($v_- = 0$, $v_+ < v_\thicksim$). Solid lines: delta shock trajectories $\delta S$, $\delta S_1$ (parabolic) and $\delta S_2$. Dashed lines: boundaries of $R_1^{(\sim)}$. Dotted lines: contact discontinuity $J_2^{(*1)}$. Black thick line: Riemann solution $\delta S$.}
    \label{fig:AsymCurCase7}
\end{figure}

Figure~\ref{fig:AsymCurCase7} shows the convergence of the wave curves of the perturbed Riemann solution to the Riemann delta shock $\delta S$ as $\epsilon \to 0$. The solid lines correspond to the delta shock trajectories: the straight segment $\delta S$ for $0 \leq t \leq t_{*1}$, the curve $\delta S_1$ for $t_{*1} \leq t \leq t_{*2}$, and the straight segment $\delta S_2$ for $t \geq t_{*2}$. The dashed lines correspond to the boundaries of the rarefaction wave $R_1^{(\sim)}$, and the dotted lines to the contact discontinuity $J_2^{(*1)}$. The black thick solid line is the Riemann delta shock $\delta S$. The parabolic trajectory of $\delta S_1$ crossing $R_1^{(\sim)}$ is clearly visible for large $\epsilon$ (blue and orange curves), while the interaction points are marked as dots. As $\epsilon$ decreases, the wave curves converge to the Riemann delta shock, and for $\epsilon \leq 0.01000$ they are virtually indistinguishable from it.

\section{Conclusions and future direction} \label{sec:con}
In this paper, we investigated the Cauchy problem with piecewise constant perturbed initial data for a newly developed nonautonomous chromatography-type system of balance laws. We established the global weak solution for the Cauchy problem by analyzing all possible interactions of classical and nonclassical waves, and this became the first result in this direction for a nonautonomous system. We further analyzed the stability of the Riemann solution, including delta shock wave, by taking the vanishing limit of the perturbation parameter. Finally, we have shown numerical evidence for our theoretical results.

A potential next step is to establish the well-posedness of the general Cauchy problem for the nonautonomous system when the initial data is a Radon measure, where choosing the suitable entropy and energy conditions plays a major role. Another potential direction is to prove the existence of the Cauchy problem when the initial data is a bounded measurable function by possibly exploiting the compensated compactness framework for a suitable state space. Moreover, it would be interesting to investigate whether the two-components nonautonomous system could be extended to three-components or possibly $m$-components ($m\ge 4$) counterpart while maintaining the structural properties of the system.


\bibliographystyle{abbrv} 
\bibliography{References}

@article{GUO_2014,
title = {The perturbed {Riemann} problem and delta contact discontinuity in chromatography equations},
journal = {Nonlinear Anal.: TMA},
volume = {106},
pages = {110-123},
year = {2014},
author = {Lihui Guo and Lijun Pan and Gan Yin},
}

@ARTICLE{shen_JMAA_2010,
   author       = "C. Shen", 
   title        = "Wave interactions and stability of the {Riemann} solutions for the chromatography equations", 
   journal      = "J. Math. Anal. Appl.", 
   volume       = "365", 
   pages        = "609–618", 
   year         = "2010", 
}

@ARTICLE{Cheng_JMAA_2011,
   author       = "Hongjun Cheng and Hanchun Yang", 
   title        = "Delta shock waves in chromatography equations", 
   journal      = "J. Math. Anal. Appl.", 
   volume       = "380", 
   pages        = "475-485", 
   year         = "2011", 
}

@ARTICLE{Marco_Mazzotti_2009,
   author       = "Marco Mazzotti", 
   title        = "Nonclassical Composition Fronts in Nonlinear Chromatography: Delta-Shock", 
   journal      = "Ind. Eng. Chem. Res.", 
   volume       = "48", 
   pages        = "7733–7752", 
   year         = "2009", 
}

@ARTICLE{Marco_Mazzotti_2010,
   author       = "Marco Mazzotti and Abhijit Tarafder and Jeroen Cornel and Fabrice Gritti and Georges Guiochond", 
   title        = "Experimental evidence of a delta-shock in nonlinear chromatography", 
   journal      = "J. Chromatogr. A", 
   volume       = "1217", 
   pages        = "2002–2012", 
   year         = "2010", 
}

@ARTICLE{Ambrosio_SIAM,
   author       = "L. Ambrosio and G. Crippa and A. Figalli and L. Spinolo", 
   title        = "Some new well-posedness results for continuity and transport equations, and applications to the chromatography system", 
   journal      = "SIAM J. Math. Anal.", 
   volume       = "41", 
   pages        = "1890–1920", 
   year         = "2009", 
}

@ARTICLE{M_Sun_AML,
   author       = "M. Sun", 
   title        = "Interactions of delta shock waves for the chromatography equations", 
   journal      = "Appl. Math. Lett.", 
   volume       = "26", 
   pages        = "631–637", 
   year         = "2013", 
}

@ARTICLE{Q_Zhang_ZAMP,
   author       = "Q. Zhang", 
   title        = "Interactions of delta shock waves and stability of {Riemann} solutions for nonlinear chromatography equations", 
   journal      = "Z. Angew. Math. Phys.", 
   volume       = "67", 
   pages        = "15", 
   year         = "2016", 
}

@article{Langmuir_1916,
author = {Langmuir, Irving},
title = {THE CONSTITUTION AND FUNDAMENTAL PROPERTIES OF SOLIDS AND LIQUIDS. PART {I}. SOLIDS.},
journal = {J. Am. Chem. Soc.},
volume = {38},
pages = {2221-2295},
year = {1916},
}

@article{glueckauf_1946,
  title={Contributions to the theory of chromatography},
  author={Gl{\"u}eckauf, E},
  journal = {Proc. R. Soc. Lond. A},
  volume={186},
  pages={35--57},
  year={1946},
}

@article{glueckauf_1949,
  title={The General Theory of Two Solutes following Non-linear Isotherms},
  author={Gl{\"u}eckauf, E},
  journal = {Discuss. Faraday Soc.},
  volume={7},
  pages={12--25},
  year={1949},
}

@book{De_boer_1953,
author="{J. H. De Boer}",
title="The dynamical character of adsorption",
publisher = "Oxford, Clarendon Press",
year = "1953",
}

@article{Rhee_1970,
  title={On the theory of multicomponent chromatography},
  author={H.-K. Rhee and R. Aris and N. R. Amundson},
  journal ={Phil. Trans. R. Soc. A, Mathematical and Physical Sciences},
  volume={267},
  pages={419-455},
  year={1970},
}

@article{DELACRUZ_JDE,
title = {On a class of nonautonomous quasilinear systems with general time-gradually-degenerate damping},
journal = {J. Differential Equations},
volume = {416},
pages = {52-81},
year = {2025},
author = {Richard {De la cruz} and Wladimir Neves},
}

@article{Rcruz_RM_WN,
title = {Riemann problem for the {Chromatography-type} system of {Langmuir} isotherm with source term},
journal = {Nonlinear Analysis: RWA},
volume = {93},
pages = {104673},
year = {2027},
author = {Richard {De la cruz} and Rakib Mondal and Wladimir Neves}
}

@article{Liu&Smoller,
title = {On the vacuum state for isentropic gas dynamic equations},
journal = {Adv. Appl. Math.},
volume = {1},
pages = {345-359},
year = {1980},
author = {T P Liu and J Smoller}
}

@article{Shen&Sun_09,
title = {Interactions of delta shock waves for the transport equations with split delta functions},
journal = {J. Math. Anal. Appl.},
volume = {351},
pages = {747-755},
year = {2009},
author = {C Shen and M Sun}
}

@article{M_Nedeljkov_08,
title = {Interactions of delta shock waves in a strictly hyperbolic system of conservation laws},
journal = {J. Math. Anal. Appl.},
volume = {344},
pages = {1143-1157},
year = {2008},
author = {M Nedeljkov and M Oberguggenberger}
}

@article{C_Shen_10_NON,
title = {Stability of the Riemann solutions for a nonstrictly hyperbolic system of conservation laws},
journal = {Nonlinear Analysis, TMA},
volume = {73},
pages = {3284-3294},
year = {2010},
author = {C Shen and M Sun}
}

\end{document}